\input epsf
\documentstyle{amsppt}
\pagewidth{5.8truein}\hcorrection{0.0in}
\pageheight{8.7truein}\vcorrection{0.0in}
\TagsOnRight
\NoRunningHeads
\NoBlackBoxes
\catcode`\@=11
\def\logo@{}
\footline={\ifnum\pageno>1 \hfil\folio\hfil\else\hfil\fi}
\topmatter
\title Symmetry classes of spanning trees of Aztec diamonds and perfect matchings of
odd squares with a unit hole
\endtitle
\author Mihai Ciucu\endauthor
\thanks Research supported in part by NSF grant DMS-0500616.
\endthanks
\affil
  School of Mathematics, Georgia Institute of Technology\\
  Atlanta, Georgia 30332-0160
\endaffil
\define\Pe{\operatorname{P}}

\abstract 
We say that two graphs are similar if their adjacency matrices are similar matrices.
We show that the square grid $G_n$ of order $n$ is similar to the disjoint union of two copies of 
the quartered Aztec diamond $QAD_{n-1}$ of order $n-1$ with the path $P_n^{(2)}$ on $n$ vertices 
having edge weights equal to~2. Our proof is based on an explicit change of basis in the vector 
space on which the adjacency matrix acts. The arguments verifying that this change of basis works 
are combinatorial. It follows in particular that the characteristic polynomials of the above
graphs satisfy the equality $\Pe(G_n)=\Pe(P_n^{(2)})\,[\Pe(QAD_{n-1})]^2$. On the one hand, this provides a
combinatorial explanation for the ``squarishness'' of the characteristic polynomial of the square
grid---i.e., that it is a perfect square, up to a factor of relatively small degree. 
On the other hand, as formulas for the characteristic polynomials of the path and the square 
grid are well known, our equality determines the characteristic polynomial of the quartered Aztec 
diamond. In turn, the latter allows computing the number of spanning trees of quartered Aztec
diamonds.

We present and analyze three more families of graphs that share the above described ``linear squarishness'' property
of square grids: odd Aztec diamonds, mixed Aztec diamonds, and Aztec pillowcases---graphs obtained 
from two copies of an Aztec diamond by identifying the corresponding vertices on their convex hulls.

We apply the above results to enumerate all the 
symmetry classes of spanning trees of the even Aztec diamonds, and all the symmetry classes not involving rotations 
of the spanning trees of odd and mixed Aztec diamonds. We also enumerate all but the base case of the symmetry 
classes of perfect matchings of odd square grids with the central vertex removed. 
In addition, we obtain a product formula for the number of spanning trees of Aztec pillowcases.

\endabstract 
\endtopmatter
\document

\def\mysec#1{\bigskip\centerline{\bf #1}\message{ * }\nopagebreak\par\bigskip}

\def\myref#1{\item"{[{\bf #1}]}"} 
 
\def\pf{{\it Proof.\ }} 

\def\epf{\hfill{$\square$}\smallpagebreak}

\def\cite#1{\relaxnext@
  \def\nextiii@##1,##2\end@{[{\bf##1},\,##2]}%
  \in@,{#1}\ifin@\def\next{\nextiii@#1\end@}\else
  \def\next{[{\bf#1}]}\fi\next}
\def\proclaimheadfont@{\smc}

\def\pf{{\it Proof.\ }}
\define\Z{{\Bbb Z}}

\define\C{{\Bbb C}}
\define\M{\operatorname{M}}
\define\NE{\operatorname{NE}}
\define\NW{\operatorname{NW}}
\define\SW{\operatorname{SW}}
\define\SE{\operatorname{SE}}
\define\No{\operatorname{N}}
\define\So{\operatorname{S}}
\define\Ea{\operatorname{E}}
\define\We{\operatorname{W}}
\define\Te{\operatorname{T}}

\define\wt{\operatorname{wt}}

\define\te{\operatorname{t}}

\define\ms{{\!\!\!\!\!\!\!\!\!\!\!\!\!\!}}
\define\mms{{\!\!\!\!\!\!\!}}

\define\twoline#1#2{\line{\hfill{\smc #1}\hfill{\smc #2}\hfill}}
\define\threeline#1#2#3{\line{{\smc #1}\hfill{\smc #2}\hfill{\smc #3}}}
\define\threelinetext#1#2#3{\line{\hfill{\smc #1}\hfill{\smc #2}\hfill{\smc #3}\hfill}}

\def\mypic#1{\epsffile{figs/#1}}



\define\Biggs{1}
\define\Cone{2}
\define\Ctwo{3}
\define\CDGT{4}
\define\KPW{5}
\define\Kn{6}
\define\Kr{7}
\define\Lov{8}
\define\Stan{9}
\define\Temp{10}


%

\bigskip
\centerline{\bf Introduction}

\bigskip
The number of spanning trees of the Aztec diamond graph $AD_n$ (see Figure 2.1 for an illustration of $AD_5$) 
was shown by Knuth \cite{\Kn} to be given by a simple explicit product involving cosines. Stanley \cite{\Stan} 
posed then the problem of determining the number of spanning trees of the quartered Aztec diamonds (Figure 2.2 
shows the quartered Aztec diamond of order five). The author found a solution in the fall of 1996, when he was
a Postdoctoral Fellow at the Mathematical Sciences Research Institute in Berkeley (see also
\cite{\Ctwo}). This solution is presented in Section 2, and was the starting point of the current paper.
In the meanwhile, a different solution has been found by Richard Kenyon, James Propp and 
David Wilson \cite{\KPW, \S6.8}. 

A useful observation in Cvetcovi\'c et al.\ \cite{\CDGT} (reproduced here as Theorem 2.2) allows one to deduce the
number of spanning trees of certain planar graphs---the graphs considered in this paper included---from the characteristic 
polynomial of a slight modification of their dual. Our approach is to determine these characteristic polynomials
by showing that the involved graphs ``reduce'' to disjoint unions of graphs whose characteristic polynomials are known.
More precisely, given a graph whose characteristic polynomial we need to find, we provide a block diagonal matrix similar
to its adjacency matrix, so that the diagonal blocks are adjacency matrices of graphs with known characteristic polynomials.
Our solution is combinatorial, in that it accomplishes this by providing simple, conceptual changes of bases.

We apply this approach four times. First, we decompose this way the square grids, thus obtaining the characteristic
polynomial of quartered Aztec diamonds. Next, we decompose the odd and mixed Aztec diamonds, thereby obtaining the
characteristic polynomials of the odd and mixed ``halved'' Aztec diamonds (see Section 3). In Section 6 we handle
similarly the ``Aztec pillowcase'' graphs (we note that these are not related to the ``Aztec pillows'' previously
introduced in the literature by James Propp).

The results of Sections 2 and 3 allow us to solve most cases of two other natural enumeration questions: the symmetry
classes of spanning trees of Azted diamonds (see Section 4), and the symmetry classes of perfect matchings of odd by
odd square grids with the center vertex removed. We conclude the paper with a section posing some open problems.

\bigskip
\centerline{\bf 1. A similarity lemma for graphs possessing an automorphism of order two}
\centerline{\bf whose fixed points form a cut set}

\bigskip
The adjacency matrix of a weighted directed graph is the matrix whose rows and columns are indexed
by the vertices and whose $(u,v)$-entry equals the weight of the directed edge from $u$ to $v$ if
there is such an edge, and 0 otherwise. If the graph is undirected, one can replace each
edge $e$ by a pair of anti-parallel directed edges of the same weight as the weight of $e$, 
and use the previous definition. 
The characteristic polynomial of the graph $G$ is the characteristic polynomial of its adjacency matrix;
we denote it by $\Pe(G;x)$ (or simply $\Pe(G)$, if we need not display its argument). 
We think of unweighted graphs as weighted graphs in which all weights are equal to 1.

Let $G=(V,E)$ be a connected, undirected, weighted graph. 
For any subset $X\subset V$, denote by $\langle X\rangle$ the weighted subgraph of $G$ induced by $X$.



Let $T$ be an automorphism of $G$ so that $T^2$ is the identity. 
Denote by $V_0$ the set of vertices of $G$ fixed by $T$. 
Assume that $\langle V\setminus V_0\rangle$ is the union of two disjoint graphs of the form
$\langle V_1\rangle$ and $\langle T(V_1)\rangle$, for some suitable $V_1\subset V$ (this happens whenever
$V_0$ is a cut set). 

Set $G^+:=\langle V_1\rangle$. Let $G':=\langle V_1\cup V_0\rangle$.
Consider the directed graph on the vertex set $V_1\cup V_0$ that has a directed
edge from $u$ to $v$ if and only if $u$ and $v$ are adjacent in $G'$. Weight this edge by twice the
weight of the edge $\{u,v\}$ if $u\in V_0$ and $v\in V_1$, and by the weight of $\{u,v\}$ 
in all remaining instances. Denote the resulting weighted directed graph by $G^-$.

We say that the graphs $G_1$ and $G_2$ are {\it similar}---and write $G_1\sim G_2$---if
their adjacency matrices are similar.

\proclaim{Lemma 1.1} Under the above assumptions, $G$ is similar to the disjoint union of 
$G^+$ with $G^-$. 


\endproclaim


\pf Let $A$ be the adjacency matrix of $G$. For each vertex $v$ of $G$ consider an indeterminate~$e_v$.
Denote by $\bold U$ the complex vector space of formal linear combinations 
$\sum_{v\in V} c_v e_v$, $c_v\in\C$.

Let $F_G$ be the linear map from $\bold U$ to itself that sends $e_v$ to $\sum_{w\in V} a_{v,w} e_w$,
where $a_{v,w}$ is the weight of the edge $\{v,w\}$ if $v$ and $w$ are adjacent, and 0 otherwise. 
Then the matrix of $F_G$ in the
basis $B:=\{e_v:v\in V\}$ is just the adjacency matrix of $G$. 

Set $B':=\{(e_v-e_{T(v)})/2:v\in V_1\}$ and $B'':=\{(e_v+e_{T(v)})/2:v\in V_1\cup V_0\}$. One readily
sees that $B'\cup B''$ is a basis of $\bold U$.

For any $v\in V_1$ we have
$$
\align
F_G((e_v-e_{T(v)})/2)&=\frac12\left[\sum_{w\in V} a_{v,w} e_w-\sum_{w\in V}a_{T(v),w} e_w\right]
\\
&=\frac12\left[\sum_{w\in V} a_{v,w} e_w-\sum_{w\in V}a_{T(v),T(w)} e_{T(w)}\right]
\\
&=\frac12\sum_{w\in V} a_{v,w}(e_w-e_{T(w)})
\\
&=\frac12\sum_{w\in V_1} a_{v,w}(e_w-e_{T(w)}),
\endalign
$$
where at the last equality we used $V=V_1\cup V_0\cup T(V_1)$, the assumption that there are no
edges between $V_1$ and $T(V_1)$, and the fact that the summand in the next to last line above
vanishes when $w\in V_0$. Thus $B'$ spans an $F_G$-invariant subspace, and the matrix 
of its restriction to it is, in basis $B'$, just the adjacency matrix of $G^+$.

Similarly, for $v\in V_1$ we have
$$
\align
F_G((e_v+e_{T(v)})/2)&=\frac12\left[\sum_{w\in V} a_{v,w} e_w+\sum_{w\in V}a_{T(v),w} e_w\right]
\\
&=\frac12\left[\sum_{w\in V} a_{v,w} e_w+\sum_{w\in V}a_{T(v),T(w)} e_{T(w)}\right]
\\
&=\frac12\sum_{w\in V} a_{v,w}(e_w+e_{T(w)})
\\
&=\sum_{w\in V_1} a_{v,w}(e_w+e_{T(w)})/2 + \sum_{w\in V_0} a_{v,w}e_w,
\endalign
$$
while for $v\in V_0$
$$
\align
F_G((e_v+e_{T(v)})/2)&=F_G(e_v)=\sum_{w\in V} a_{v,w} e_w
\\
&=\sum_{w\in V_1} a_{v,w} e_w + \sum_{w\in T(V_1)} a_{v,w} e_w + \sum_{w\in V_0} a_{v,w} e_w 
\\
&=\sum_{w\in V_1} a_{v,w} e_w + \sum_{w\in V_1} a_{v,T(w)} e_{T(w)} + \sum_{w\in V_0} a_{v,w} e_w 
\\
&=\sum_{w\in V_1} a_{v,w} (e_w+ e_{T(w)}) + \sum_{w\in V_0} a_{v,w} e_w.
\endalign
$$ 
Therefore, $B''$ also spans an $F_G$-invariant subspace, and the matrix with respect to $B''$ of 
the restriction of $F_G$ to this $F_G$-invariant subspace is precisely the adjacency matrix of the 
weighted directed graph $G^-$.

Thus the matrix of $F_G$ in the basis $B'\cup B''$ is block-diagonal with the two blocks equal to the
adjacency matrices of $G^+$ and $G^-$, respectively, and the claim follows. \epf

\medskip
For the remainder of this paper, given a graph $G=(V,E)$ we will denote by ${\bold U}_G$ the complex 
vector space of formal linear combinations $\sum_{v\in V} c_v e_v$, $c_v\in\C$, the $e_v$'s
being indeterminates, and by $F_G$ the linear
map from ${\bold U}_G$ to itself that sends $e_v$ to $\sum_{w\in V} a_{v,w} e_w$, where $a_{v,w}$ 
is the $(v,w)$-entry of the adjacency matrix of $G$.

\mysec{2. The quartered Aztec diamond}

The Aztec diamond graph of order $n$, denoted $AD_n$, is the subgraph of the grid $(\Z+1/2)^2$
induced by the vertices $(x,y)$ with $|x|+|y|\leq n$ (Figure 2.1 shows $AD_5$). 
The {\it quartered Aztec diamond} $QAD_n$ is the subgraph of $AD_n$ induced by the vertices 
in its southeastern quarter ($QAD_5$ is pictured in Figure 2.2).

\topinsert
\twoline{\mypic{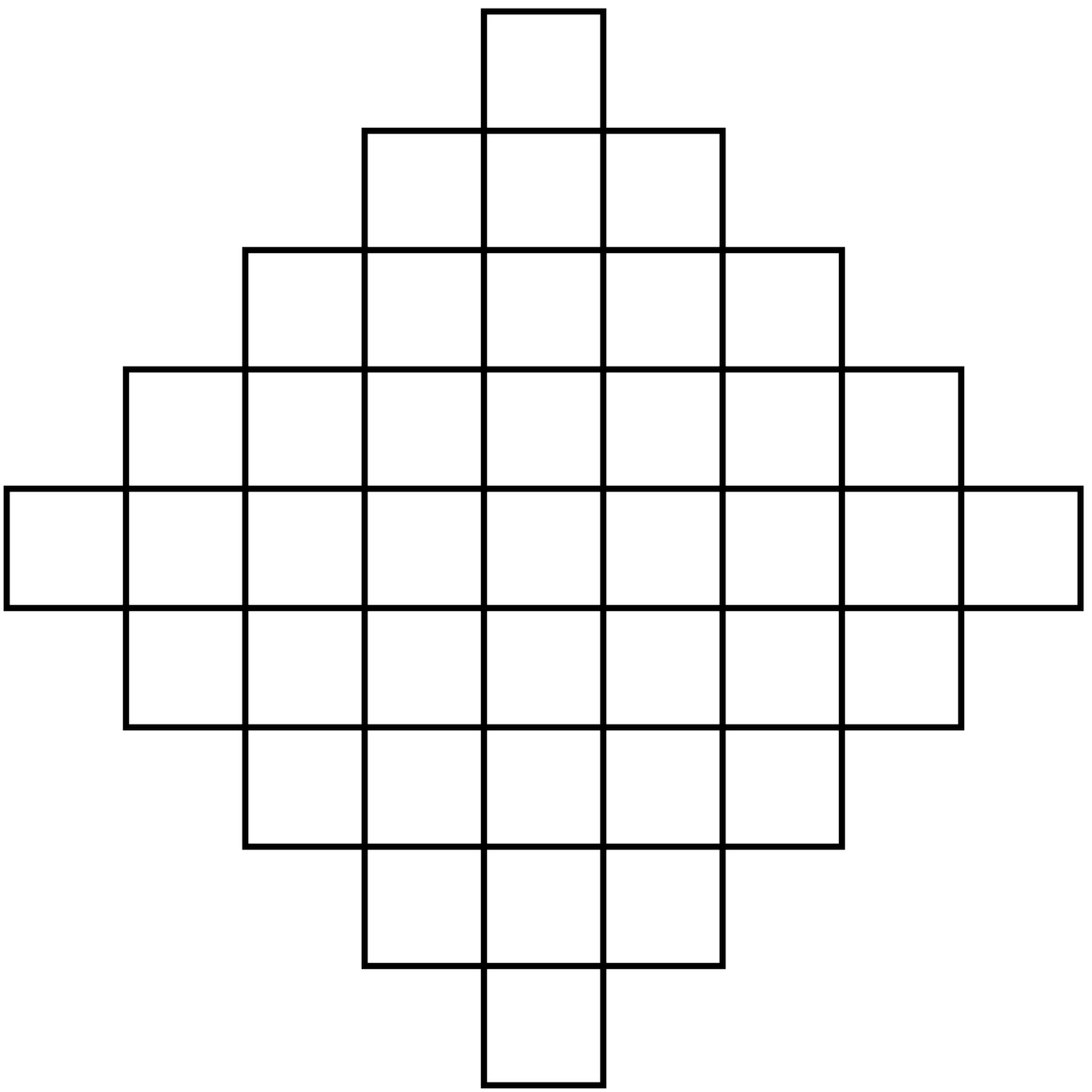}}{\mypic{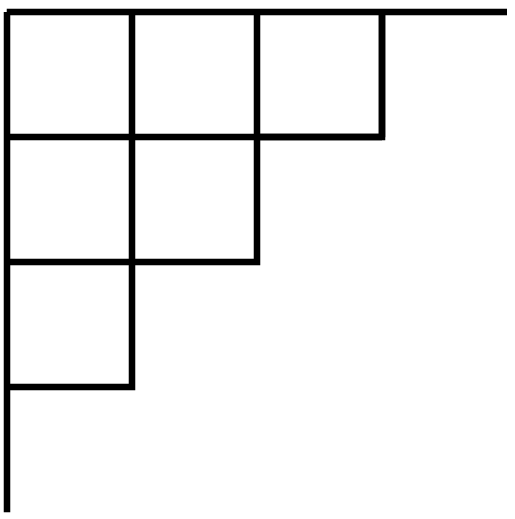}}
\twoline{Figure~2.1. {\rm  The Aztec diamond $AD_5$.}}
{Figure~2.2. {\rm The quartered Aztec diamond $QAD_5$.}}
\endinsert

\topinsert
\twoline{\mypic{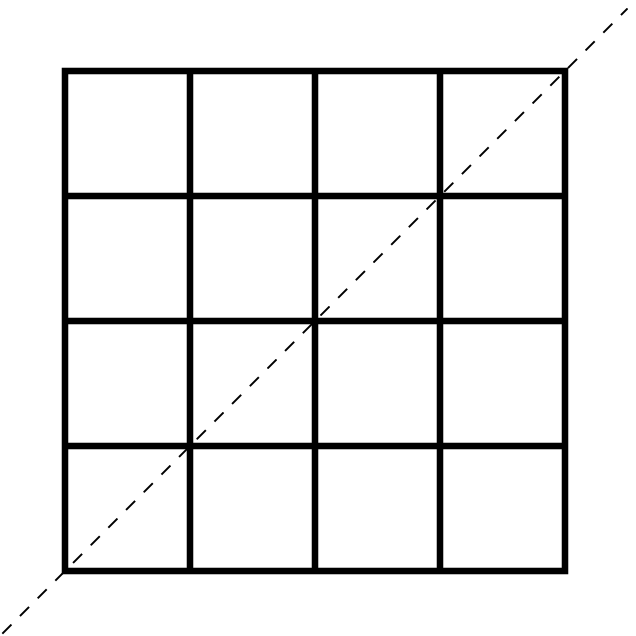}}{\mypic{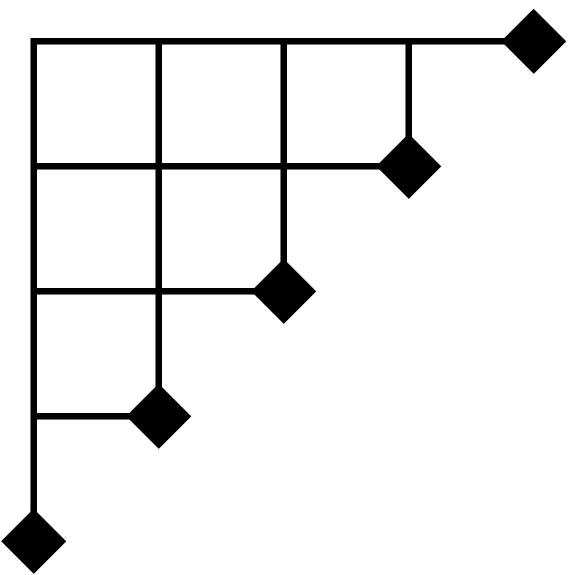}}
\twoline{Figure~2.3. {\rm  The square grid $G_5$.}}
{Figure~2.4. {\rm The weighted directed graph $\overline{QAD}_5$.}}
\endinsert

Denote by $G_n$ the $n\times n$ grid graph (illustrated in Figure 2.3 for $n=5$), and by $P_n^{(2)}$ 
the path on $n$ vertices having edge weights equal to 2.

\proclaim{Theorem 2.1} $G_n$ is similar to the disjoint union of two copies of $QAD_{n-1}$ with 
$P_n^{(2)}$:
$$
G_n\sim QAD_{n-1}\ \dot{\cup}\  QAD_{n-1}\ \dot{\cup}\ P_n^{(2)}.\tag2.1
$$


\endproclaim


The proof is presented after Corollary 2.3 below. 

Theorem 2.1 implies in particular that the characteristic polynomial of $G_n$ equals the product of
the characteristic polynomials of the graphs on the right hand side of (2.1). 

It is well known (see e.g.\ \cite{\Lov, problem 1.29}) 
that the eigenvalues of the path $P_n$ (with edge weights equal to 1) are
$$
\left\{2\cos\frac{\pi}{n+1},2\cos\frac{2\pi}{n+1},\dotsc,2\cos\frac{n\pi}{n+1}\right\};\tag2.2
$$
the characteristic polynomial of $P_n^{(2)}$ readily follows from this.

On the other hand, the square grid $G_n$ is the so-called tensor sum of $P_n$ with itself (see e.g. 
\cite{\Kr}), and thus by \cite{\Kr, Theorem I} its eigenvalues are given by all possible
sums of two of the numbers in (2.2):
$$
 2\cos\frac{j\pi}{n+1}+2\cos\frac{k\pi}{n+1},\ \ \ 1\leq j,k\leq n.\tag2.3
$$
It follows that
$$
\Pe(QAD_{n-1};x)=\prod_{1\leq i<j\leq n}\left(x-2\cos\frac{j\pi}{n+1}-2\cos\frac{k\pi}{n+1}\right).
\tag2.4
$$

The following result of \cite{\CDGT} provides a useful connection between characteristic polynomials and 
the number of spanning trees of certain planar graphs. For completeness, we include the proof. 
For a graph $G$ we denote by $\te(G)$ the number of its spanning trees.

\proclaim{Theorem 2.2 \cite{\CDGT}}
Let $G$ be a connected planar graph all of whose bounded faces have $r$ sides. 
Let $G'$ be the graph obtained
from the planar dual of $G$ by deleting the vertex corresponding to the infinite face. Then
$$
\te(G)=\Pe(G';r).
$$
\endproclaim

\pf By a well known fact, $G$ has the same number of spanning trees as its 
planar dual $G^\perp$. By the Matrix Tree Theorem (see e.g., \cite{\Biggs}), 
$\te(G^\perp)$ is equal to the determinant of 
the matrix obtained from the negative Laplacian of $G^\perp$ (i.e., the difference 
between the diagonal matrix of its vertex degrees and its adjacency matrix) by deleting the row and 
column indexed by the vertex of $G^\perp$ corresponding to the infinite face of $G$. 
The hypothesis implies that the latter matrix is precisely the evaluation of the 
characteristic polynomial of $G'$ at $r$. \epf

\proclaim{Corollary 2.3} The number of spanning trees of $QAD_{n}$ equals
$$
\prod_{1\leq j<k\leq n-1}\left(4-2\cos\frac{j\pi}{n}-2\cos\frac{k\pi}{n}\right).\tag2.5
$$

\endproclaim

\pf This is a direct consequence of Equation (2.4), Theorem 2.2, and the fact that the graph obtained 
from the planar dual of $QAD_{n}$ by deleting the vertex corresponding to the infinite face is
isomorphic to $QAD_{n-2}$. \epf

\medskip
{\it Proof of Theorem 2.1.} Apply Lemma 1.1 to the square grid $G_n$, choosing $T$ to be the symmetry
across one of its diagonals (see Figure 2.3). The resulting graph $G^+$ is isomorphic
to $QAD_{n-1}$. The resulting weighted directed graph $G^-$ is obtained from $QAD_{n}$ by
replacing all its edges by pairs of anti-parallel directed edges of weight 1, and then changing to 2
the weights of all directed edges originating from the $n$ vertices on the hypotenuse of the convex
hull of $QAD_{n}$; denote the latter by $\overline{QAD}_{n}$ (Figure 2.4 shows $\overline{QAD}_{5}$;
edges between non-marked vertices mean that there are directed edges of weight 1 between both 
corresponding ordered pairs of vertices; an edge between an unmarked vertex $v$ and a marked vertex $w$
signifies a directed edge from $v$ to $w$ of weight 1, and a directed edge from $w$ to $v$ of weight 2).
We obtain
$$
G_n\sim QAD_{n-1}\ \dot{\cup}\ \ \overline{QAD}_{n}.\tag2.6
$$
Let ${\bold U}={\bold U}_{\overline{QAD}_{n}}$ and $F=F_{\overline{QAD}_{n}}$ be as defined at the
end of Section 1.

As proved in Lemma 2.6 below, the union of the set $B$ of $n(n-1)/2$ vectors in the statement 
of Lemma 2.4 with the set $B'$ of $n$ vectors in the statement of Lemma 2.5 forms a basis of 
${\bold U}$. 
By Equation (2.9), the subspace spanned by $B$ is $F$-invariant. Lemma 2.5 proves that the subspace spanned 
by $B'$ is $F$-invariant as well. Thus the matrix of the linear map $F$ in the basis $B\cup B'$ is a
block diagonal matrix with two blocks. By (2.9), the first of these is precisely
the adjacency matrix of $QAD_{n-1}$. By (2.13), the second block equals the adjacency matrix 
of the path $P_n^{(2)}$. Thus
$$
\overline{QAD}_{n}\sim QAD_{n-1}\ \dot{\cup}\  P_n^{(2)}.\tag2.7
$$
This and (2.6) imply (2.1). \epf

\medskip
In order to present the lemmas employed in the above proof, we
coordinatize the vertices of $\overline{QAD}_{n}$ using matrix-style coordinates, assigning the top
left vertex coordinates $(1,1)$. Let $e_{ij}$ be the indeterminate corresponding to vertex
$(i,j)$ in ${\bold U}={\bold U}_{\overline{QAD}_{n}}$. Thus $\{e_{ij}:i,j\geq1,i+j\leq n+1\}$
is a basis of 
${\bold U}$, and the matrix of $F=F_{\overline{QAD}_{n}}$ in this basis is the adjacency matrix of
$\overline{QAD}_{n}$.

\proclaim{Lemma 2.4} For $i,j\geq 1$, $i+j\leq n$, define
$$
v_{ij}:=e_{i-1,j}-e_{i,j-1}+e_{i+1,j}-e_{i,j+1}\tag2.8
$$
$($if any index is outside the range $[1,n]$, we omit the corresponding term$)$.



Then for all $i,j\geq 1$, $i+j\leq n$, we have
$$
F(v_{ij})=v_{i-1,j}+v_{i,j-1}+v_{i+1,j}+v_{i,j+1},\tag2.9
$$
where by convention undefined $v_{kl}$'s on the right hand side $($these occur when $k$ and $l$ fail 
to satisfy both $k,l\geq 1$ and $k+l\leq n$$)$ are omitted.



\endproclaim

\topinsert
\twoline{\mypic{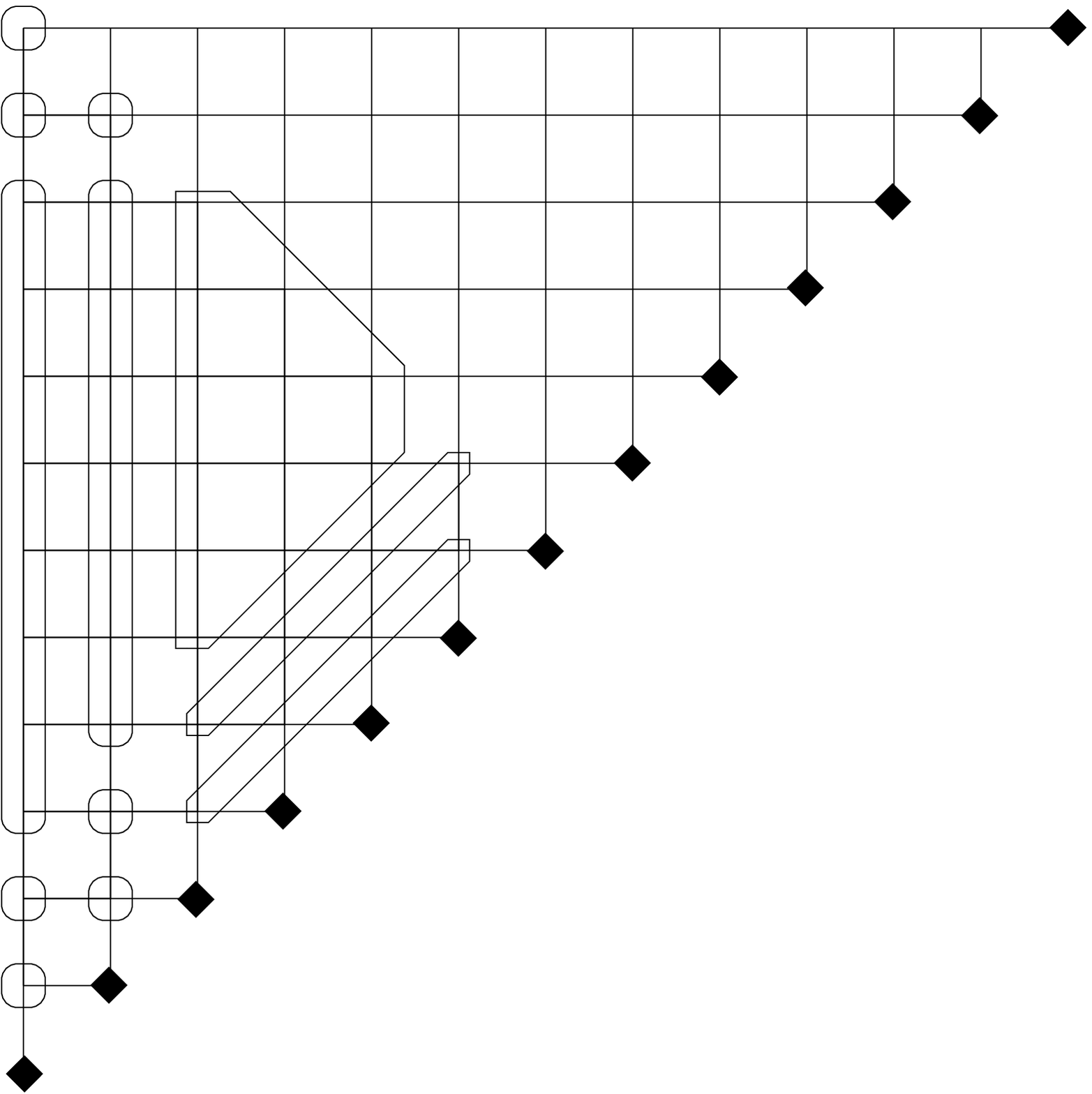}}{\mypic{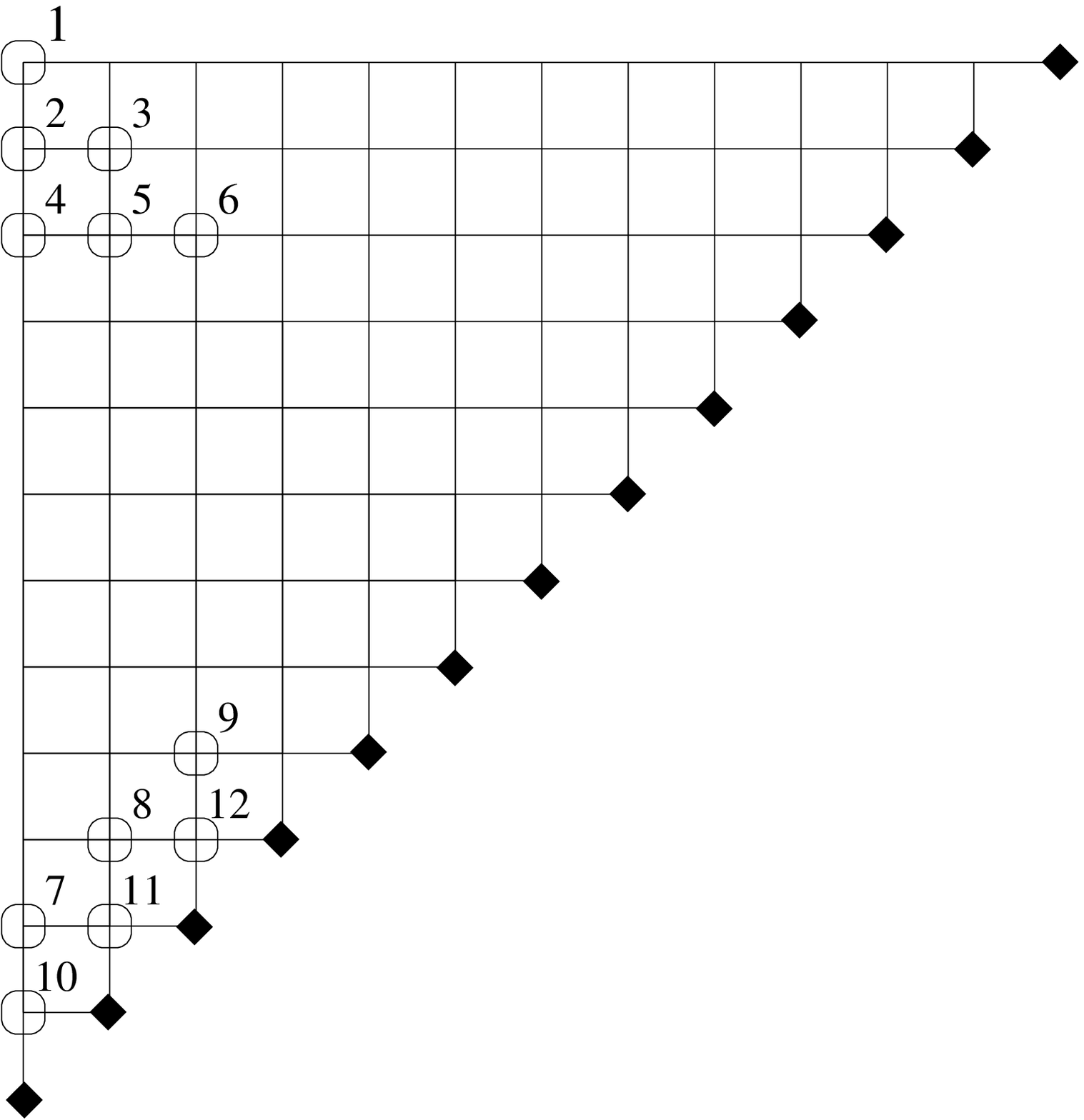}}
\twoline{Figure~2.5{\rm (a). The twelve cases.}}
{Figure~2.5{\rm (b). Twelve representatives.}}
\endinsert

\topinsert
\centerline{\mypic{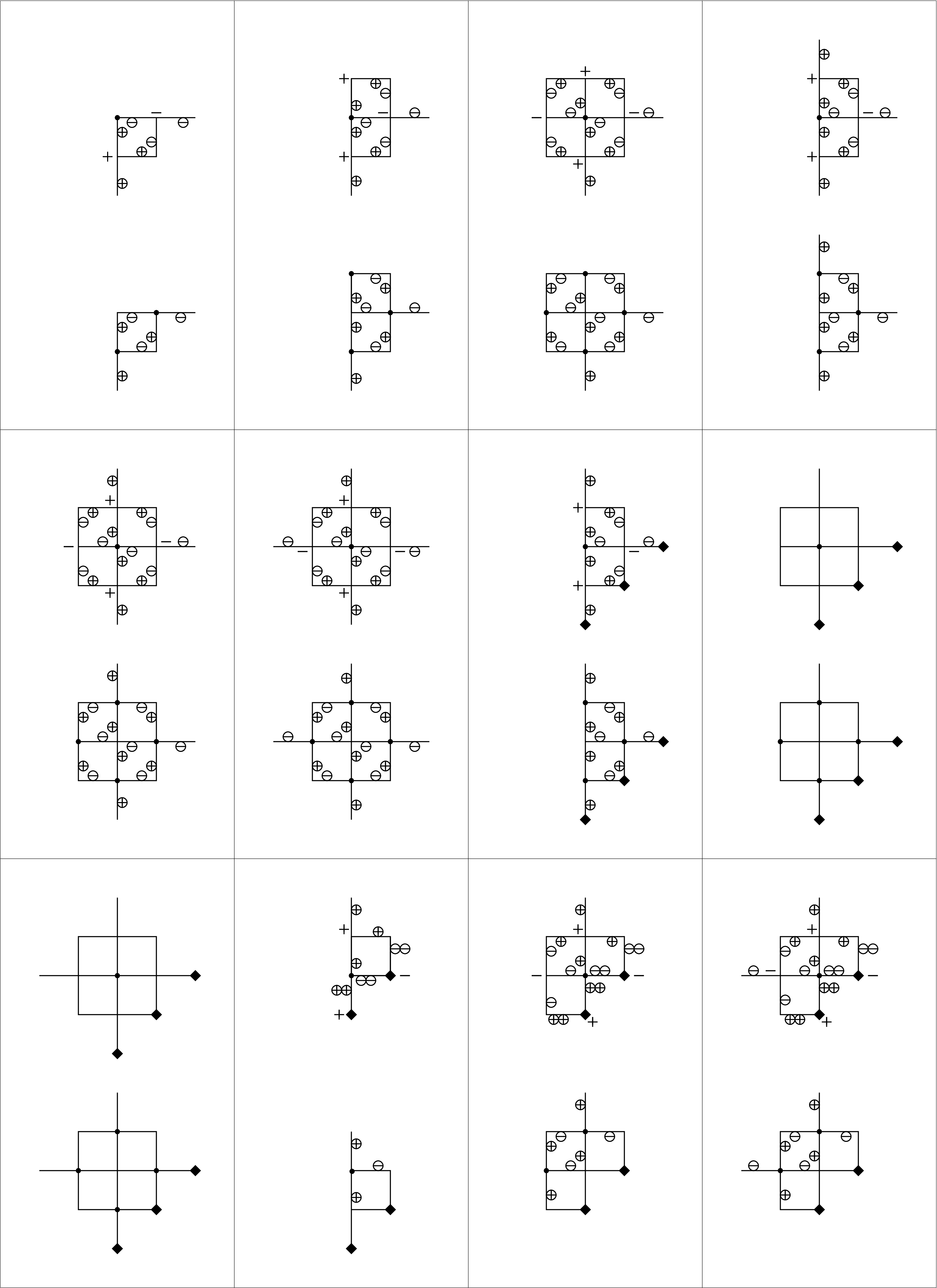}}
\centerline{{\smc Figure~2.6.} {\rm Checking the twelve cases.}}
\endinsert

\pf By the definition of $F$, we have
$$
F(e_{ij})=e_{i-1,j}+e_{i,j-1}+e_{i+1,j}+e_{i,j+1},\tag2.10
$$
for all $i,j\geq 1$, $i+j\leq n$, where we omit the terms on the right hand side corresponding to
undefined $e_{k,l}$'s, i.e. to those index values that fail to satisfy both $k,l\geq1$ and
$k+l\leq n+1$. 

Call the union of the neighborhoods of all neighbors of a vertex in a graph the {\it second neighborhood} 
of that vertex.
Then by equations (2.8), (2.10), and the linearity of $F$, both sides of (2.9) are linear combinations of
$e_{kl}$'s corresponding to vertices $(k,l)$ from the second neighborhood of vertex~$(i,j)$. 

It is apparent from the definition of $\overline{QAD}_n$ that, up to reflection in its symmetry axis,
there are only twelve translationally
distinct second neighborhoods of its vertices; the induced partition of the set of vertices is
shown in Figure 2.5(a) (pairs of vertices of $\overline{QAD}_n$ that are mirror images with respect
to its symmetry axis are understood to be in the same class of the partition). Because of this,
it is enough to check that (2.9) holds for one representative of each of the
twelve classes---for instance those indicated in Figure 2.5(b).

These twelve checkings are accomplished by Figure 2.6; the order of the panels corresponds to the
numbering of the representative vertices shown in Figure 2.5(b). 

In each panel there are two copies of the second neighborhood of the vertex $(i,j)$ in question. 
The vertex $(i,j)$ is indicated by the black dot. 
The top copy performs the computation of the left hand side of (2.9), using (2.8), (2.10), and the 
linearity of $F$. 
Consider for instance the top half of the sixth panel. By (2.8) and the linearity of $F$,
$$
F(v_{ij})=F(e_{i-1,j})-F(e_{i,j-1})+F(e_{i+1,j})-F(e_{i,j+1}).
$$
Each term $\pm F(e_{k,l})$ on the right hand side above is depicted in the figure by placing the
corresponding uncircled plus or minus sign next to the vertex $(k,l)$.
By (2.10), we have
$$
\align
F(e_{i-1,j})-F(e_{i,j-1})+F(e_{i+1,j})-F(e_{i,j+1})=&(e_{i-2,j}+e_{i-1,j-1}+e_{i,j}+e_{i-1,j+1})
\\
-&(e_{i-1,j-1}+e_{i,j-2}+e_{i+1,j-1}+e_{i,j})
\\
+&(e_{i,j}+e_{i+1,j-1}+e_{i+2,j}+e_{i+1,j+1})
\\
-&(e_{i-1,j+1}+e_{i,j}+e_{i+1,j+1}+e_{i,j+2}).
\endalign
$$
The sixteen $\pm e_{k,l}$'s on the right hand side are marked on the figure by placing circled plus 
or minus signs next to vertex $(k,l)$. Then the coefficient of any $e_{k,l}$ in $F(v_{ij})$ is
visually apparent from the top half of the sixth panel: Simply ``add up'' the circled signs, counting
a circled plus as 1, and a circled minus as $-1$ (the coefficient is clearly zero unless vertex $(k,l)$ is shown 
in the figure).

To continue with our example, turn now to the bottom copy of the sixth panel in Figure 2.6. It performs the computation of the right
hand side of (2.9), in the following sense. By (2.8), we have
$$
\align
v_{i-1,j}+v_{i,j-1}+v_{i+1,j}+v_{i,j+1}=&(e_{i-2,j}-e_{i-1,j-1}+e_{i,j}-e_{i-1,j+1})
\\
+&(e_{i-1,j-1}-e_{i,j-2}+e_{i+1,j-1}-e_{i,j})
\\
+&(e_{i,j}-e_{i+1,j-1}+e_{i+2,j}-e_{i+1,j+1})
\\
+&(e_{i-1,j+1}-e_{i,j}+e_{i+1,j+1}-e_{i,j+2}).
\endalign
$$
Each $\pm e_{k,l}$'s above is depicted in the figure by placing a circled plus 
or minus sign next to vertex $(k,l)$. This affords a visual computation of the result, by the same
adding of circled signs as above.

One checks by inspection that in the top half of the sixth panel in Figure 2.6, each vertex gets the same sum of circled signs
as the corresponding vertex in the bottom half. This establishes Equation (2.9) for vertices of type 6.

The remaining cases are checked similarly by the rest of Figure 2.6. (The calculations
required in panels 8 and 9 are identical to those in panels 5 and 6, respectively; thus we did not repeat 
them.) The double
circled plusses and minuses in the last three panels indicate contributions of $\pm2 e_{kl}$, and appear
due to the fact that the directed edge from a marked to an unmarked vertex of $\overline{QAD}_n$ has
weight 2 (marked vertices are indicated by black diamonds).
Note that in the bottom panels there are no contributions
coming from the diamond vertices, because by the convention stated after Equation (2.9), the undefined
terms on the right hand side of the latter are omitted (and by the statement after Equation (2.8), $v_{ij}$ is not
defined if $(i,j)$ is a diamond vertex). \epf

\proclaim{Lemma 2.5} Let
$$
c_{ij}:=\left\{\matrix 2, &&i,j\geq1, i+j\leq n\\1, &&i,j\geq1, i+j=n+1,\endmatrix\right.\tag2.11
$$
and define
$$
w_k:=\sum_{ {{\scriptstyle n-k+2\leq i+j\leq n+1}
\atop {\scriptstyle i-j\in\{n-k,n-k-2,n-k-4,\dotsc,k-n\}} }}
c_{ij}e_{ij},\tag2.12
$$
for $k=1,\dotsc,n$ $($for $n=9$ these are represented in Figures $2.7(i)$--$(ix)$$)$. Then 
$$
F(w_k)=2w_{k-1}+2w_{k+1}\tag2.13
$$
for all $k=1,\dotsc,n$, where by convention on the right hand side we omit terms of index
outside the range $[1,n]$.

\endproclaim

\topinsert
\threeline{\mypic{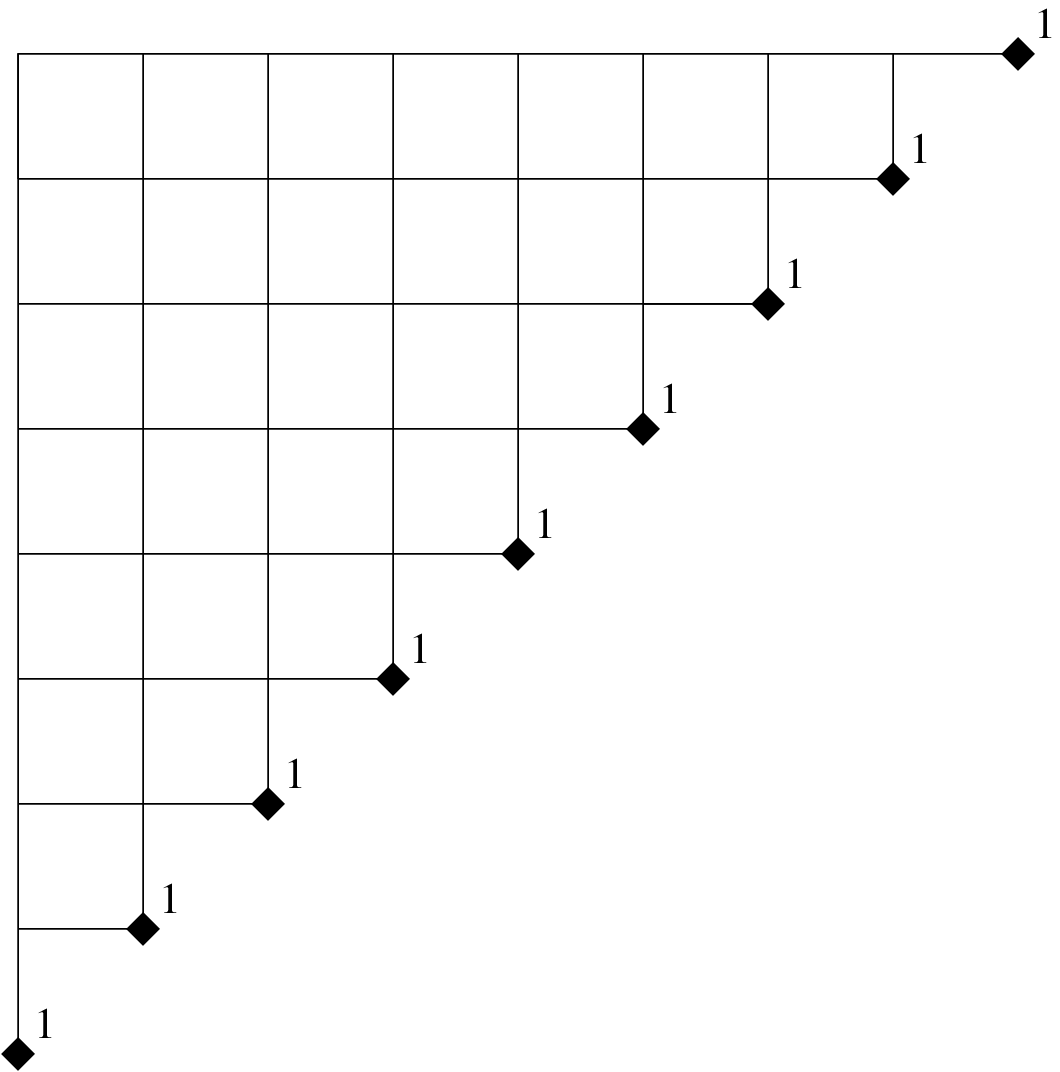}}{\mypic{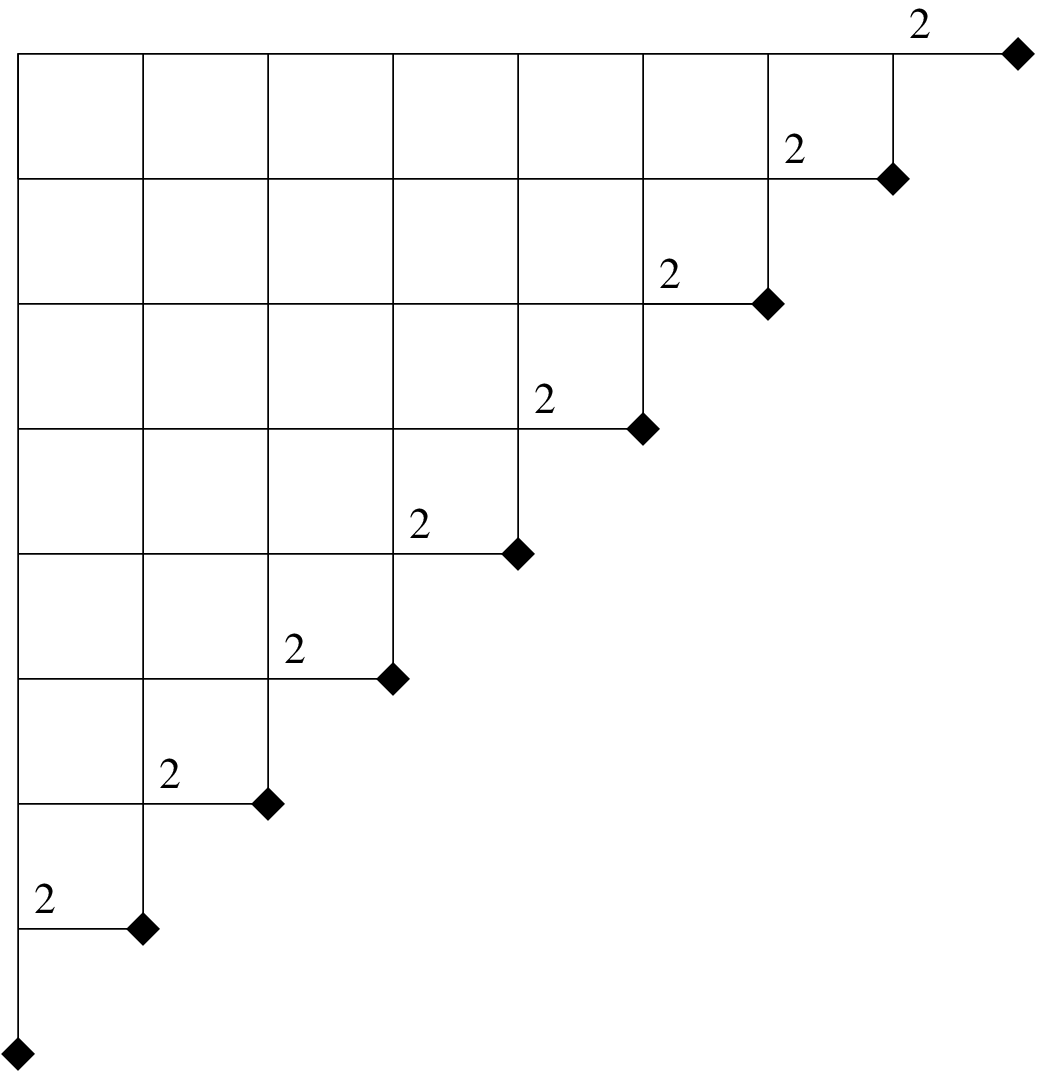}}{\mypic{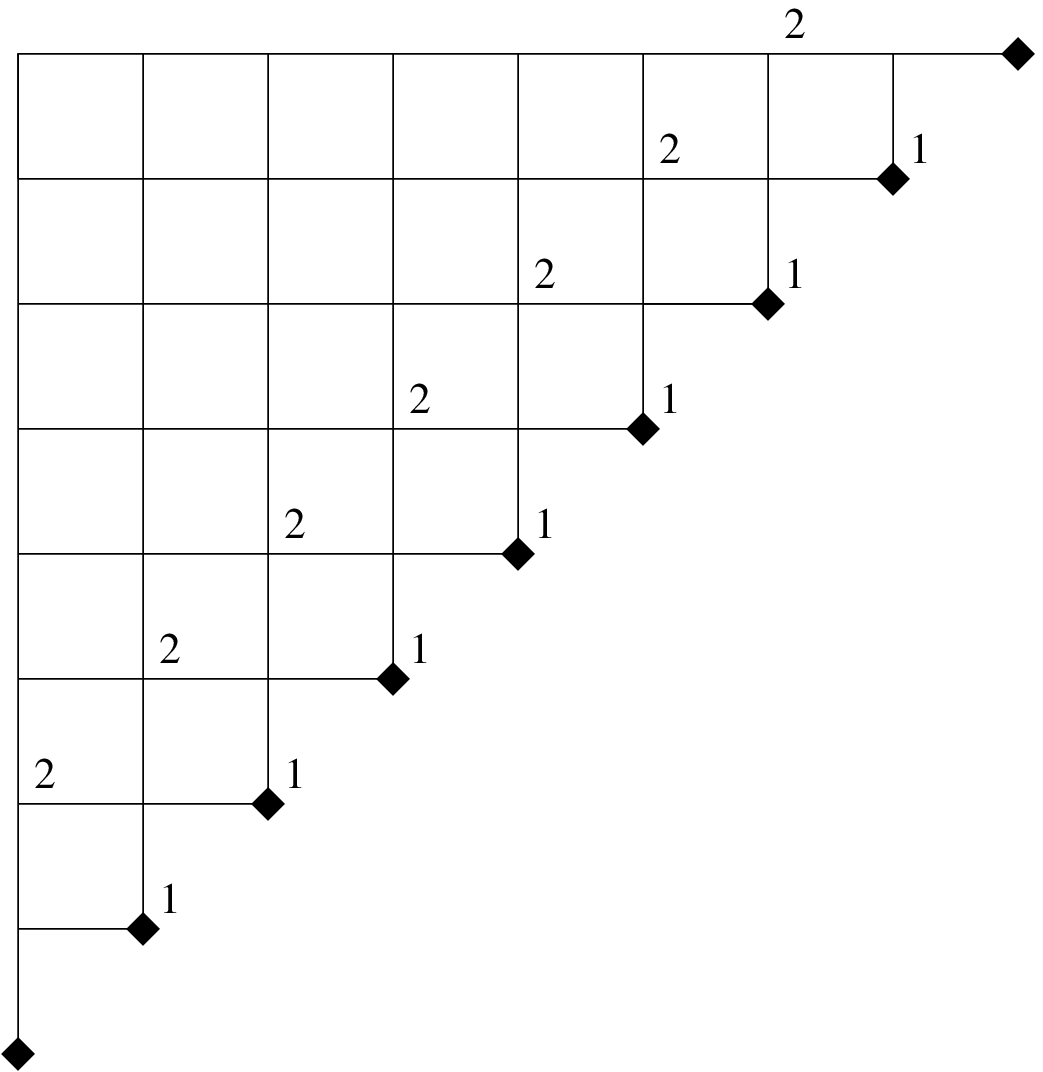}}
\threelinetext{Figure~2.7{\rm $(i)$.\phantom{aaaaaaaaaa}}}
{Figure~2.7{\rm $(ii)$.\phantom{aaaaaaaaaaa}}}{Figure~2.7{\rm $(iii)$.}}

\threeline{\mypic{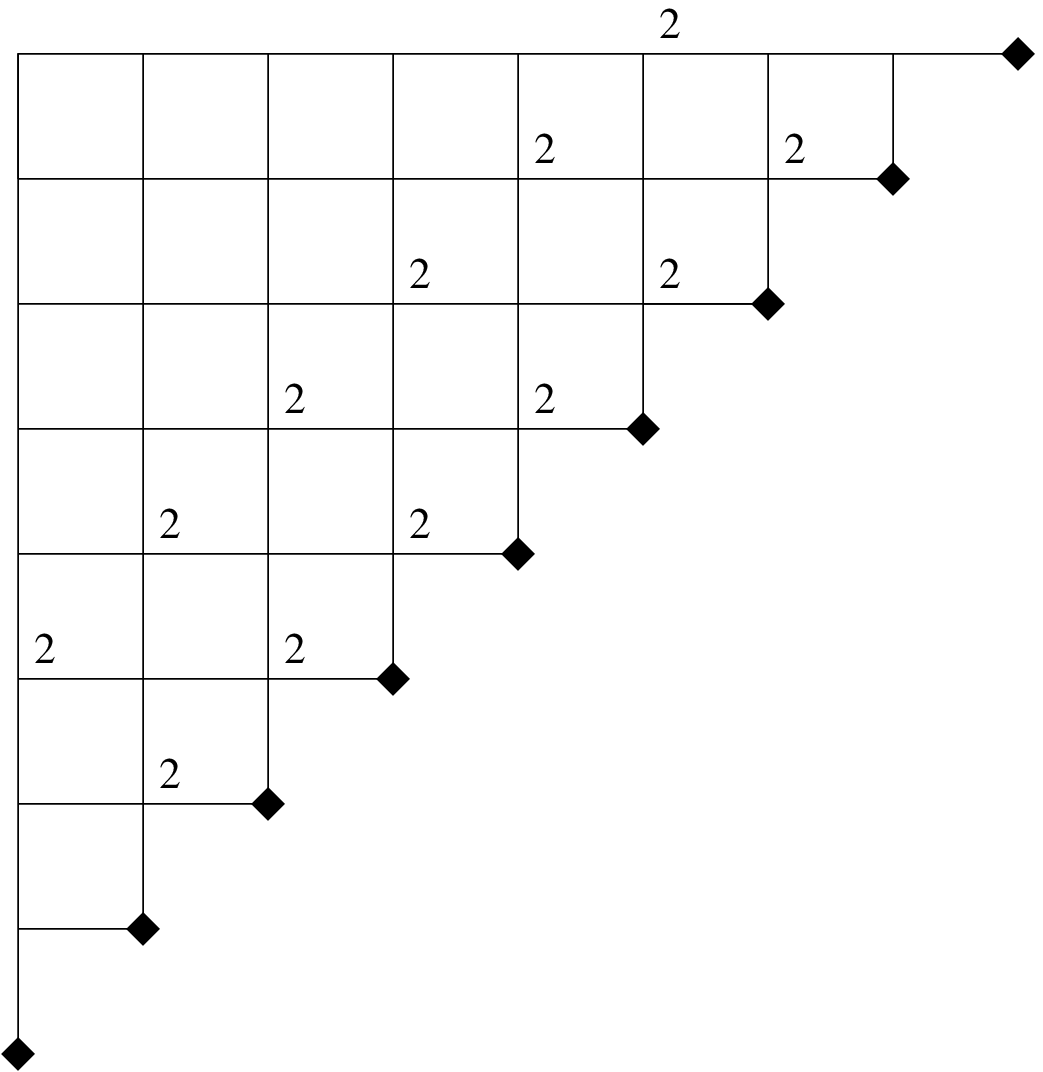}}{\mypic{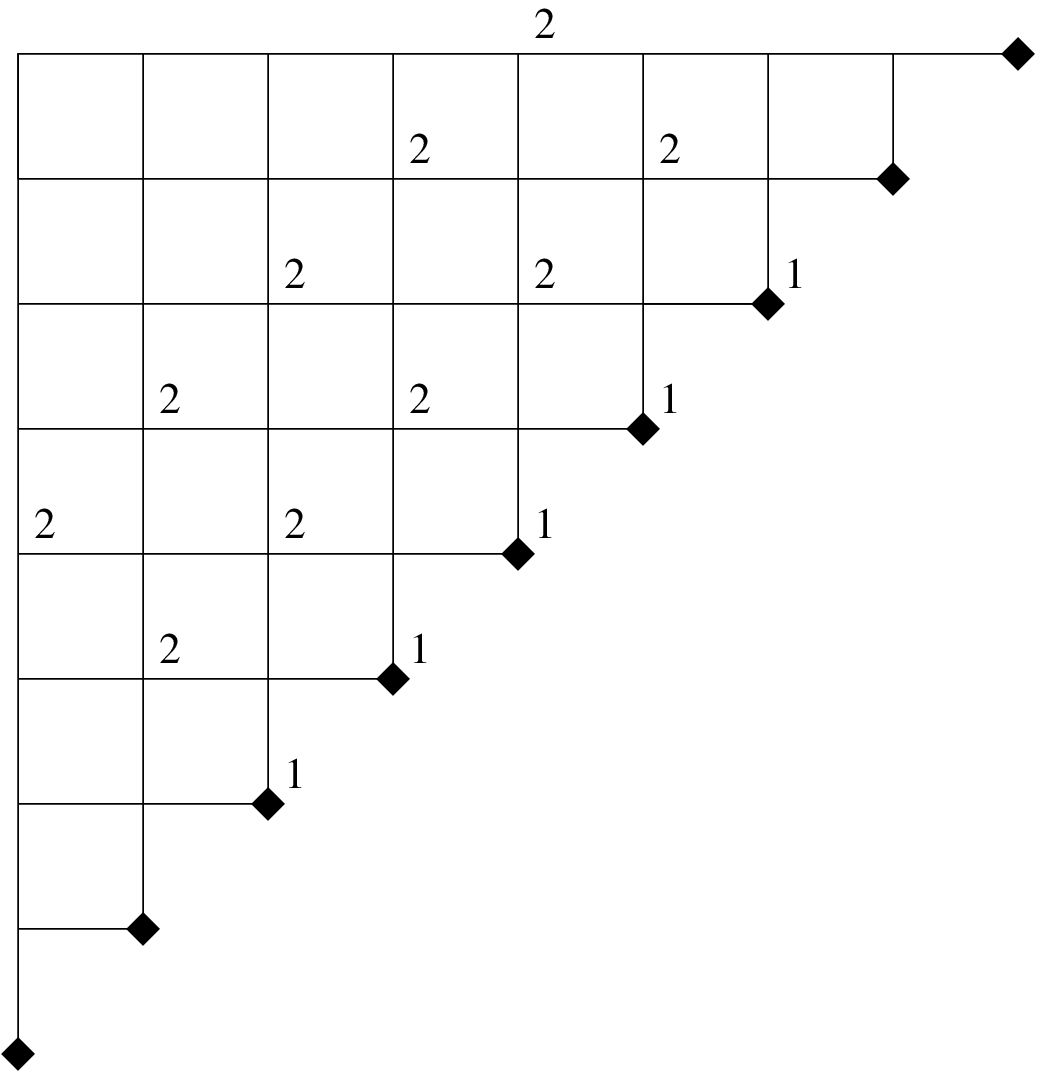}}{\mypic{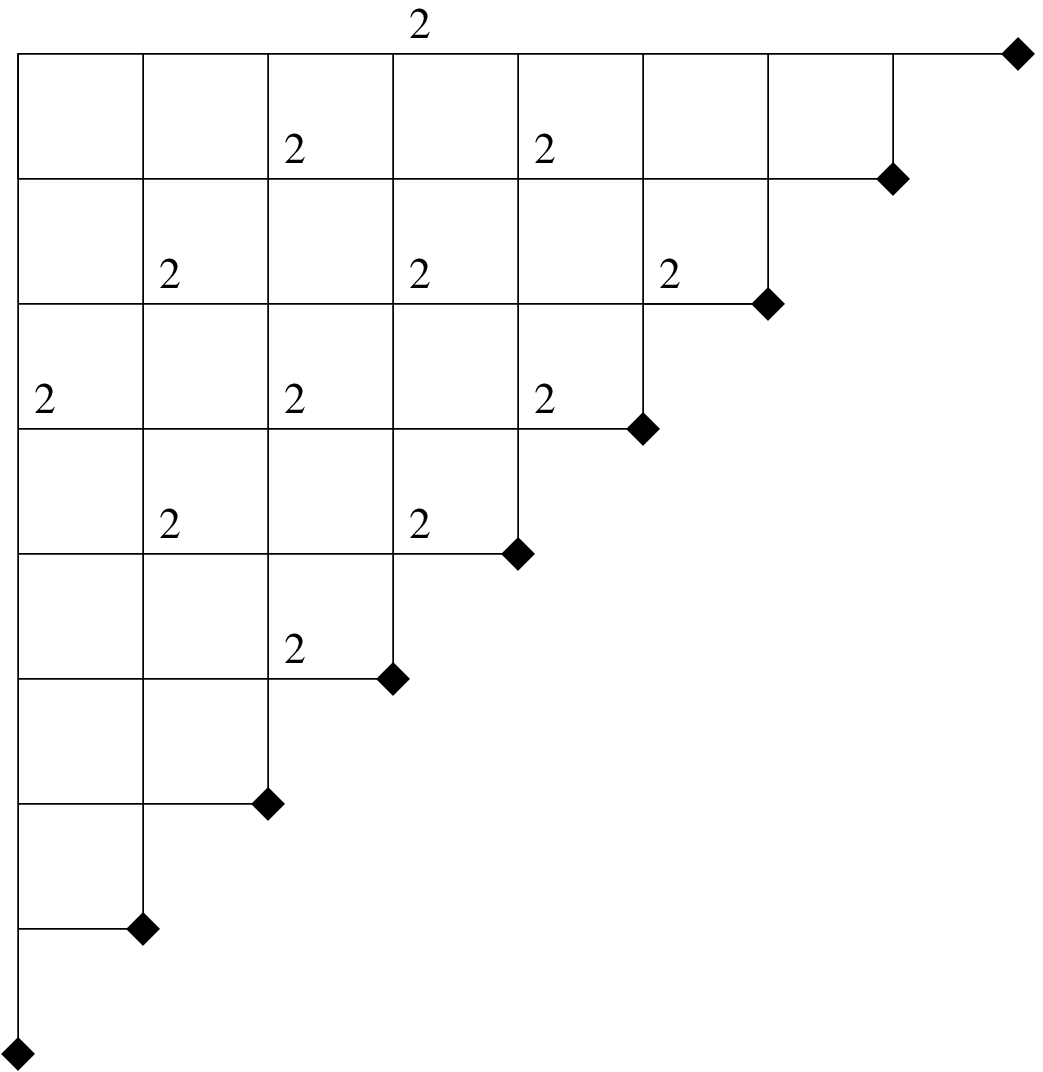}}
\threelinetext{Figure~2.7{\rm $(iv)$.\phantom{aaaaaaaaaa}}}
{Figure~2.7{\rm $(v)$.\phantom{aaaaaaaaaaa}}}{Figure~2.7{\rm $(vi)$.}}

\threeline{\mypic{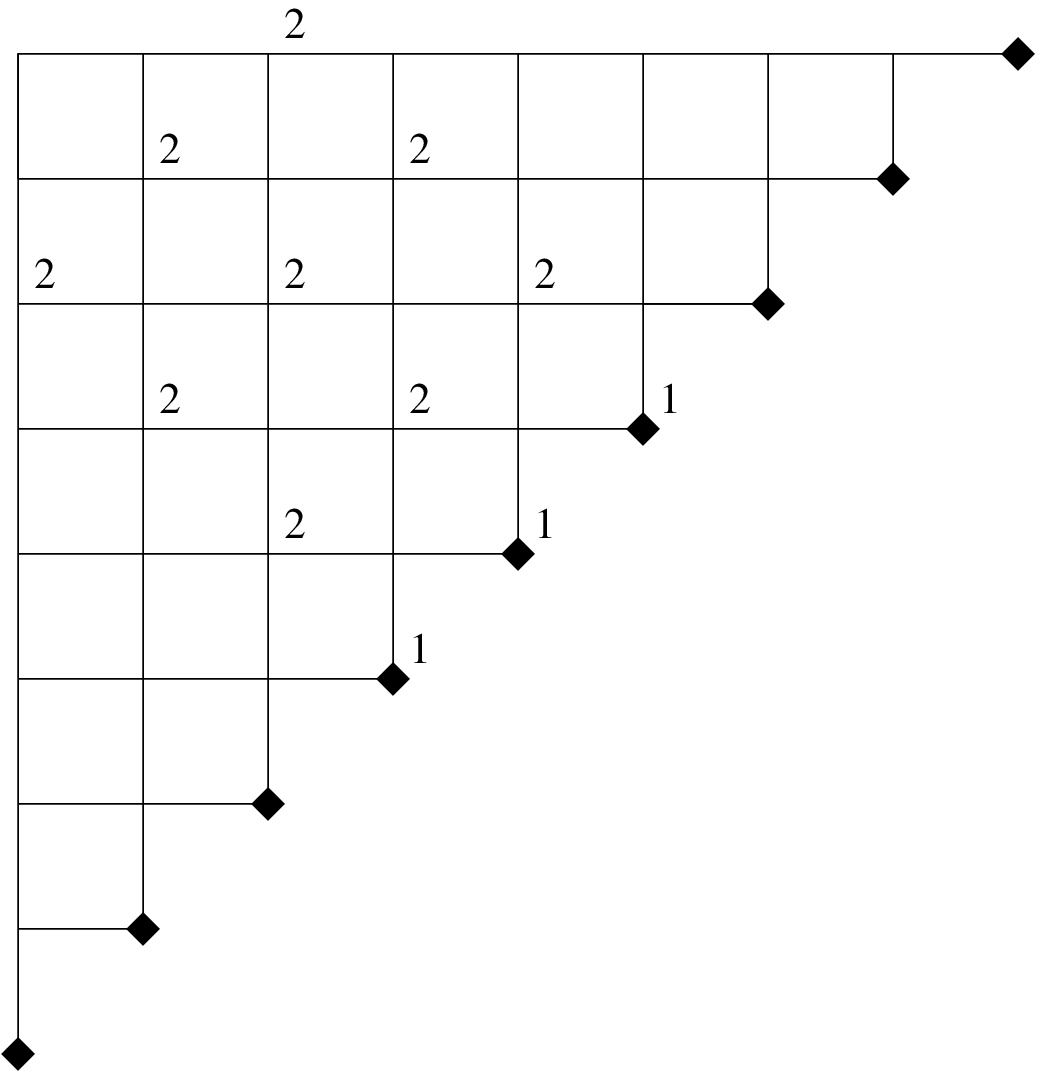}}{\mypic{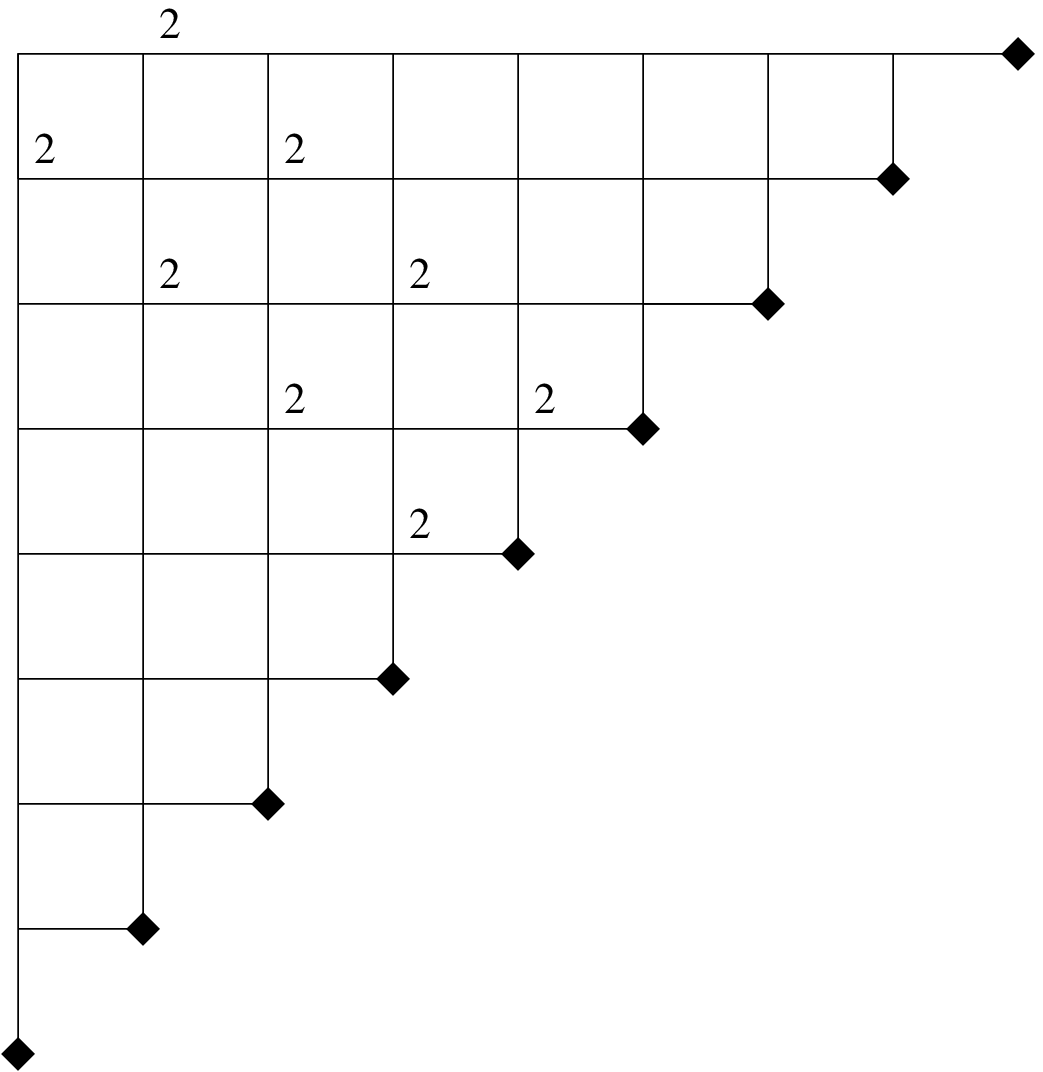}}{\mypic{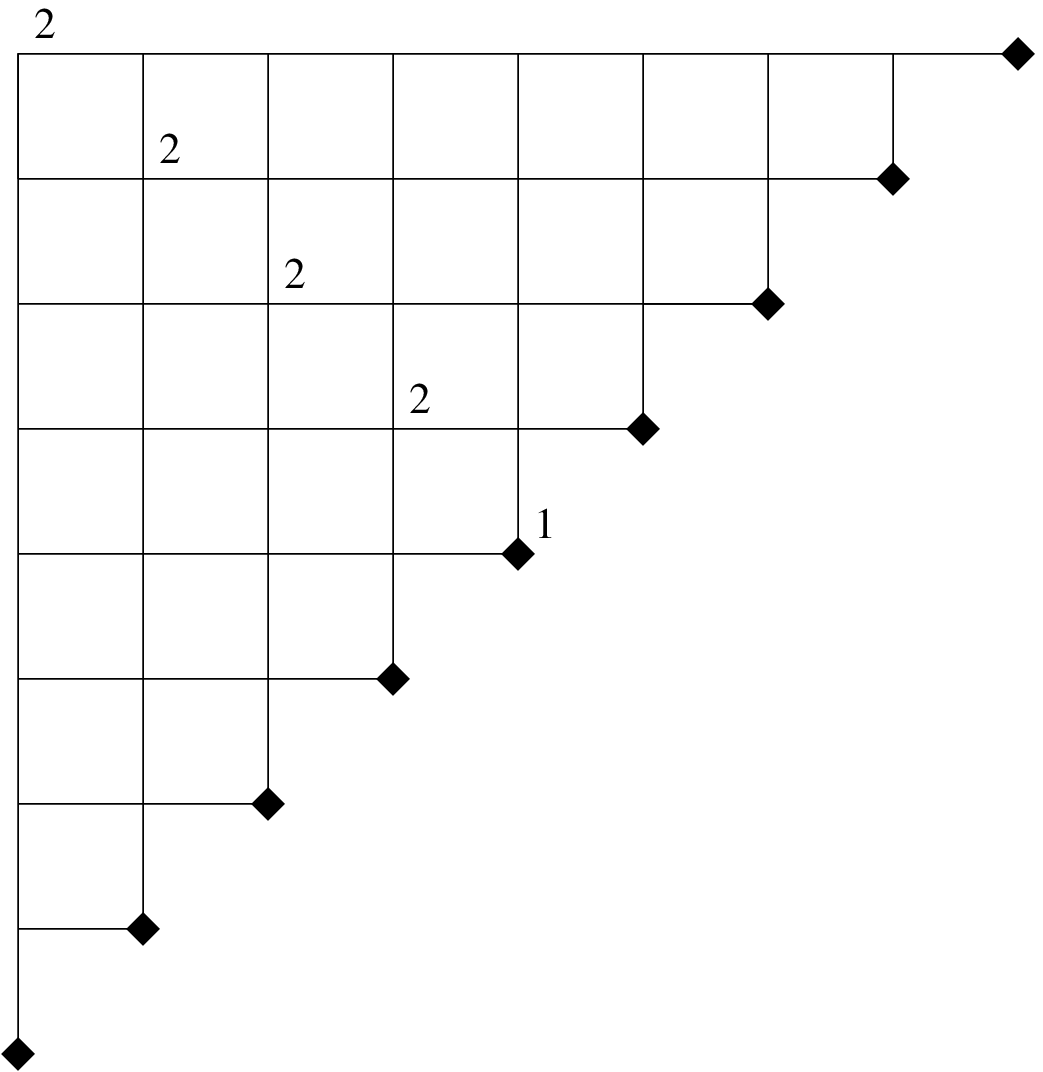}}
\threelinetext{Figure~2.7{\rm $(vii)$.\phantom{aaaaaaaaaa}}}
{Figure~2.7{\rm $(viii)$.\phantom{aaaaaaaaaaa}}}{Figure~2.7{\rm $(ix)$.}}

\threeline{\mypic{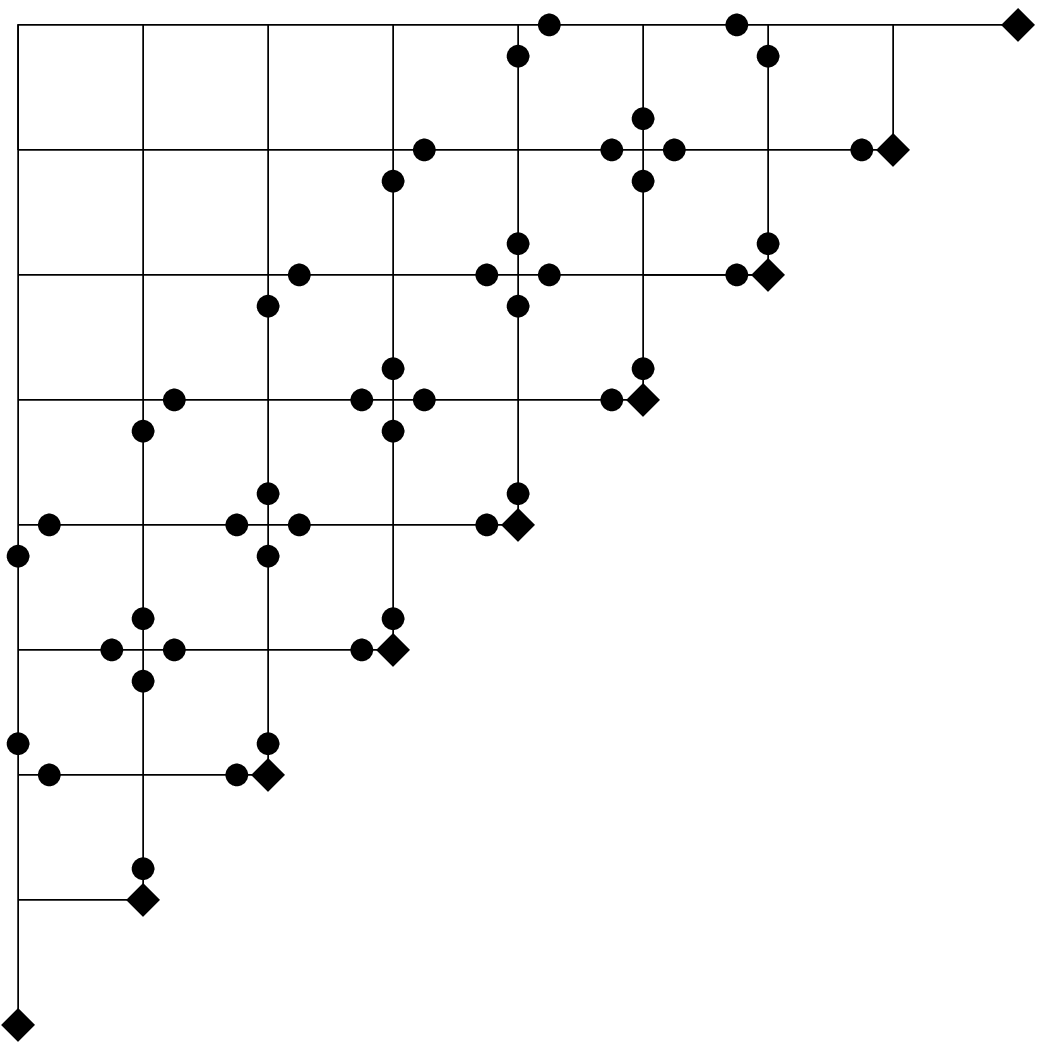}}{\mypic{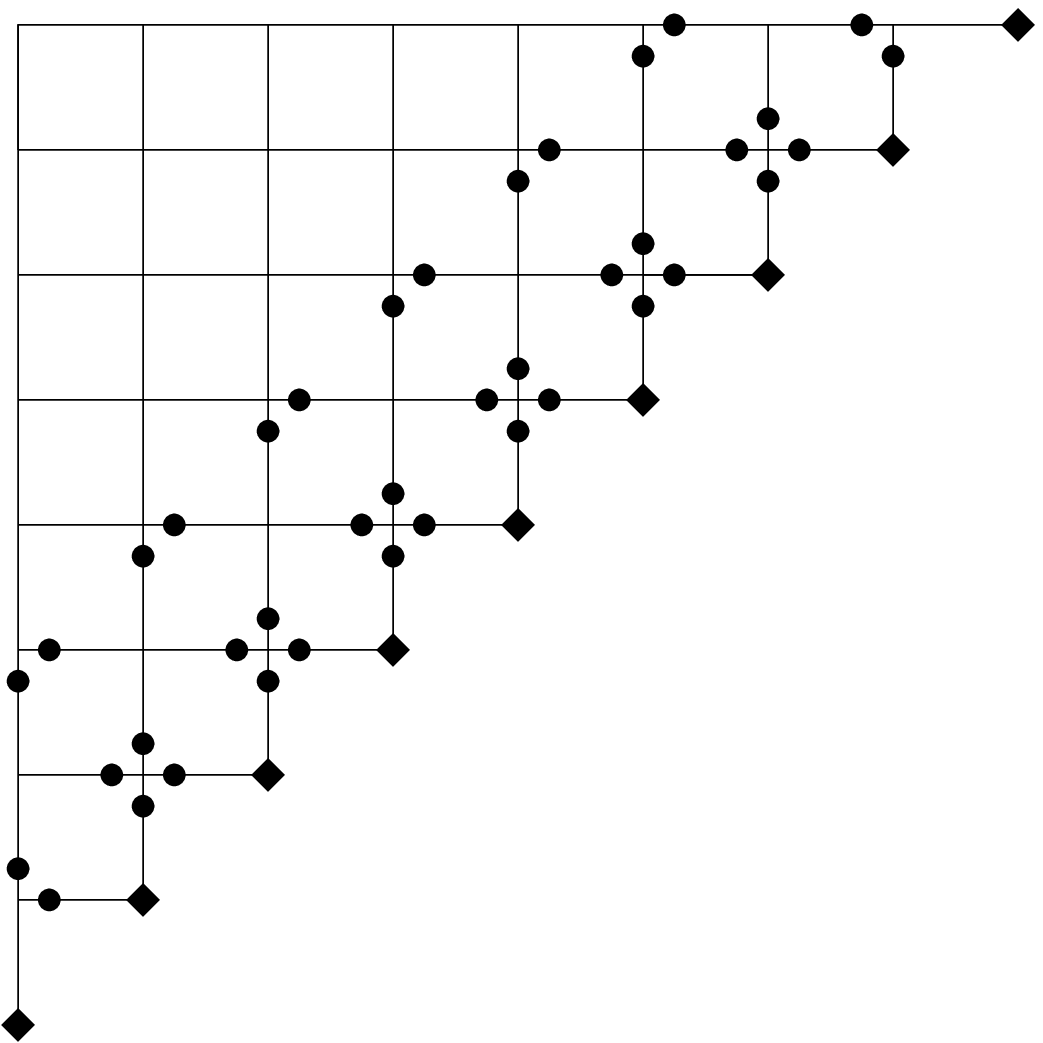}}{\mypic{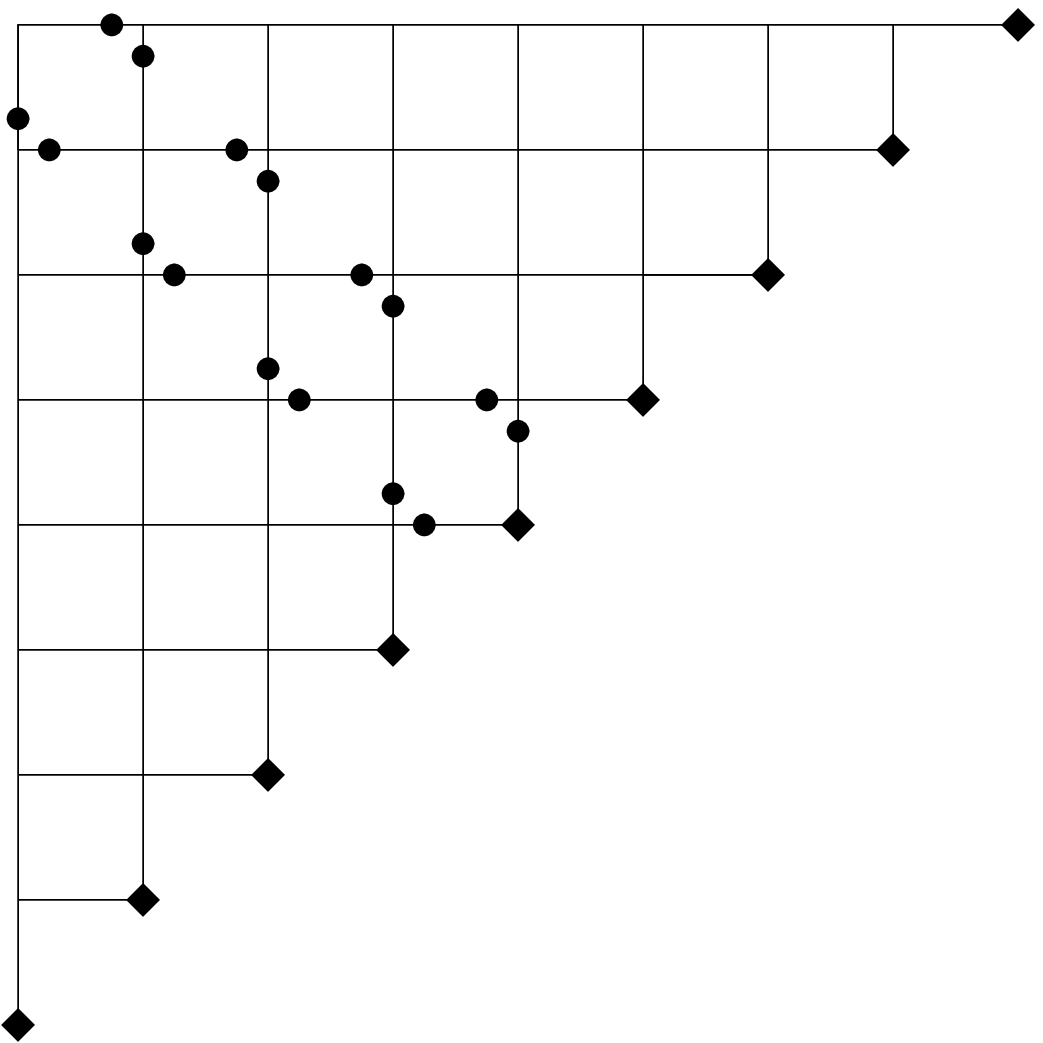}}
\threelinetext{Figure~2.7{\rm $(x)$.\phantom{aaaaaaaaaa}}}
{Figure~2.7{\rm $(xi)$.\phantom{aaaaaaaaaaa}}}{Figure~2.7{\rm $(xii)$.}}

\endinsert

\pf This follows again by the linearity of $F$ and the fact that its action on $e_{ij}$ is 
the weighted sum of the indeterminates corresponding to the neighbors of vertex $(i,j)$ in 
$\overline{QAD}_n$. Figure 2.7 corresponds to
$n=9$, but all the features needed to check the general case are present in it. Panels $(i)$--$(ix)$
show the vectors $w_1,\dotsc,w_9$, respectively: the vertices that contribute indeterminates with 
non-zero coefficients are indicated by the value of that coefficient next to them. 

Panel $(x)$ shows 
$\frac12 F(w_4)$---one dot next to a vertex $(k,l)$ represents the 
contribution of one unit to the coefficient of $e_{kl}$. It is apparent that this vector is
the same as the sum of the vectors in panels $(iii)$ and $(v)$. This checks (2.13) for $k=4$.

Similarly, panels $(ii)$, $(iv)$, and $(xi)$ check (2.13) for $k=3$, while the case $k=n=9$ follows
by panels $(viii)$ and $(xii)$. \epf

\topinsert
\twoline{\mypic{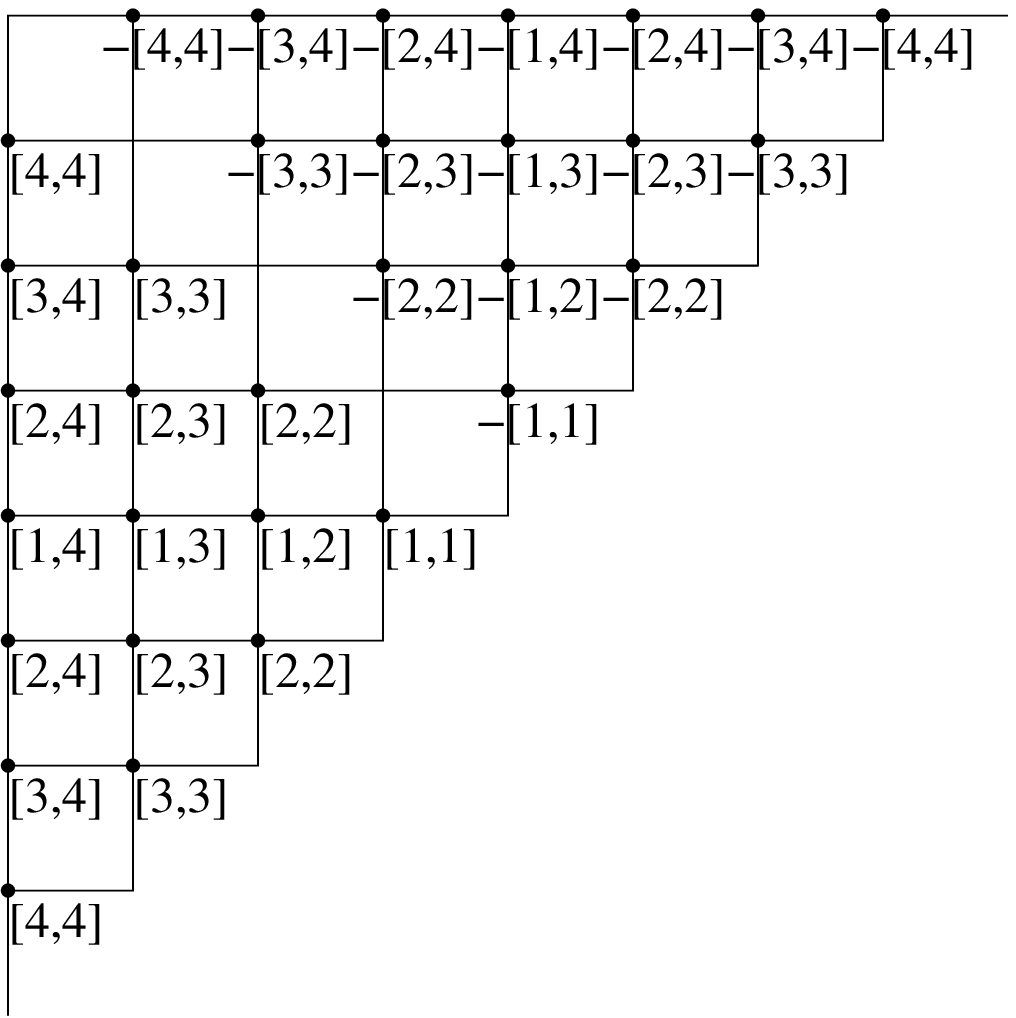}}{\mypic{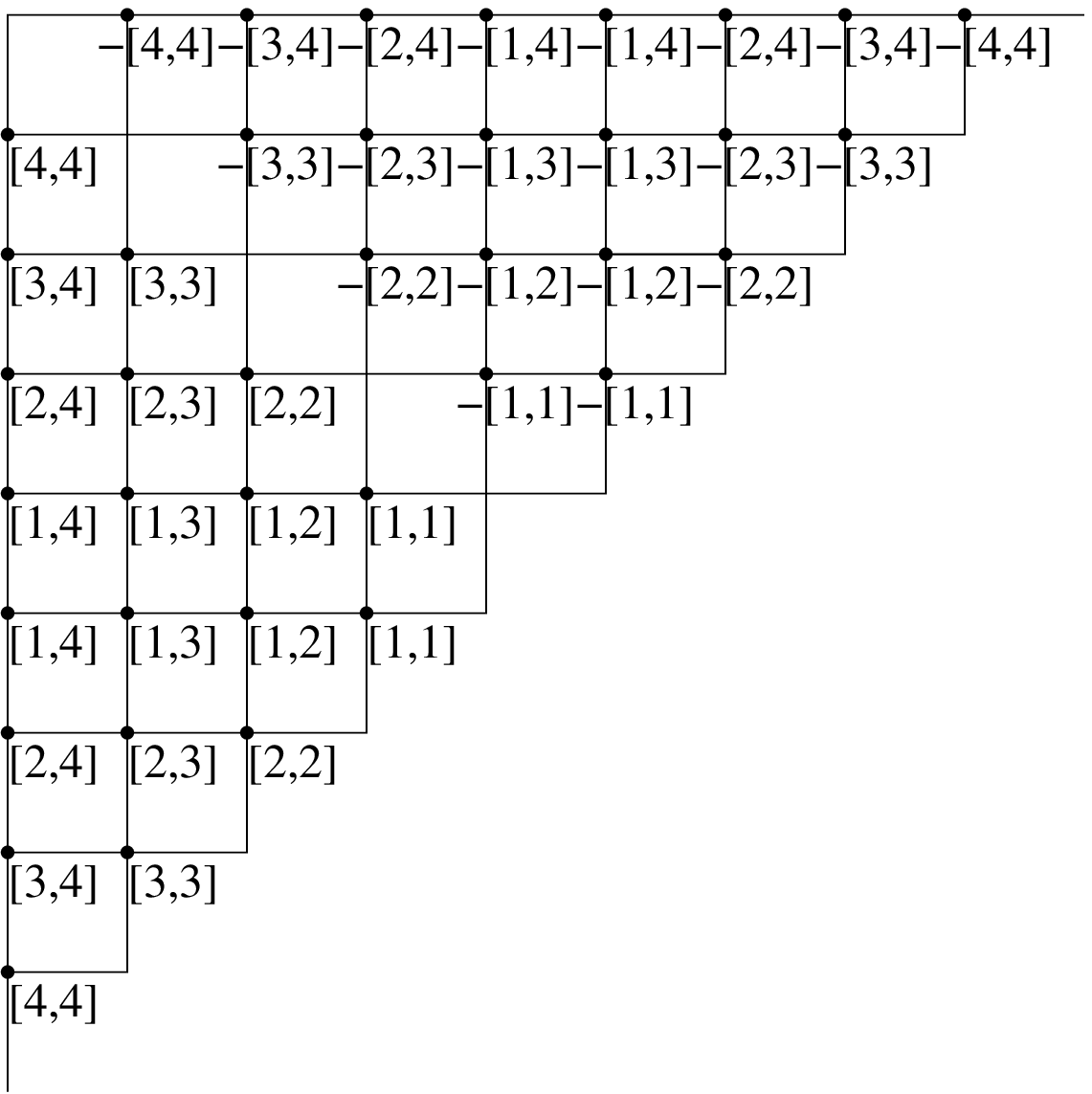}}
\twoline{Figure~2.8{\rm (a)}. {\rm The weight function $\wt_9$.}}
{Figure~2.8{\rm (b)}. {\rm The weight function $\wt_{10}$.}}
\endinsert

\topinsert
\centerline{\mypic{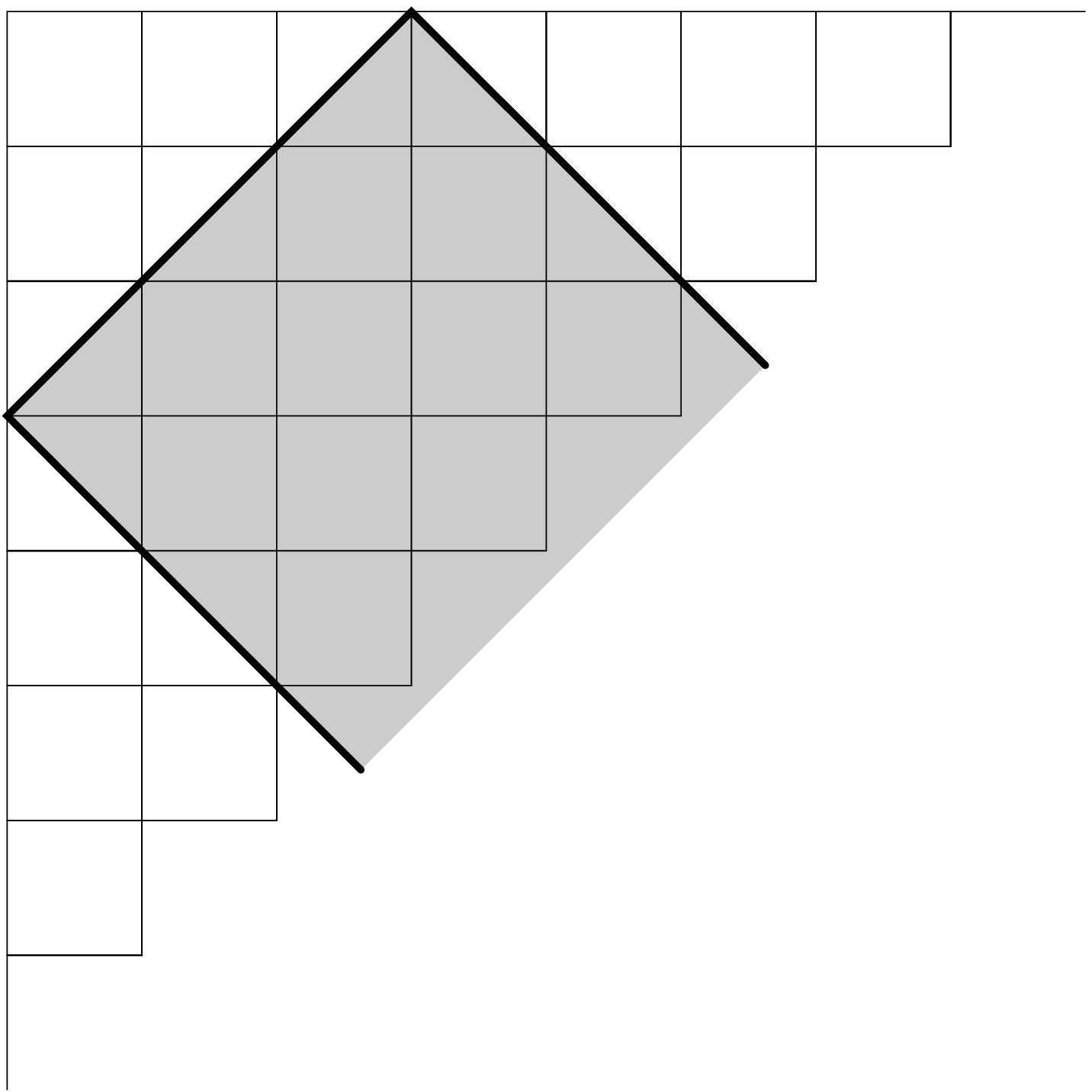}}
\centerline{Figure~2.9. {\rm $C_4$.}}
\endinsert

\proclaim{Lemma 2.6} The vectors $(2.8)$ and $(2.12)$ form a basis of $\bold U$.
\endproclaim

\pf Since there are $n(n-1)/2$ vectors $v_{ij}$ in (2.8) and $n$ vectors $w_k$ in (2.12), it suffices to 
show that $e_{ij}$ is in their span for all $i,j\geq1$, $i+j\leq n+1$. 
We show first that the
vectors $u_k:=\sum_{i+j=k+1} e_{ij}$,  $k=1,\dotsc n$, are in this span. 

Indeed, for $1\leq k\leq l\leq n$ define
$$
\spreadmatrixlines{3\jot}
[k,l]:=\left\{\matrix (2k-1)+(2k+1)+\cdots+(2l-1),&&\ \ \ \text{\rm if $n$ is odd,}\\
                      \!\!\!\!\!\!(2k)+(2k+2)+\cdots+(2l),&&\ \ \ \text{\rm if $n$ is even.}
\endmatrix\right.
$$
Let $\wt_n$ be the weight function on the vertices of $\overline{QAD}_n$ whose nonzero values are
given by the patterns shown in Figure 2.8---Figure 2.8(a) shows the pattern for odd values of $n$, and
Figure 2.8(b) the pattern for even $n$ (vertices of nonzero weight are marked by dots).

Define the {\it closed half-strip} $C_k$ to be the region described, in the standard rectangular coordinate system 
with origin at the northwest corner of $\overline{QAD}_n$, by
$\{(x,y): -(k-1)\leq x+y\leq k-1, x-y\geq k-1\}$ (for $n=9$, $C_4$ is illustrated in Figure 2.9).
Color the vertices of $\overline{QAD}_n$ in chessboard fashion. 
Denote the set of vertices of $\overline{QAD}_n$ by $V$.

We claim that for $k=1,\dotsc,n-2$ one has\footnote{ For a vertex $x$ of $\overline{QAD}_n$ we write 
$v_x$ for $v_{ij}$, where $(i,j)$ are the coordinates of $x$.}  
$$
\sum_{{\scriptstyle x\in C_{k+1}\cap V}\atop{\scriptstyle  \text{\rm $x$ of same color as vertex (k,1)}}} 
\wt_n(x)v_x - kw_{n-1-k}+(k+2)w_{n+1-k}
=
(2n+4)(e_{k,1}+e_{k-1,2}+\cdots+e_{1,k}).\tag2.14
$$

Indeed, suppose that $n$ is odd. It is then readily seen that the definition (2.8) of the 
vectors $v_{ij}$ and the definition of the weight $\wt_n$ imply that when expressing the sum on the left 
hand side of (2.14) in terms of the vectors $e_{ij}$, the nonzero coefficients are given by the pattern 
shown in Figure~2.10 (the two pictures illustrate the cases of odd and even $k$). On the other hand,
when expressing the combination of the last two terms on the left hand side of (2.14) as a linear 
combination of the $e_{ij}$'s, the nonzero coefficients are given by the patterns indicated by 
Figure~2.11. It is apparent from these two figures that (2.14) holds. The case of even $n$ is justified
analogously.

By (2.14), the vectors $u_1,\dotsc,u_{n-2}$ are seen to be in the span of the vectors (2.8) and (2.12).
Clearly, $u_{n-1}=\frac{1}{2}w_2$ and $u_n=w_1$. Thus all $u_k$'s, $k=1,...,n$ are in the span of the 
vectors (2.8) and (2.12), as claimed at the beginning of this proof.

We now show by induction on $i+j\geq2$ that $e_{ij}$ is in the span $S$ of the vectors (2.8) and 
$\{u_1,\dotsc,u_n\}$.

If $i+j=2$, we must have $i=j=1$, and the claim follows as $e_{11}=u_1$. Assume that $e_{ij}$ is
in $S$ for all $i+j\leq k$. Then the definition (2.8) of $v_{k-i,i}$ implies that 
$e_{k-i,i+1}-e_{k-i+1,i}$ is in $S$, for $i=1,\dotsc,k-1$. Thus $e_{k-i,i+1}=e_{k,1}+s_i$ 
with $s_i\in S$, for  $i=1,\dotsc,k-1$. Therefore $u_k=e_{k,1}+e_{k-1,2}+\cdots+e_{1,k}=ke_{k,1}+s$, with
$s\in S$, implying $e_{k,1}\in S$. Since $e_{k-i,i+1}=e_{k,1}+s_i$, it follows that 
$e_{k-i,i+1}\in S$ for all $i=0,\dotsc,k-1$, and the induction step is complete. This completes the proof.
\epf

%

\medskip
\flushpar
{\smc Remark 2.7.} Theorem 2.1 shows considerably more than the fact that the characteristic polynomial
of the graph on the left hand side of (2.1) is equal to the product of the characteristic polynomials
of the graphs on the right hand side of (2.1): it shows that the corresponding adjacency matrices are
similar (i.e., they have the same Jordan form). The same remark applies to the results of
Section 3.

\topinsert
\twoline{\mypic{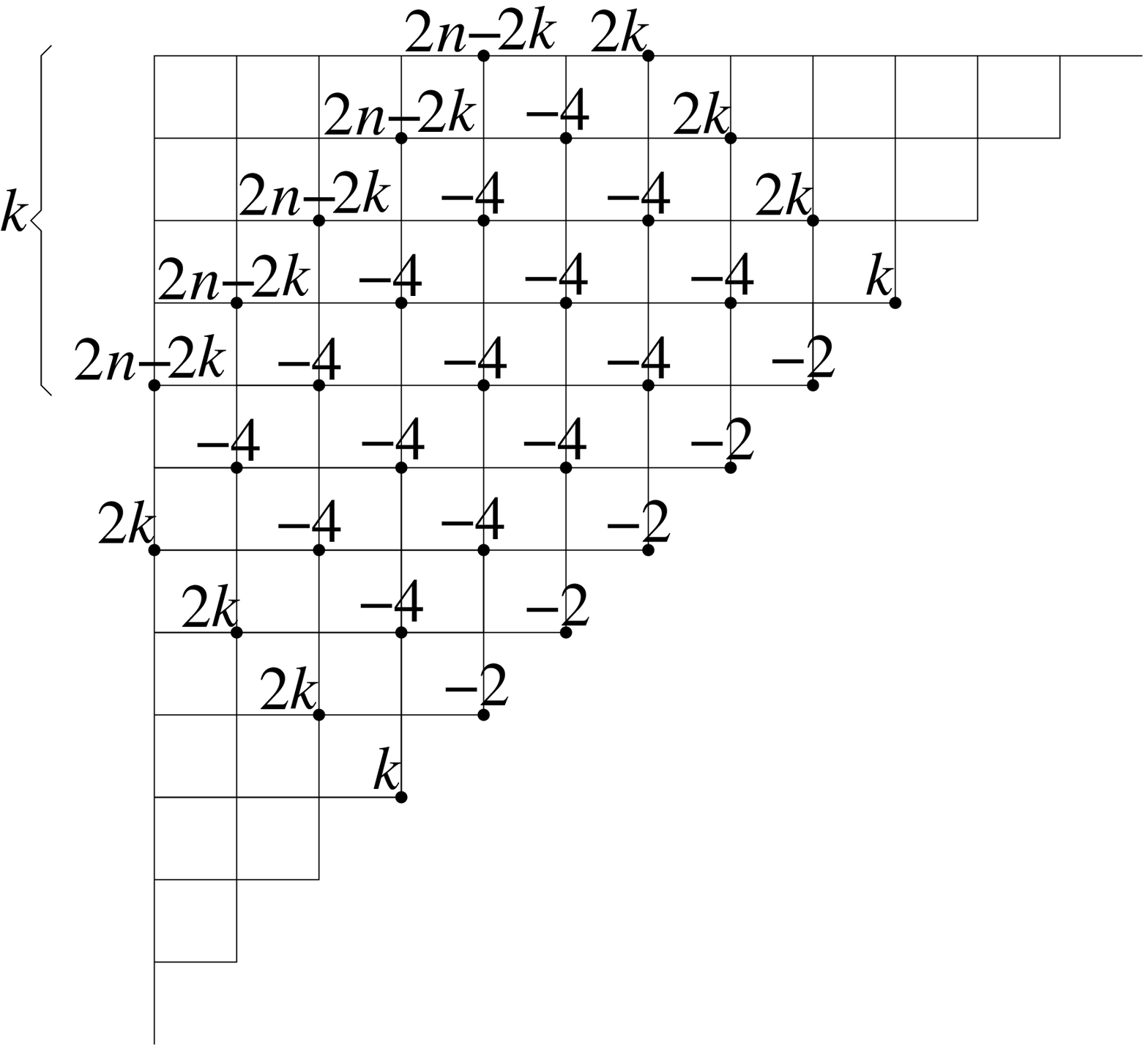}}{\mypic{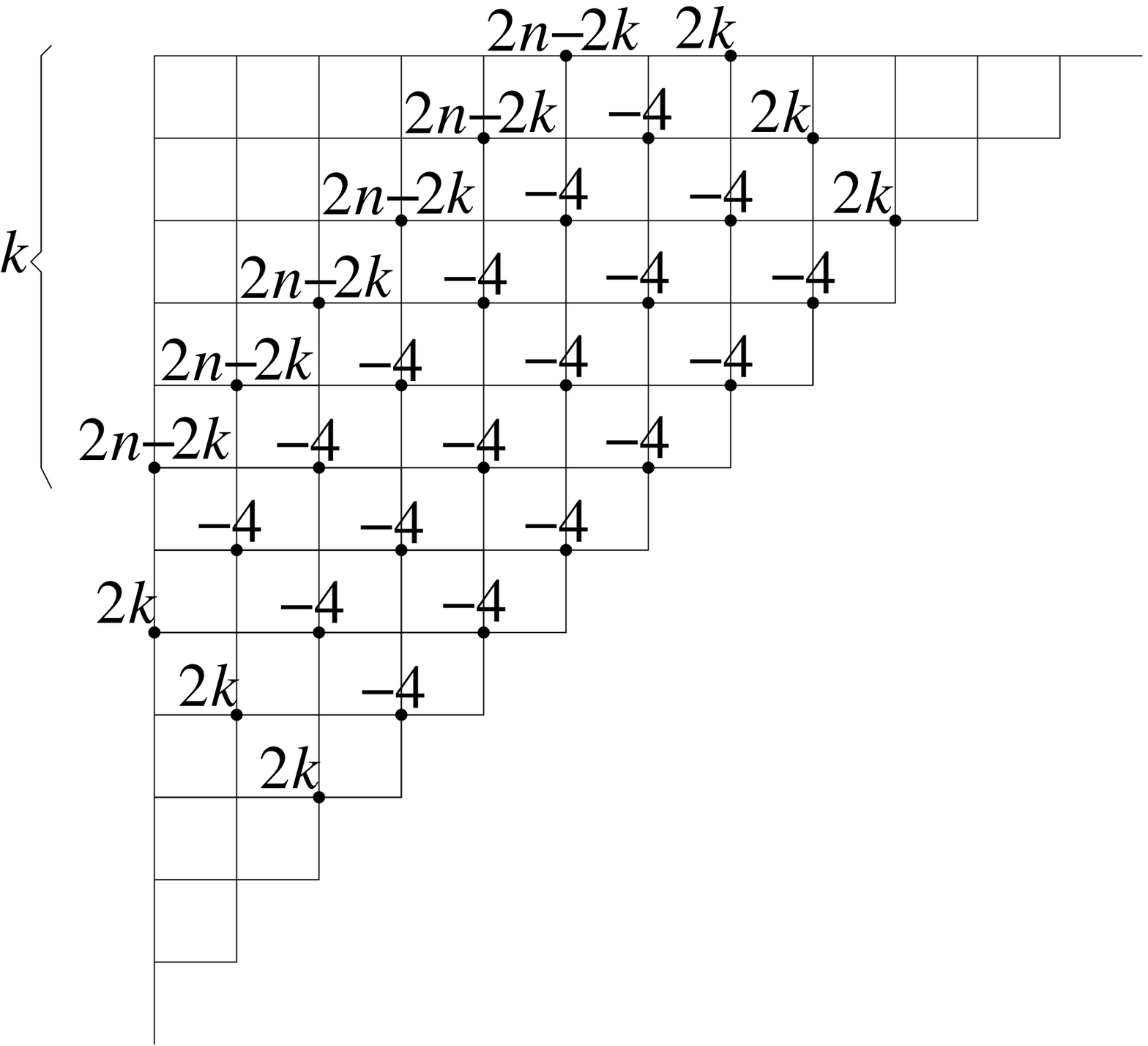}}
\centerline{Figure~2.10. {\rm (a). The sum in (2.14) for $n$ and $k$ odd. (b). Same for $n$ odd and 
$k$ even.}}
\endinsert

\topinsert
\twoline{\mypic{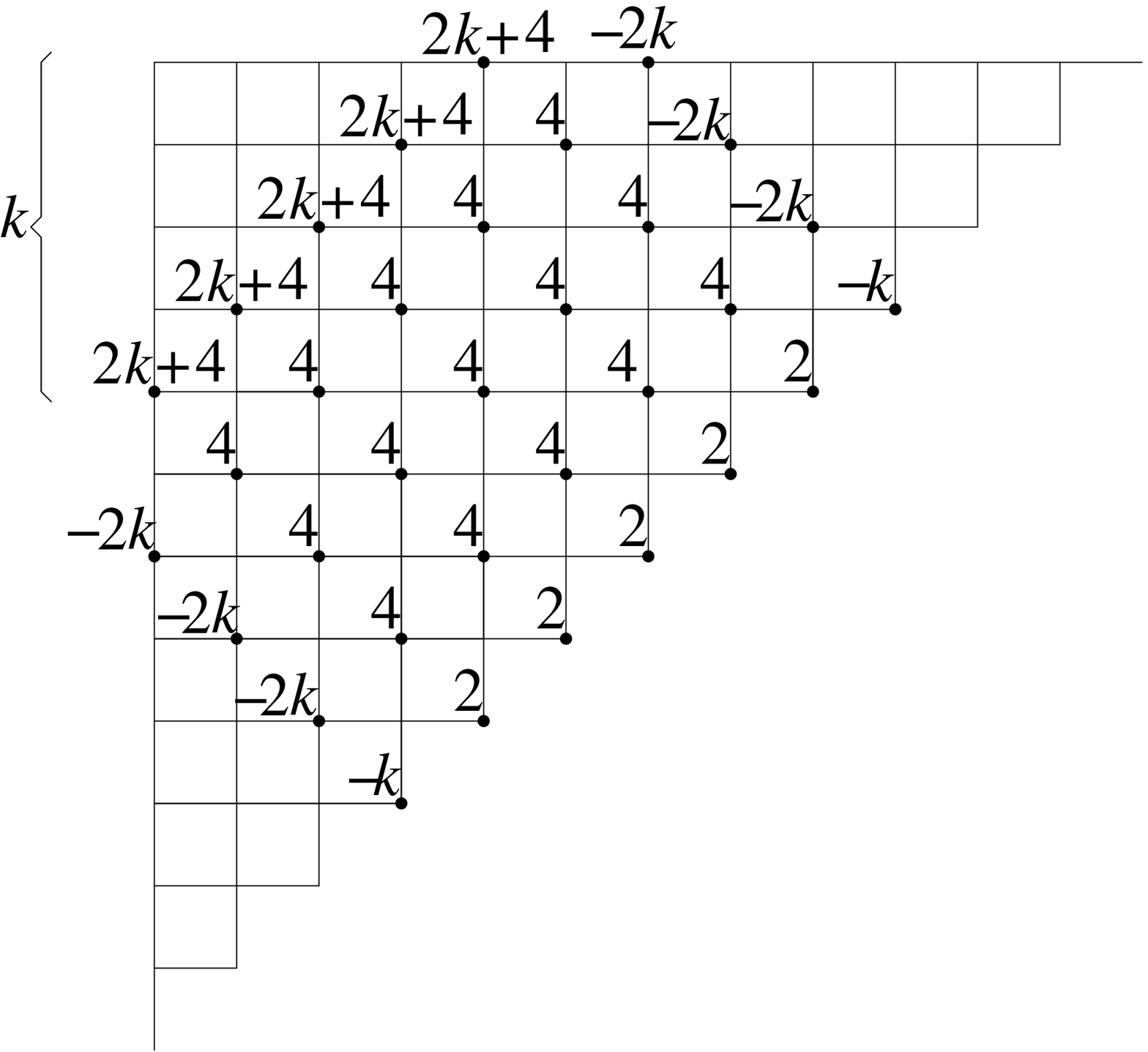}}{\mypic{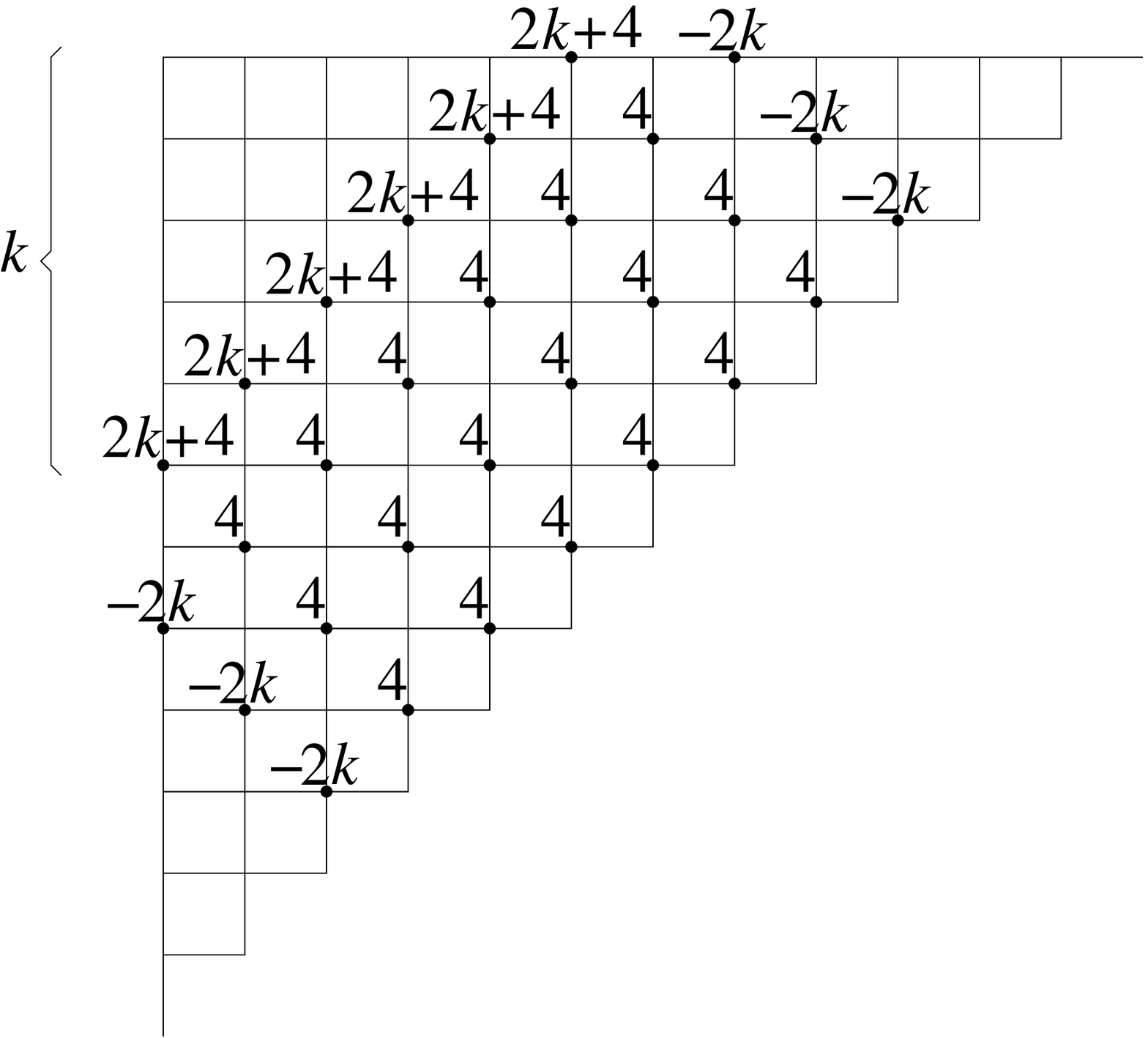}}
\centerline{Figure~2.11. {\rm (a). The sum of the last two terms on the left hand side of (2.14) for $n$ 
and $k$ odd.}}
\centerline{{\rm (b). Same for $n$ odd and $k$ even.}}
\endinsert



\mysec{3. The halved Aztec diamonds}

One way to define the Aztec diamond $AD_n$ is to say that it is the graph whose vertices are the
white unit squares of a $(2n+1)\times(2n+1)$ chessboard with black corners, two vertices being
adjacent if they correspond to white unit squares that share a corner. The analogous graph on the
black unit squares of this chessboard is called the {\it odd Aztec diamond} of order $n$ and is denoted
$OD_n$ (Figure 3.1 shows $OD_5$).

Similar considerations on a $2n\times 2n$  chessboard lead to graphs on the white and black unit
squares that are isomorphic. We denote them by $MD_n$ (``mixed'' diamonds; $MD_5$ is pictured in
Figure 3.2).

Define $HOD_n$ to be the subgraph of $OD_n$ induced by the black unit squares on or above a diagonal of
the $(2n+1)\times(2n+1)$ chessboard. 
Similarly, let $HMD_n$ be the subgraph of $MD_n$ induced by the black vertices on or above the
black diagonal of the $2n\times 2n$ board ($HOD_5$ and $HMD_5$ are shown in Figures 3.3 and 3.4,
respectively).

Let $Q_n$ be the graph obtained from the path $P_n$ by including a loop of weight 2 at each vertex.
Let $Q'_n$ be graph obtained similarly from $P_n$, but by weighting each loop by $-2$.
The graph $R_n$ is obtained from $Q_n$ by changing the weight of the loop at the last vertex to 1, 
and $R'_n$ is obtained from $Q'_n$ by changing the weight of the loop at the last vertex to $-1$.

\topinsert
\twoline{\mypic{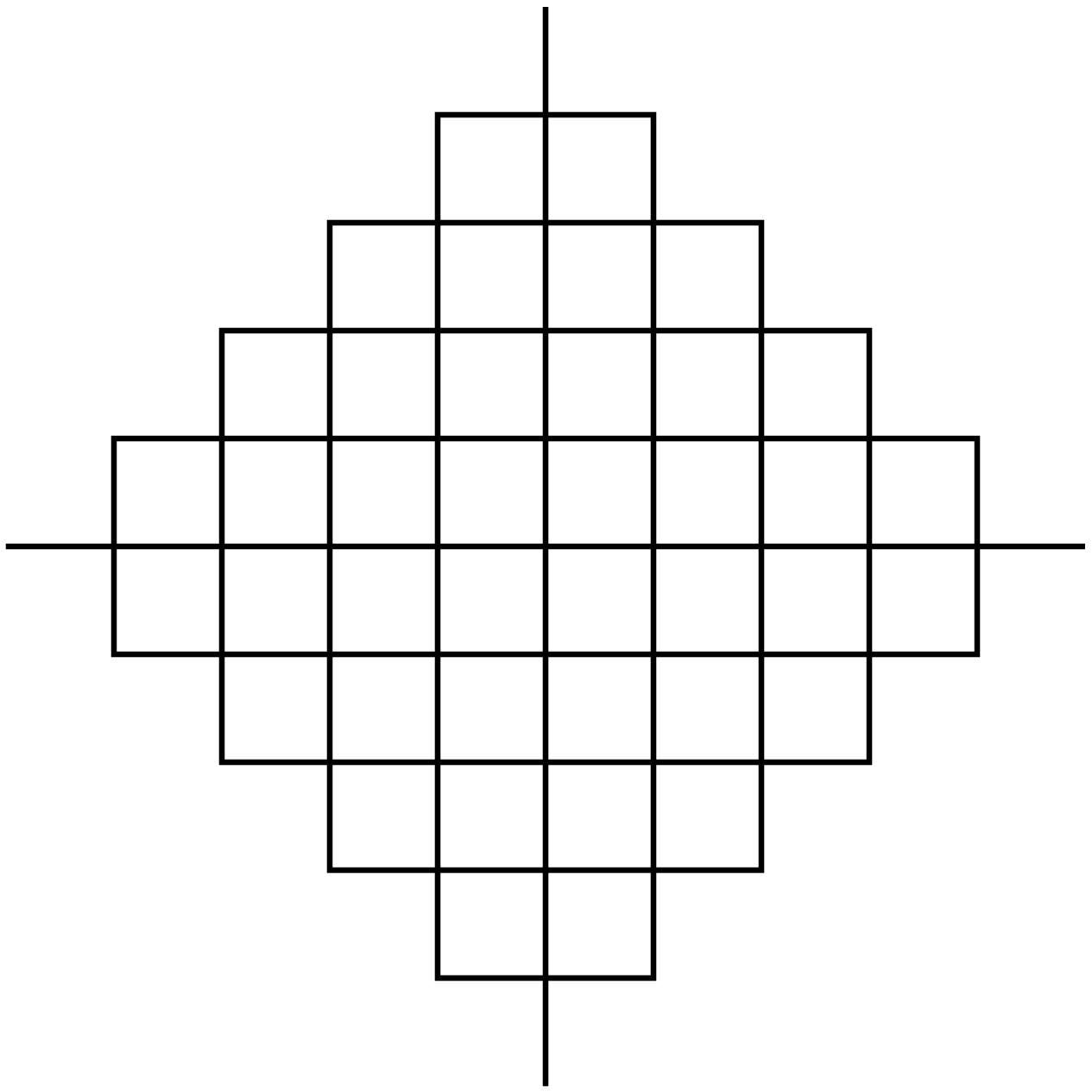}}{\mypic{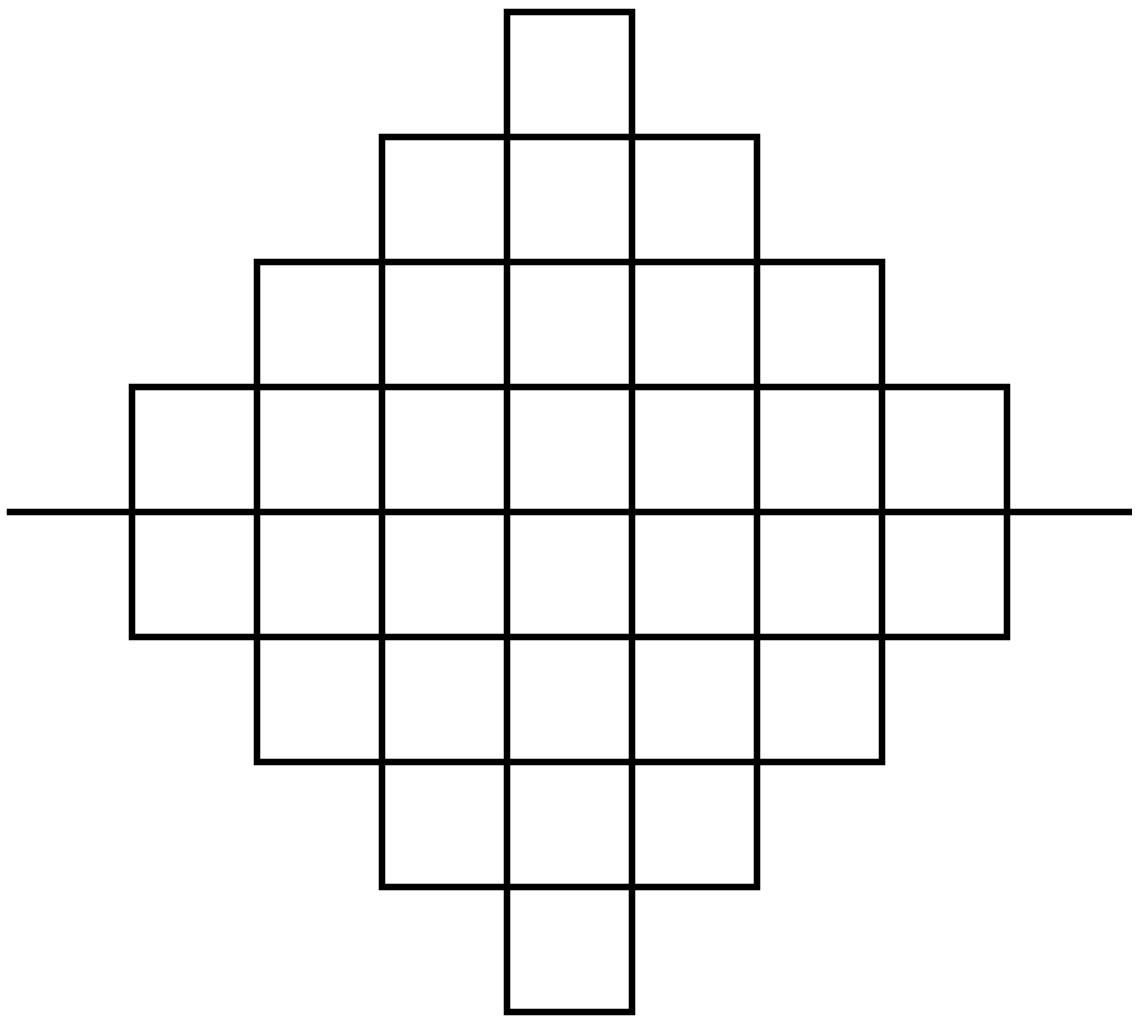}}
\twoline{Figure~3.1. {\rm $OD_5$.}}
{Figure~3.2. {\rm $MD_5$.}}

\bigskip
\twoline{\mypic{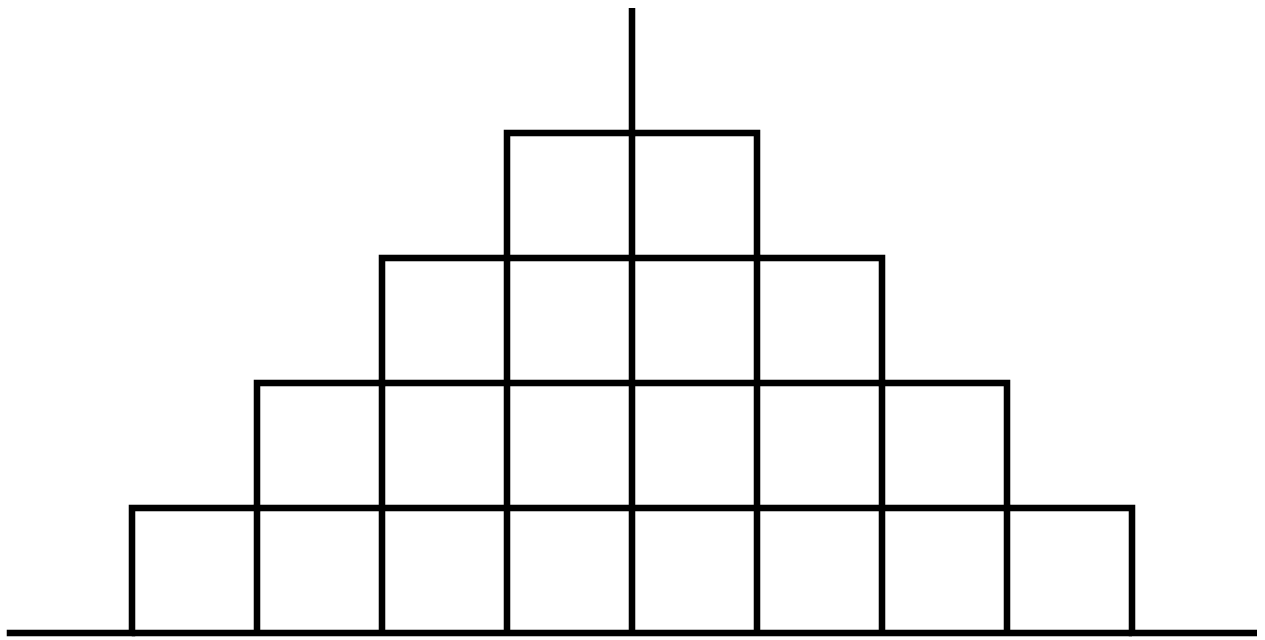}}{\mypic{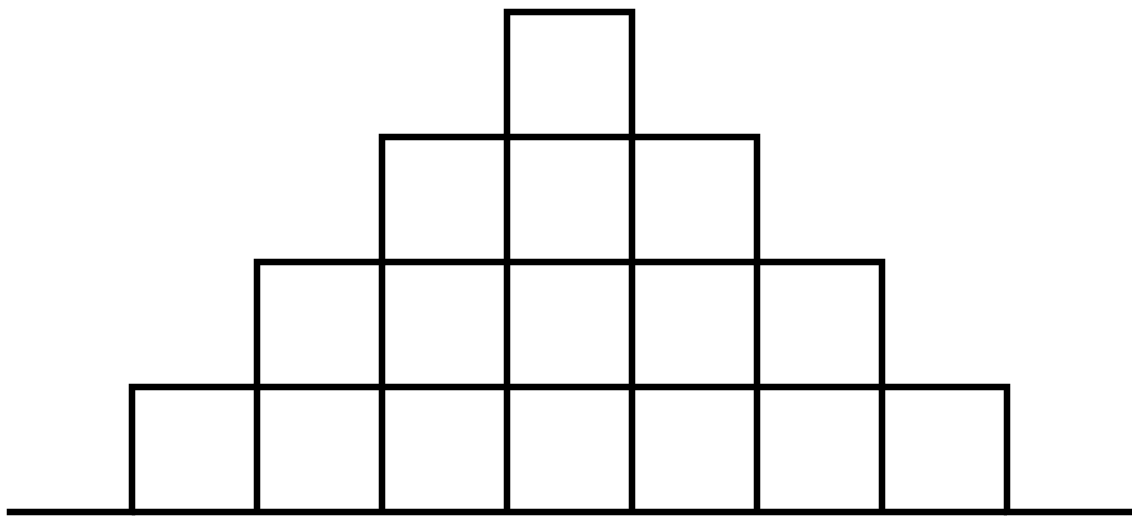}}
\twoline{Figure~3.3. {\rm $HOD_5$.}}
{Figure~3.4. {\rm $HMD_5$.}}

\bigskip
\twoline{\mypic{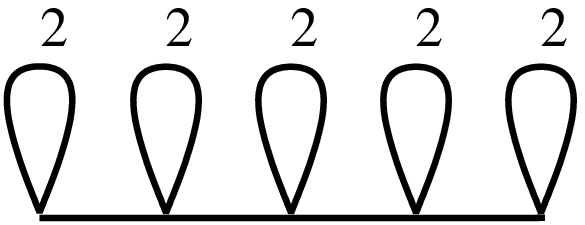}}{\mypic{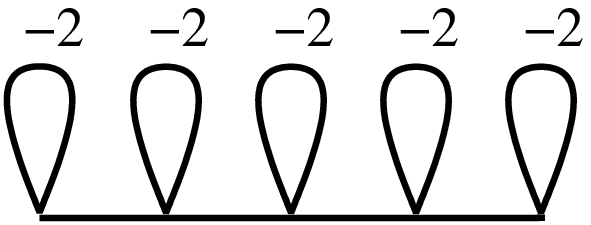}}
\twoline{Figure~3.5. {\rm $Q_5$.}}
{Figure~3.6. {\rm $Q'_5$.}}

\bigskip
\twoline{\mypic{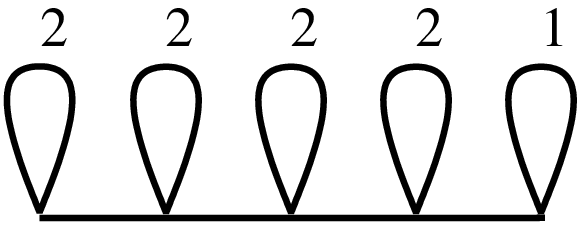}}{\mypic{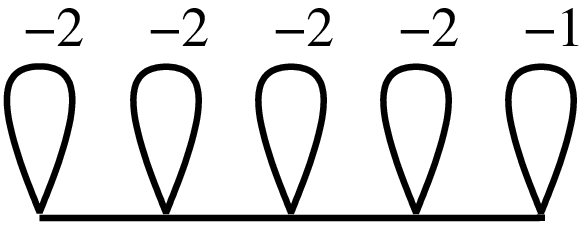}}
\twoline{Figure~3.7. {\rm $R_5$.}}
{Figure~3.8. {\rm $R'_5$.}}

\endinsert

\proclaim{Theorem 3.1} We have
$$
MD_n \sim HMD_{n-1}\ \dot{\cup}\ HMD_{n-1}\ \dot{\cup}\ R_n\ \dot{\cup}\ R'_n.\tag3.1
$$


\endproclaim


\pf Apply Lemma 1.1 to $MD_n$, choosing the automorphism $T$ to be the reflection across its
symmetry axis. The resulting graph $G^+$ is isomorphic to $HMD_{n-1}$. The resulting weighted
directed graph $G^-$, which we denote $\overline{HMD}_n$, is obtained 
from $HMD_n$ by: (1) marking the $2n$ bottommost vertices, 
(2) replacing each edge between two marked or two unmarked vertices by two anti-parallel directed edges 
of weight 1, and (3) replacing each edge between an unmarked vertex $v$ and a marked vertex $w$ by a 
directed edge $(v,w)$ weighted 1 and a directed edge $(w,v)$ weighted 2 (Figure 3.9 illustrates this
when $n=5$). We obtain
$$
MD_n\sim HMD_{n-1}\ \dot{\cup}\ \ \overline{HMD}_{n}.\tag3.2
$$

\topinsert
\centerline{\mypic{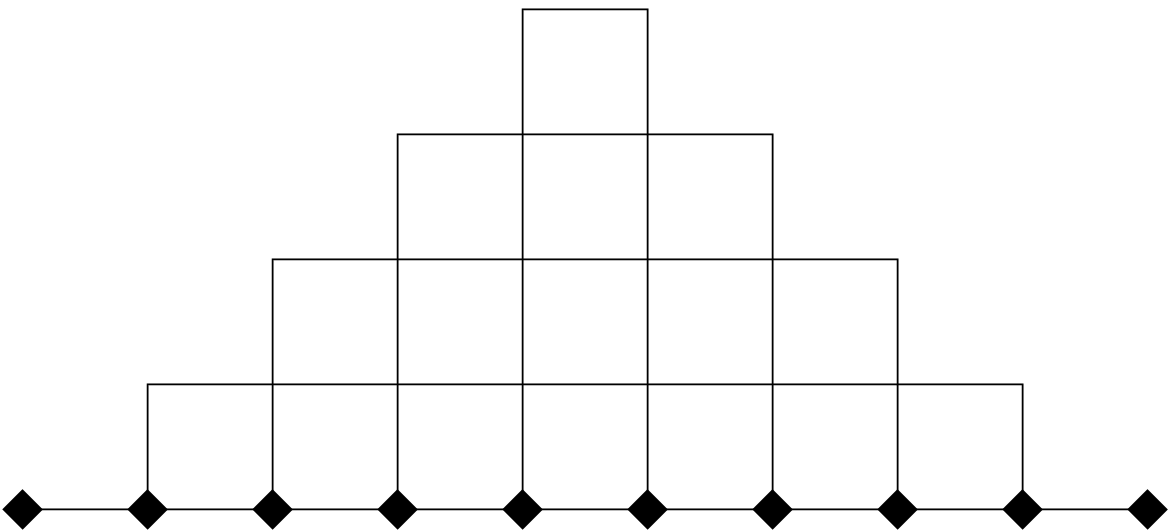}}
\centerline{{\smc Figure~3.9.} {\rm $\overline{HMD}_5$.}}
\endinsert

Let $v$ be a vertex of $\overline{HMD}_{n}$. Proceed from $v$ along a ray in the
northeast direction. If no other vertex of $\overline{HMD}_{n}$ is on this ray, set ${\NE(v)}:=\infty$; 
otherwise, let $w$ be the first encountered vertex and set ${\NE(v)}:=w$. Define $\NW(v)$, ${\SW(v)}$, 
and ${\SE(v)}$ analogously, via the vertices closest to $v$ in the northwestern, southwestern, and 
southeastern directions, respectively.

For each vertex $v$ of $\overline{HMD}_{n}$, let $e_v$ be an indeterminate.
Denote by $V_1$ the set of unmarked vertices of $\overline{HMD}_{n}$. Let $B_1\subset V_1$ consist of
the $n-1$ vertices on the northeastern side of the convex hull of $V_1$, and $B_2\subset V_1$ of
the $n-1$ vertices on the northwestern side of the latter. Define
$$
\spreadmatrixlines{3\jot}
f_v:=\left\{\matrix -e_{\NE(v)}+e_{\NW(v)}-e_{\SW(v)}+e_{\SE(v)}, && v\in V_1\setminus(B_1\cup B_2) \\
                    e_v+e_{\NW(v)}-e_{\SW(v)}+e_{\SE(v)}, && v\in B_1 \\
                   -e_v-e_{\NE(v)}-e_{\SW(v)}+e_{\SE(v)}, && v\in B_2,
\endmatrix\right.\tag3.3
$$
where by definition $e_\infty:=0$. It will turn out that these particular special definitions for the vertices
in $B_1$ and $B_2$ are just the right ones for obtaining the desired decomposition of $\overline{HMD}_{n}$.

Regard (3.3) as vectors in ${\bold U}={\bold U}_{\overline{HMD}_{n}}$.
A case analysis analogous to that in the proof of Lemma 2.4 shows that the action of 
$F=F_{\overline{HMD}_{n}}$ on the vectors (3.3) is given by
$$
F(f_v)=f_{\No(v)}+f_{\We(v)}+f_{\So(v)}+f_{\Ea(v)},\tag3.4
$$
where $\No(v)$, $\We(v)$, $\So(v)$, and $\Ea(v)$ are defined in analogy to ${\NE(v)}$---via 
rays from $v$ in the directions of the four cardinal points---and $f_\infty:=0$.

The set $V$ of vertices of $\overline{HMD}_{n}$ naturally splits into $n$ levels---denote them, from
bottom to top, by $L_1,\dotsc,L_n$. The same set can also be regarded as consisting of $2n$ columns
of vertices---denote them, from left to right, by $C_1,\dotsc,C_{2n}$. We define subsets 
$S_1,\dotsc,S_n\subset V$ as follows. 

Color the vertices of $V$ black and white in chessboard fashion so that $C_1$ is black; let $V_B$ be
the set of black vertices. 
For $1\leq k\leq n$, set
$$
\align
&\ \ \ 
S_k:=\left(V_B\cap(C_1\cup C_2\cup\cdots\cup C_k\cup C_{2n-k+1}\cup C_{2n-k+2}\cup\cdots\cup C_{2n})\right)
\\
&\ \ \ \ \ \ \ \ \ \ \ \ \ \ \ \ \ \ \ \ \
     \cup
     \left((C_k\cup C_{k+1}\cup\cdots\cup C_{2n-k+1})\cap(L_k\cup L_{k-2}\cup L_{k-4}\cup\cdots )\right)\tag3.5
\endalign
$$
(for $n=6$, these are pictured in Figure 3.10).

\topinsert
\threeline{\mypic{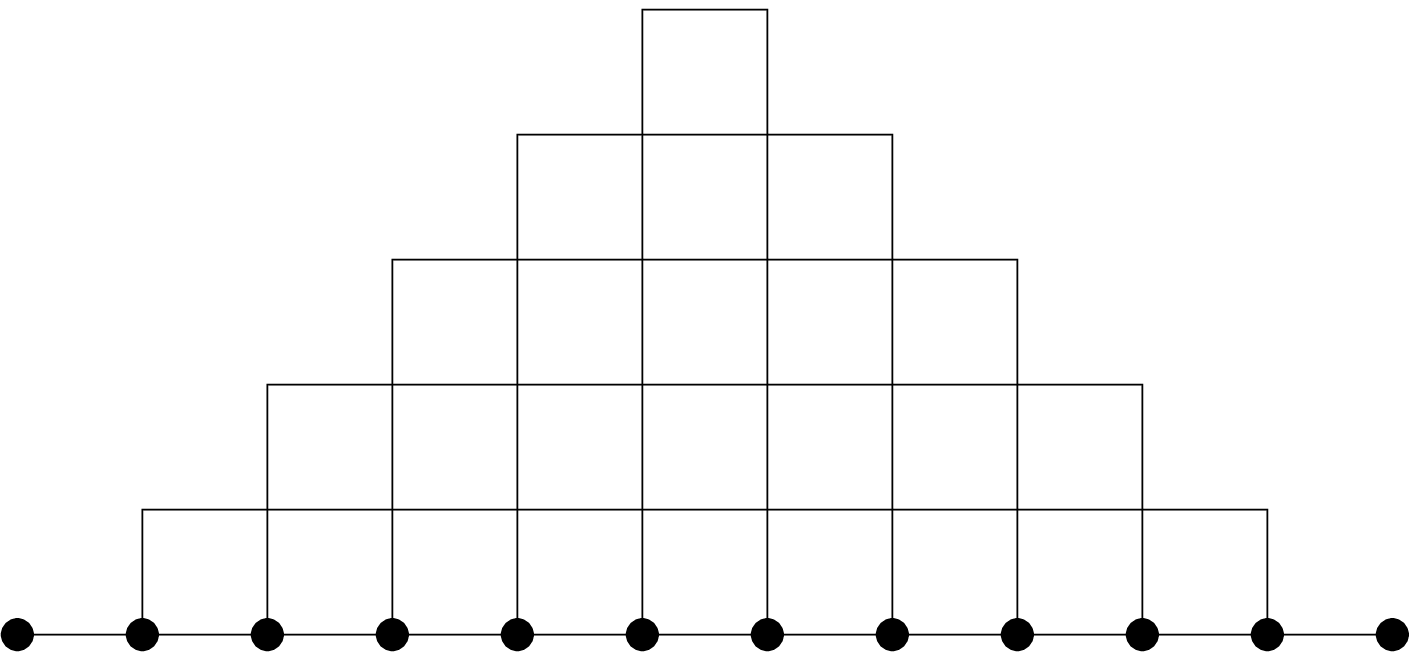}}{\mypic{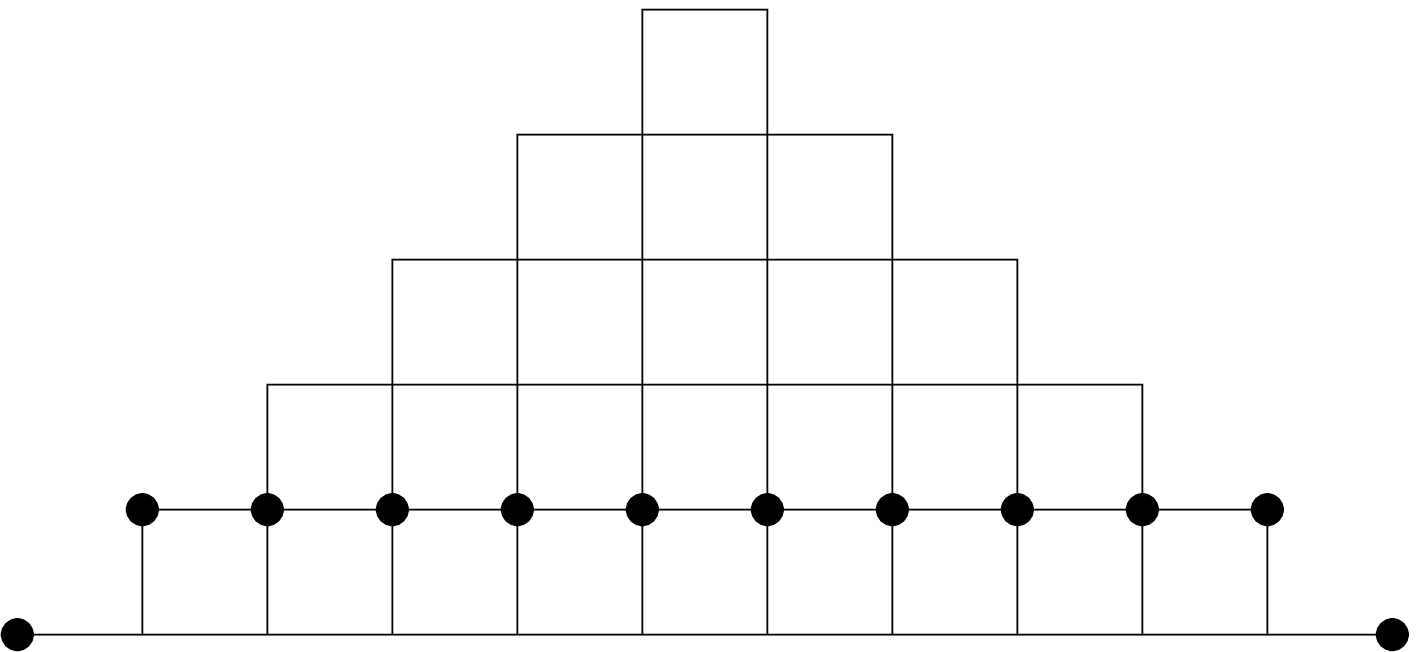}}{\mypic{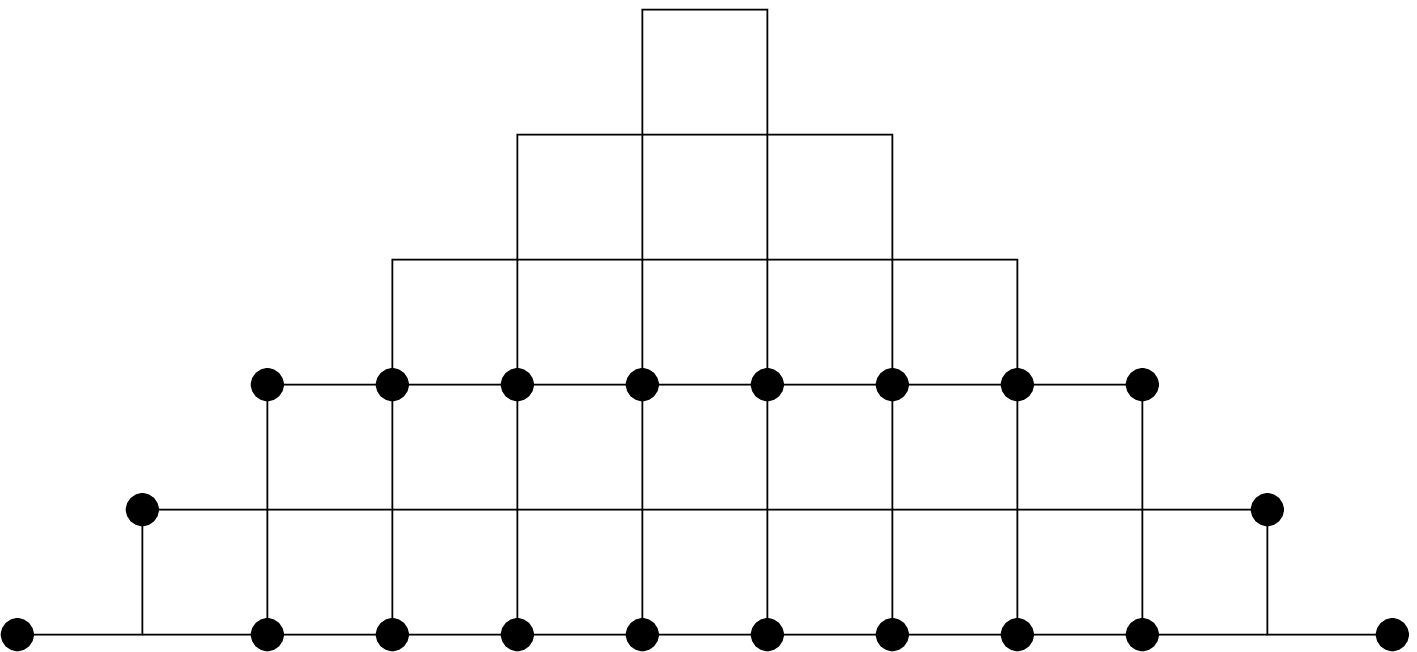}}

\bigskip
\threeline{\mypic{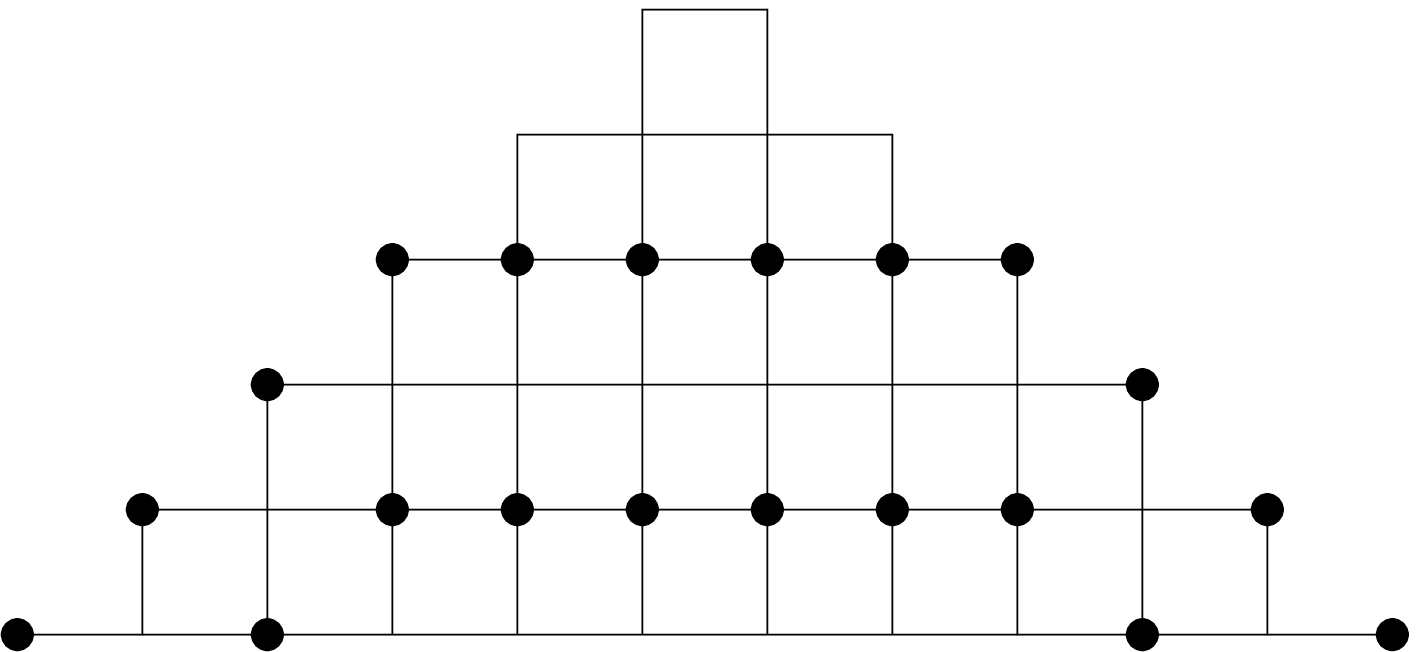}}{\mypic{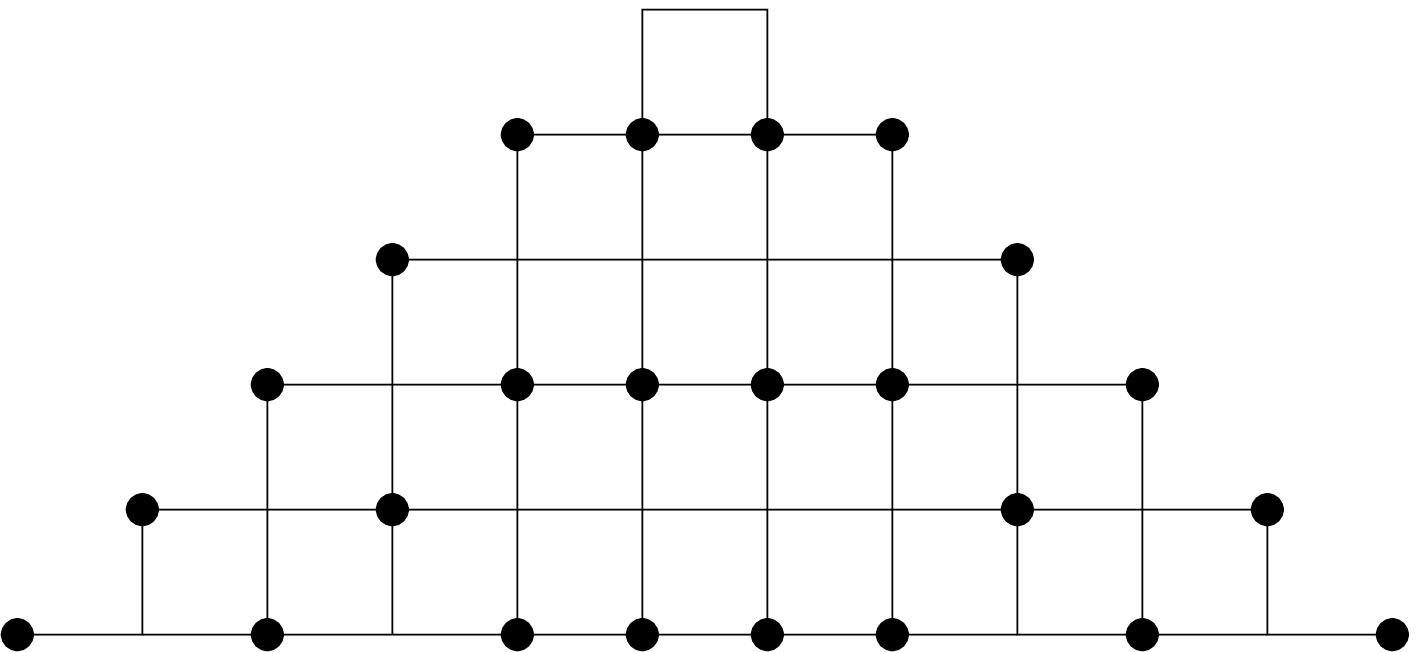}}{\mypic{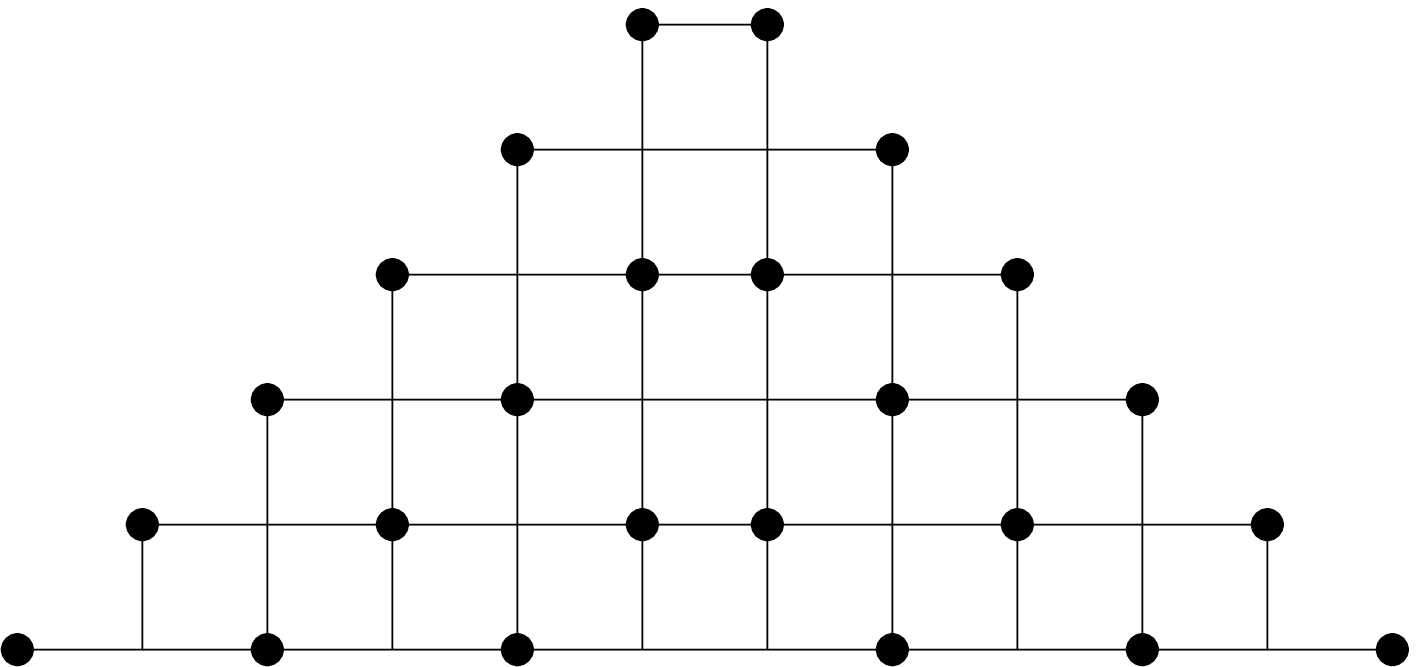}}
\bigskip
\centerline{{\smc Figure~3.10.} {\rm The sets $S_k$ for $n=6$.}}
\endinsert

Let $\wt$ be the weight function on $V$ that assigns 1 to the vertices in $L_1$, and $(-1)^{k-1}\cdot2$ to
the vertices in $L_k$, if $k=2,\dotsc,n$. Let $\wt'(v):=\wt(v)$ if $v$ is black, and $\wt'(v):=-\wt(v)$ 
if $v$ is white.

Define, for $k=1,\dotsc,n$,
$$
g_k:=(-1)^{k-1}\sum_{v\in S_k}\wt(v) e_v\tag3.6
$$ 
and
$$
g'_k:=\sum_{v\in S_k}\wt'(v) e_v.\tag3.7
$$ 

An argument similar to the one that proved Lemma 2.5 shows that the action of $F$ on these vectors
is given by
$$
F(g_1)=2g_{1}+g_2,\tag3.8
$$
$$
F(g_k)=g_{k-1}+2g_k+g_{k+1},\ \ \ 2\leq k\leq n-1,\tag3.9
$$
$$
F(g_n)=g_{n-1}+g_n,\tag3.10
$$
and
$$
F(g'_1)=-2g'_{1}+g'_2,\tag3.11
$$
$$
F(g'_k)=g'_{k-1}-2g'_k+g'_{k+1},\ \ \ 2\leq k\leq n-1,\tag3.12
$$
$$
F(g'_n)=g'_{n-1}-g'_n.\tag3.13
$$

Furthermore, by Lemma 3.2, the span of the vectors (3.3), (3.6) and (3.7) contains
$u_k=\sum_{v\in L_k} e_v$ and $u'_k=\sum_{v\in L_k} c_v e_v$, $k=1,\dotsc,n$, where $c_v$ equals 
1 or $-1$ according as $v$ is black or white. 
Then a simple inductive argument shows that each
$e_v$, $v\in V$ is contained in the span of the union of the vectors (3.3) with the $u_k$'s and the
$u'_k$'s. This implies that the union of the vectors (3.3), (3.6) and (3.7) forms a basis of ${\bold U}$.

However, by (3.4) and (3.8)--(3.13), the matrix of $F$ in this basis is a block diagonal matrix consisting
of three blocks. One of them, corresponding to the rows and columns indexed by the vectors (3.3), is
by (3.4) the same as the adjacency matrix of $HMD_{n-1}$. The other two, by~(3.8)--(3.13), are the
same as the adjacency matrices of $R_n$ and $R'_n$, respectively. This implies~(3.1).~\epf

\bigskip

\proclaim{Lemma 3.2} The span of the vectors $\{f_v:v\in V_1\}$, $\{g_k:k=1,\dotsc,n\}$, and $\{g'_k:k=1,\dotsc,n\}$
given by $(3.3)$, $(3.6)$, and $(3.7)$ contains
$u_k=\sum_{v\in L_k} e_v$ and $u'_k=\sum_{v\in L_k} c_v e_v$, $k=1,\dotsc,n$, where $c_v$ equals 
$1$ or $-1$ according as $v$ is black or white. 
\endproclaim

\topinsert
\centerline{\mypic{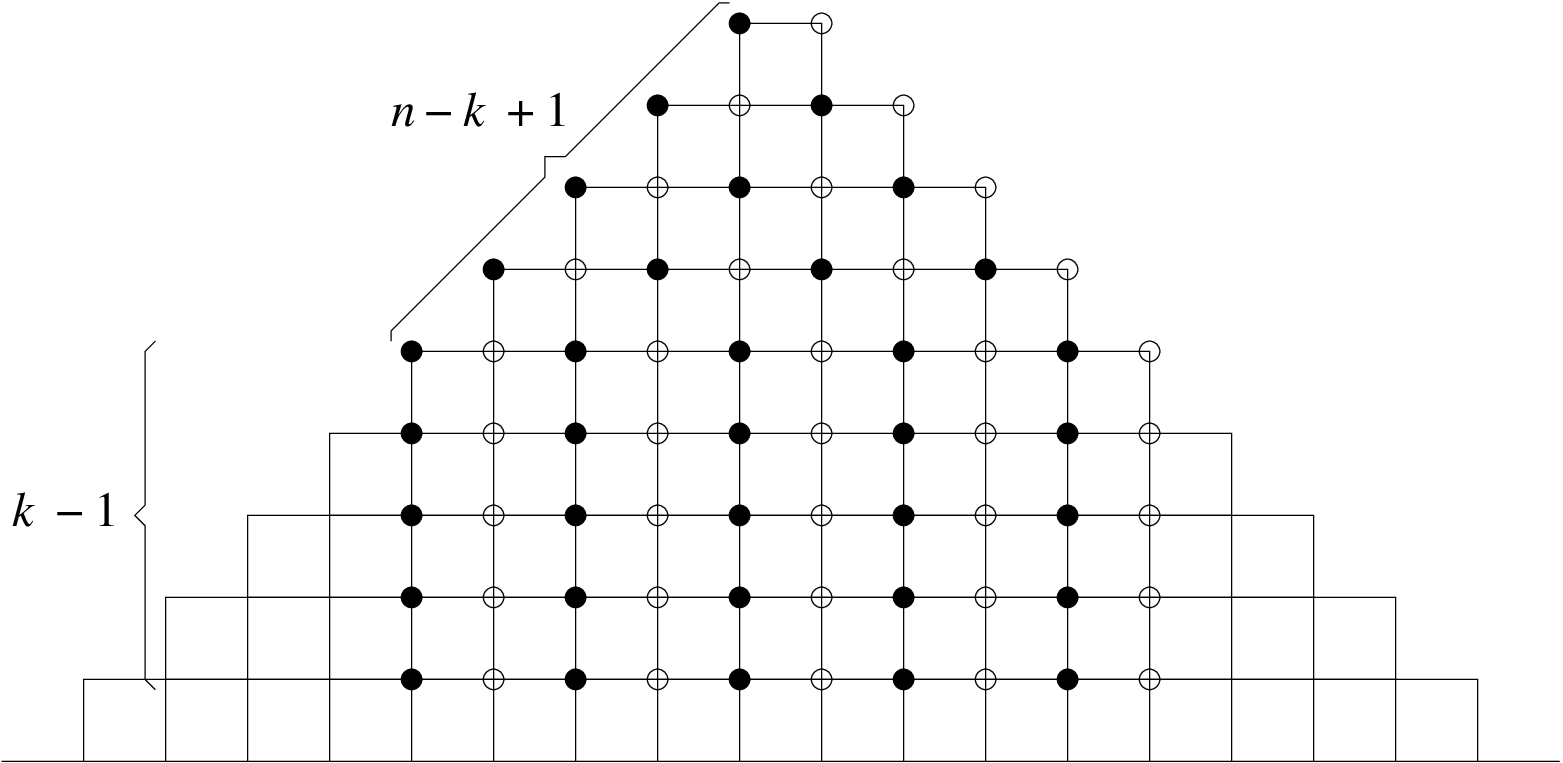}}
\medskip
\centerline{{\smc Figure~3.11.} {\rm The support of $\alpha_{n,k}$ (circled vertices) and $\beta_{n,k}$ 
(dotted vertices).}}
\endinsert

\pf
Suppose $n$ is even. For $2\leq k\leq n$, let $\alpha_{n,k}$ be the function on $V$ that is zero on the vertices not 
circled in Figure 3.11,
and whose value at a circled vertex is given by the corresponding entry in the following array (for briefness we
denote $a=n-k$):
$$
\spreadmatrixlines{4\jot}
\matrix 
\ && && && && && &&{\scriptstyle k-1}\ms&& && && && && && &&
\\
 && && && && &&{\scriptstyle -(k-1)}\ms&& &&{\scriptstyle -2(k-1)}\ms&& && && && && && 
\\
 && && && &&
{\scriptstyle k-1}\ms&& &&{\scriptstyle 2(k-1)}\ms&& &&{\scriptstyle 3(k-1)}\ms&& && && && &&
\\
 && && && && \vdots\ms && && && && \vdots\ms && && && && &&  
\\
 && && {\scriptstyle k-1}\ms&& &&{\scriptstyle 2(k-1)}\ms&&  &&{\scriptstyle 3(k-1)}\ms&& &&{\scriptstyle \cdots}\ms&& && 
{\scriptstyle (a-1)(k-1)}\ms && && &&
\\
 && {\scriptstyle -(k-1)}\ms&& &&{\scriptstyle -2(k-1)}\ms&&  &&{\scriptstyle -3(k-1)}\ms&& &&{\scriptstyle \cdots}\ms&& && 
{\scriptstyle -(a-1)(k-1)}\ms && &&{\scriptstyle -a(k-1)}\ms && && 
\\
{\scriptstyle k-1}\ms&& &&{\scriptstyle 2(k-1)}\ms&&  &&{\scriptstyle 3(k-1)}\ms&& &&{\scriptstyle \cdots}\ms&& && 
{\scriptstyle (a-1)(k-1)}\ms && &&{\scriptstyle a(k-1)}\ms && &&{\scriptstyle (a+1)(k-1)}\ms &&  
\\
{\scriptstyle k-2}\ms&& &&{\scriptstyle 2(k-2)}\ms&&  &&{\scriptstyle 3(k-2)}\ms&& &&
{\scriptstyle \cdots}\ms&& &&{\scriptstyle (a-1)(k-2)}\ms&& &&{\scriptstyle a(k-2)}\ms&& &&{\scriptstyle (a+1)(k-2)}\ms&&
\\
 && && &&{\scriptstyle \vdots}\ms&& && && && && && &&{\scriptstyle \vdots}\ms&& && && 
\\
{\scriptstyle 1\cdot2}\ms&& &&{\scriptstyle 2\cdot2}\ms&& &&{\scriptstyle 3\cdot2}\ms&& &&
{\scriptstyle \cdots}\ms&& &&{\scriptstyle (a-1)\cdot2}\ms&& &&{\scriptstyle a\cdot2}\ms&& &&{\scriptstyle (a+1)\cdot2}\ms&&
\\
{\scriptstyle 1}\ms&& &&{\scriptstyle 2}\ms&& &&{\scriptstyle 3}\ms&& &&
{\scriptstyle \cdots}\ms&& &&{\scriptstyle a-1}\ms&& &&{\scriptstyle a}\ms&& &&{\scriptstyle a+1}\ms&&
\endmatrix
$$
(read the pattern above from bottom to top; for $k$ even the pattern ends on top as indicated; for $k$ odd the 
alternation of signs along the rows of the top triangular portion of the array causes the entries in the top three rows above 
to be the negatives of the shown ones).

Let $\beta_{n,k}$ be the function on $V$ that is zero on the vertices not dotted in Figure 3.11,
and defined at each dotted vertex $v$ by $\beta_{n,k}(v)=-\alpha_{n,k}(v')$, where $v'$ is the reflection of $v$ across the
vertical symmetry axis of $HMD_n$ (note that the supports of $\alpha_{n,k}$ and $\beta_{n,k}$ are disjoint).

It is straightforward to check that for $k=2,\dotsc,n$ one has
$$
\spreadlines{3\jot}
\align
&
\sum_{v\in V_1}\left(\alpha_{n,k}(v)+\beta_{n,k}(v)\right)f_v   
+\left[-(n-k+1)g_{k-2}-g_{k-1}+(n-k+2)g_k\right]=(2n+1)\sum_{v\in L_k}e_v.
\\
&
\endalign
$$
Furthermore, if $c_v$ is plus or minus one according as $v$ is black or white, one similarly checks that
$$
\spreadlines{3\jot}
\align
&
\sum_{v\in V_1}c_v\left(\alpha_{n,k}(v)+\beta_{n,k}(v)\right)f_v   
+(-1)^k\left[(n-k+1)g'_{k-2}+g'_{k-1}-(n-k+2)g'_k\right]
\\
&\ \ \ \ \ \ \ \ \ \ \ \ \ \ \ \ \ \ \ \ \ \ \ \ \ \ \ \ \ \ \ \ \ \ \ \ \ \ \ \ \ \ \ \ \ \ \ \ \ \ 
=(2n+1)\sum_{v\in L_k}c_ve_v,
\endalign
$$
for $k=2,\dotsc,n$. Since by definition $u_1=g_1$ and $u'_1=g'_1$, the above equalities prove the statement for $n$ even.

The above equalities hold without change also for odd $n$, provided we define $\alpha_{n,k}$ and $\beta_{n,k}$ to be the
negatives of their values above. \epf

\proclaim{Theorem 3.3} We have
$$
OD_n\sim HOD_{n-1}\ \dot{\cup}\  HOD_{n-1}\ \dot{\cup}\  
P_1\ \dot{\cup}\  Q_n\ \dot{\cup}\  Q'_n.\tag3.14
$$

\endproclaim

\pf We proceed in a way analogous to the proof of Theorem 3.1. Lemma 1.1 implies
$$
OD_n\sim HOD_{n-1}\ \dot{\cup}\ \ \overline{HOD}_n,\tag3.15
$$
where $\overline{HOD}_n$ is the weighted directed graph obtained from $HOD_n$ by marking its 
bottommost $2n+1$ vertices, replacing each edge of it between two unmarked or two marked vertices
by a pair of opposite arcs of weights 1, and each edge between an unmarked vertex $v$ and a marked
vertex $w$ by an arc $(v,w)$ of weight 1 and an arc $(w,v)$ of weight 2 (Figure 3.12 shows 
$\overline{HOD}_4$).

\topinsert
\centerline{\mypic{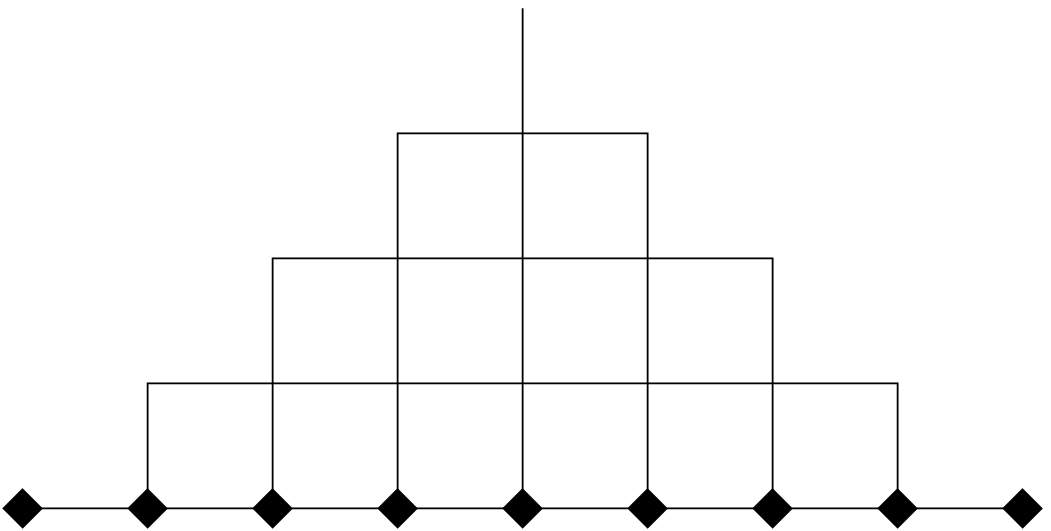}}
\centerline{{\smc Figure~3.12.} {\rm $\overline{HOD}_4$.}}
\endinsert

Let $V$ be the set of vertices of $\overline{HOD}_n$, and let $V_1$ be the set of its unmarked vertices.
For each $v\in V$, let $e_v$ be an indeterminate. For $v\in V_1$, define $f_v$ by (3.3), the same formulas
we used in the proof of Theorem 3.1---except for the case when $v$ is the topmost vertex in $V_1$, when we define
$f_v:=-e_{\SW(v)}+e_{\SE(v)}$. Then, as an argument similar to the one in the proof of Theorem 3.1
readily checks, (3.4) holds for all $v\in V_1$.

The subsets of vertices and the weights on $V$ that will provide us with the vectors we need to complete
$\{f_v:v\in V_1\}$ to a basis of ${\bold U}={\bold U}_{\overline{HOD}_n}$ are now defined slightly
differently than in the proof of Theorem 3.1. Define the $g_k$'s and $g'_k$'s to be the linear
combinations of the $e_v$'s indicated by the patterns shown in Figure 3.13(a) and (b), respectively (these
figures correspond to $n=4$).

\topinsert
\twoline{\mypic{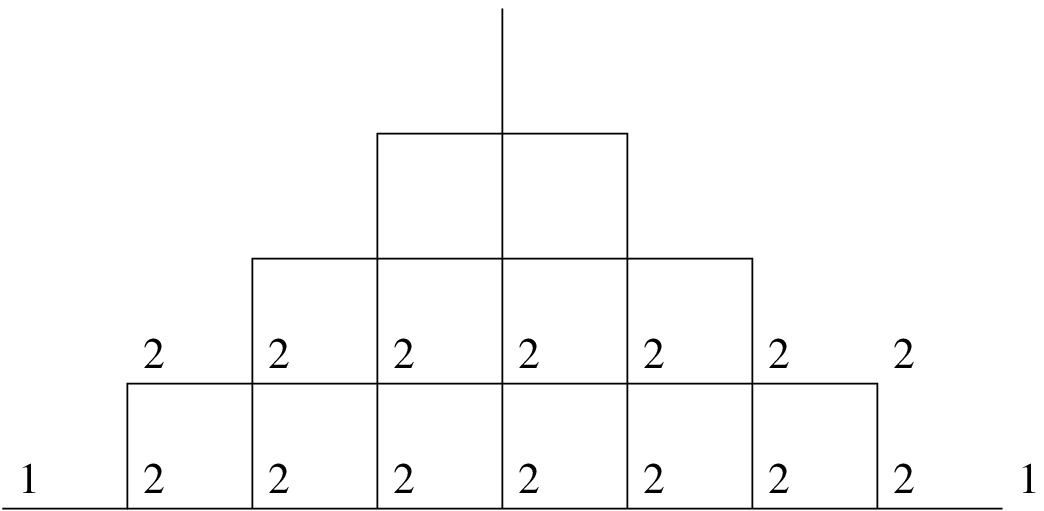}}{\mypic{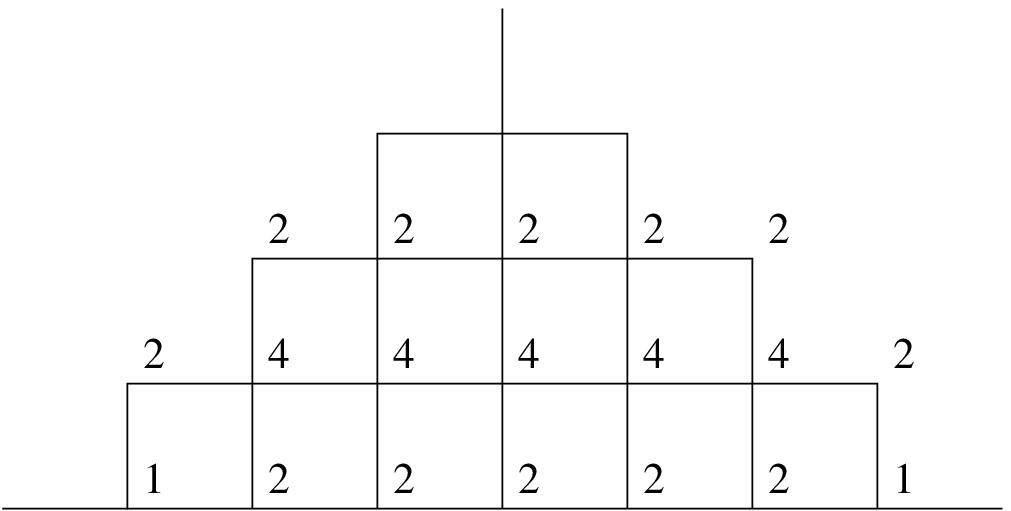}}
\bigskip
\twoline{\mypic{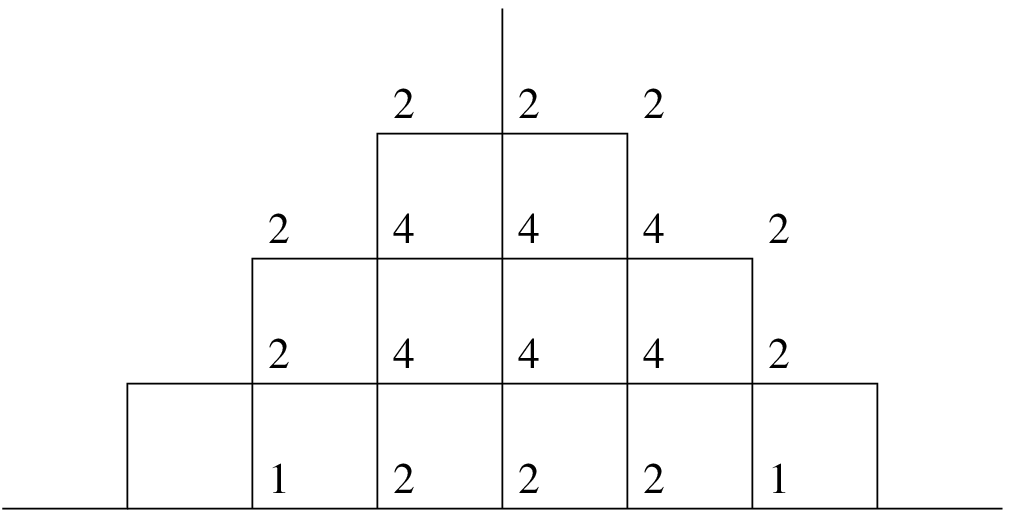}}{\mypic{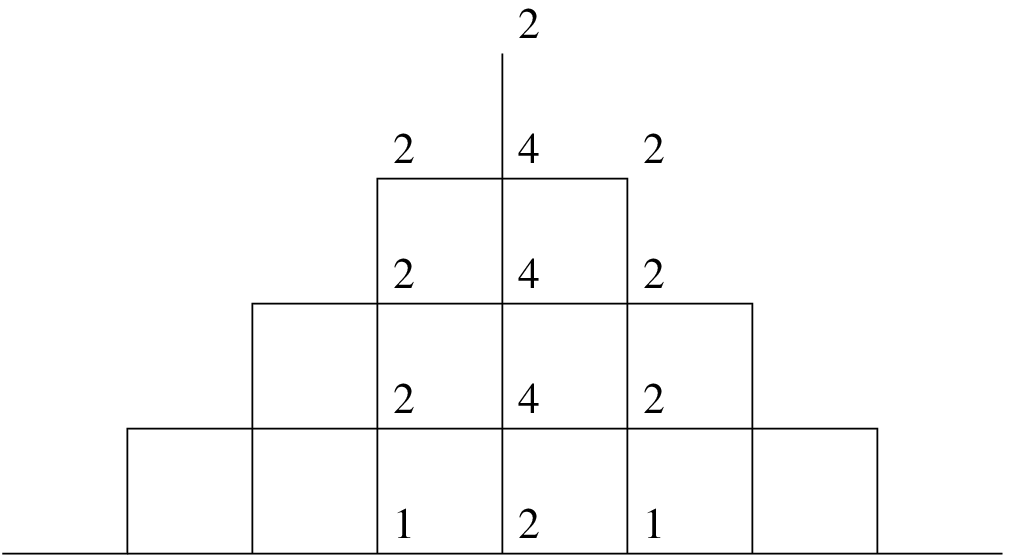}}
\bigskip
\centerline{\smc Figure 3.13{\rm (a). The $g_k$'s for $n=4$.}}

\bigskip
\twoline{\mypic{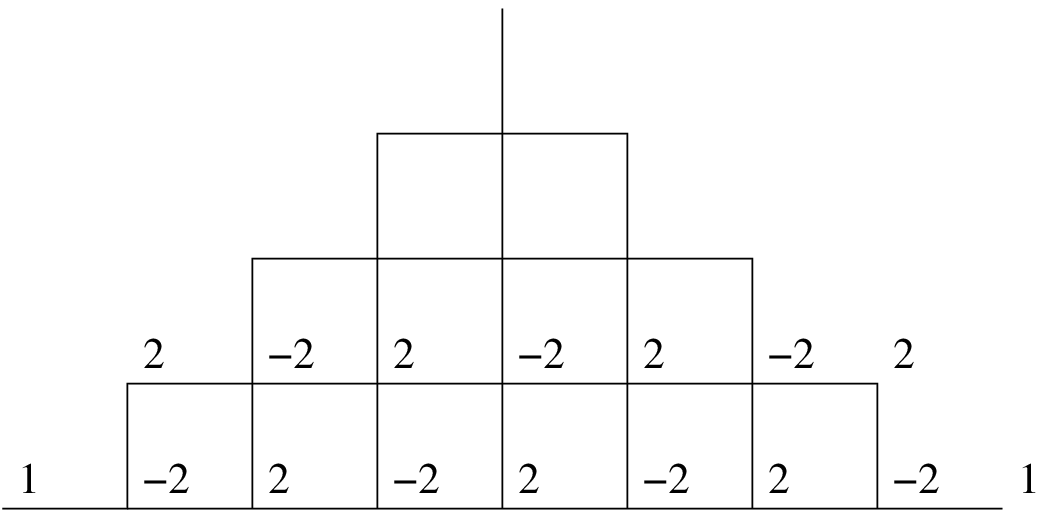}}{\mypic{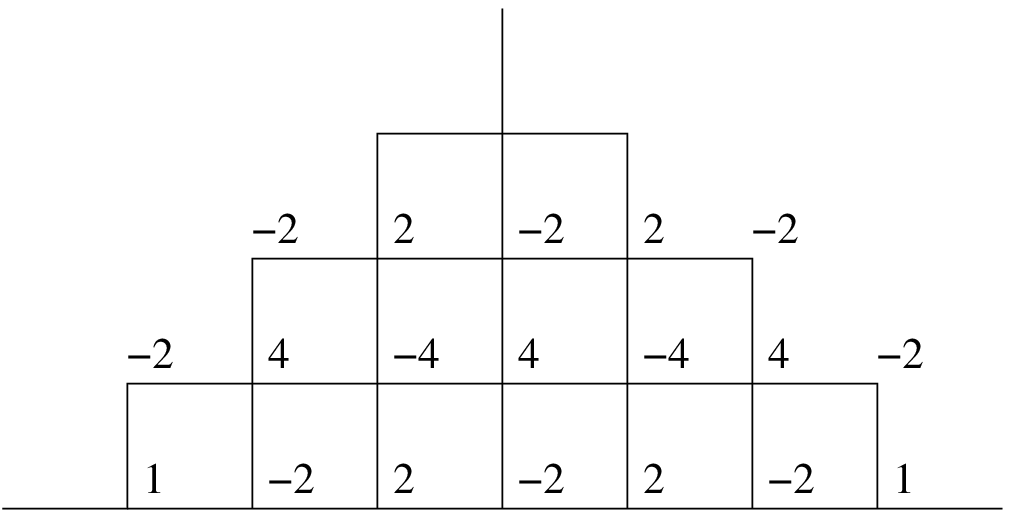}}
\bigskip
\twoline{\mypic{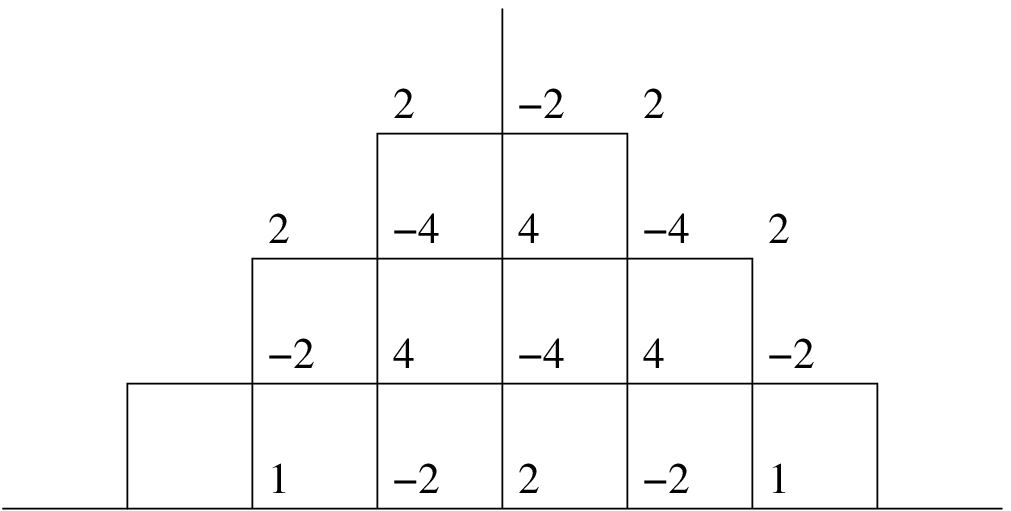}}{\mypic{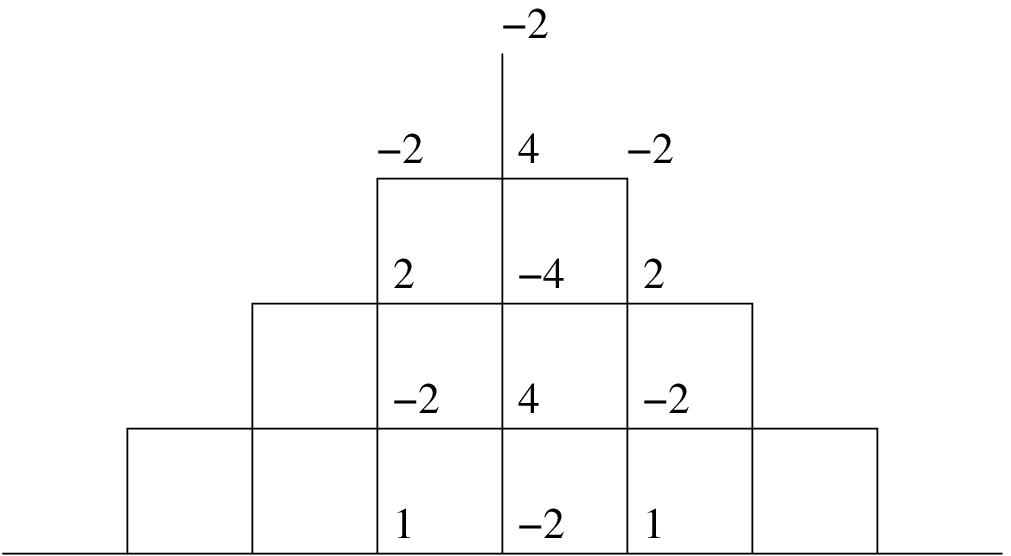}}
\bigskip
\centerline{\smc Figure 3.13{\rm (b). The $g'_k$'s for $n=4$.}}

\endinsert



An argument similar to the one given in Lemma 2.5 shows that the action of $F=F_{\overline{HOD}_n}$ 
on these vectors is given by
$$
F(g_1)=2g_{1}+g_2,\tag3.16
$$
$$
F(g_k)=g_{k-1}+2g_k+g_{k+1},\ \ \ 2\leq k\leq n-1,\tag3.17
$$
$$
F(g_n)=g_{n-1}+2g_n,\tag3.18
$$
and
$$
F(g'_1)=-2g'_{1}+g'_2,\tag3.19
$$
$$
F(g'_k)=g'_{k-1}-2g'_k+g'_{k+1},\ \ \ 2\leq k\leq n-1,\tag3.20
$$
$$
F(g'_n)=g'_{n-1}-2g'_n.\tag3.21
$$
Finally, the vector $h$ that equals the alternating sum of $e_v$'s over the set of the $n$ vertices
of $\overline{HOD}_n$ on the northwestern side of the convex hull of its vertex set (the topmost $e_v$ having
coefficient +1)
is readily seen to satisfy $F(h)=0$.


An inductive argument similar to the one used in the proof of Lemma 2.6 shows that
for all $v\in V$, $e_v$ is in the span of the vectors $\{f_v:v\in V_1\}$, $\{\tilde{u_k}:k=1,\dotsc,n\}$, $\{\tilde{u}'_k:k=1,
\dotsc,n\}$ and $h$ (the only part that requires separate justification is the checking of the base case; this
is provided by Lemma 3.5). Thus, by Lemma 3.4, the span of the vectors $\{f_v:v\in V_1\}$, $\{g_k:k=1,\dotsc,n\}$, $\{g'_k: 
k=1,\dotsc,n\}$, and $h$ contains $e_v$ for all $v\in V$, and hence these vectors form a basis of $\bold U$.

By the above formulas describing the action of $F$ on them it follows that the matrix of $F$ in this 
basis is a block-diagonal matrix
consisting of four blocks: one the same as the adjacency matrix of $HOD_{n-1}$, by the fact that (3.4)
holds for $V_1$; the next two equal to the adjacency matrices of the path-like graphs $Q_n$ and $Q'_n$,
respectively, by (3.16)--(3.21); and the last equal to the $1\times1$ block consisting of a single 0, 
i.e., the adjacency matrix of the path $P_1$.
Together with (3.15) this implies (3.14). \epf


\proclaim{Lemma 3.4} The span of the vectors $\{f_v:v\in V_1\}$, $\{g_k:k=1,\dotsc,n\}$, $\{g'_k:k=1,\dotsc,n\}$, and $h$ 
in the proof of Theorem 3.3 contains the vectors
$\tilde{u}_k$ and $\tilde{u}'_k$ defined by the patterns in Figure 3.14, for $k=1,\dotsc,n$. 

\endproclaim

\topinsert
\twoline{\mypic{3-12a1.eps}}{\mypic{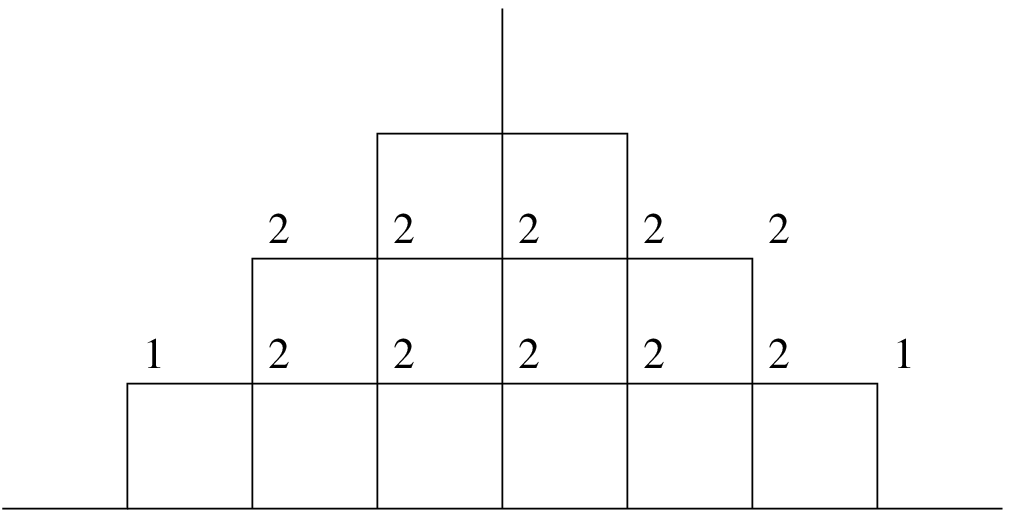}}
\bigskip
\twoline{\mypic{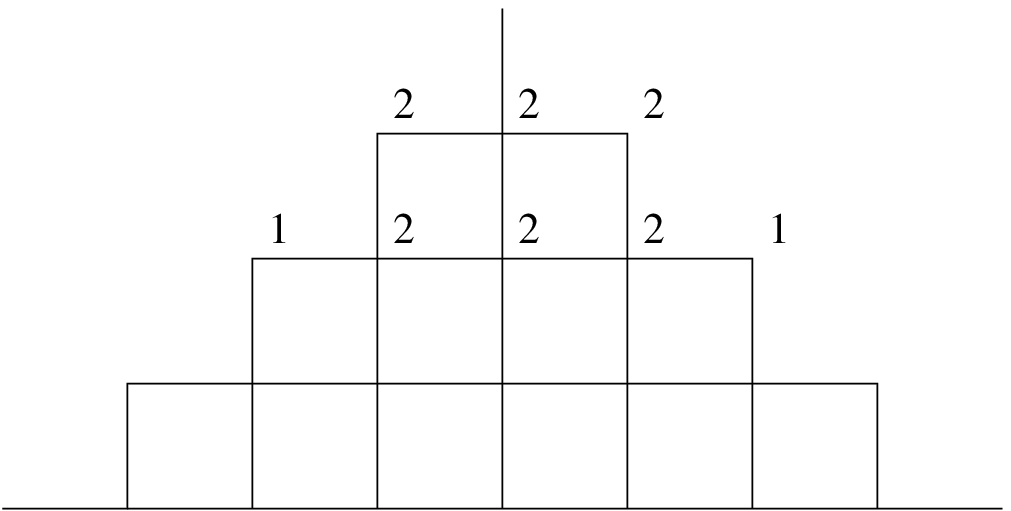}}{\mypic{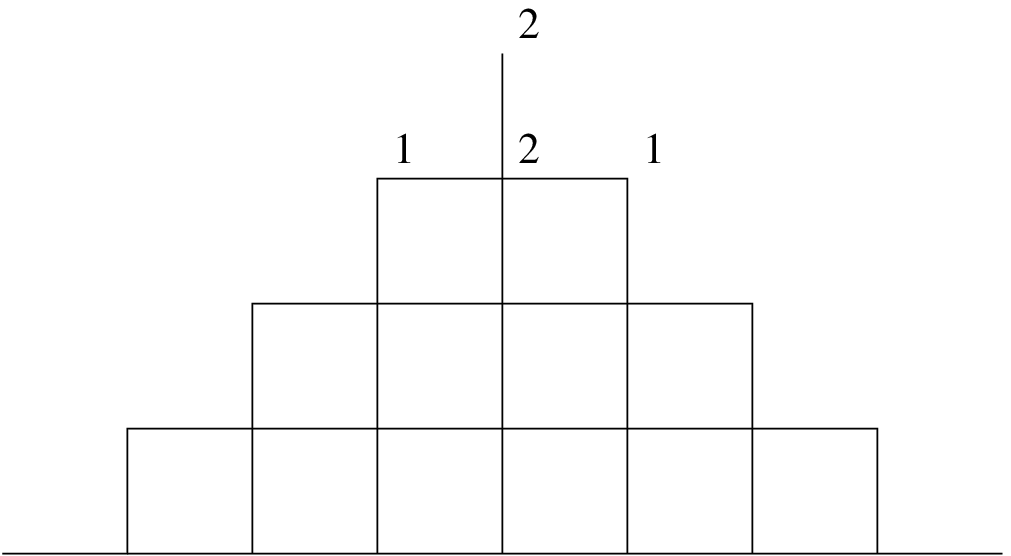}}
\bigskip
\centerline{\smc Figure 3.14{\rm (a). The $\tilde{u}_k$'s for $n=4$.}}

\bigskip
\twoline{\mypic{3-12b1.eps}}{\mypic{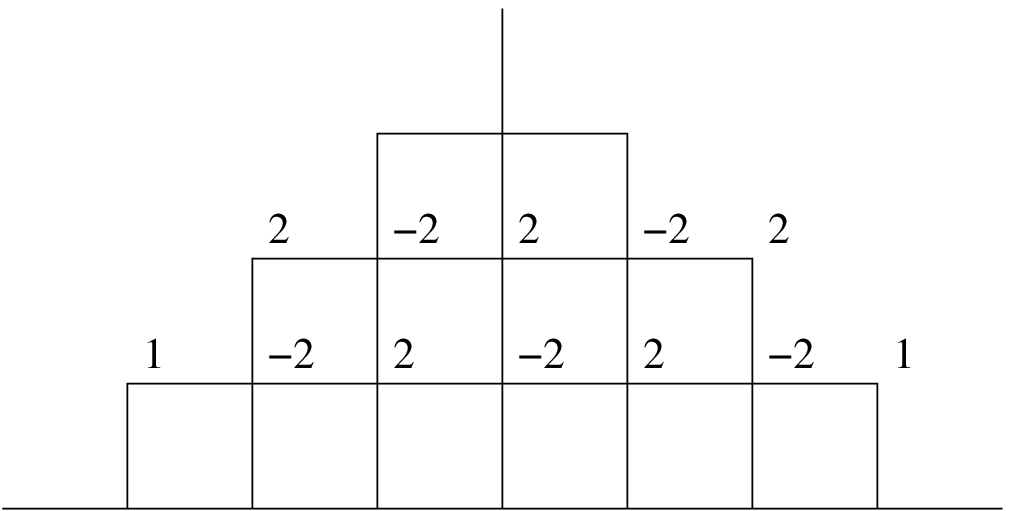}}
\bigskip
\twoline{\mypic{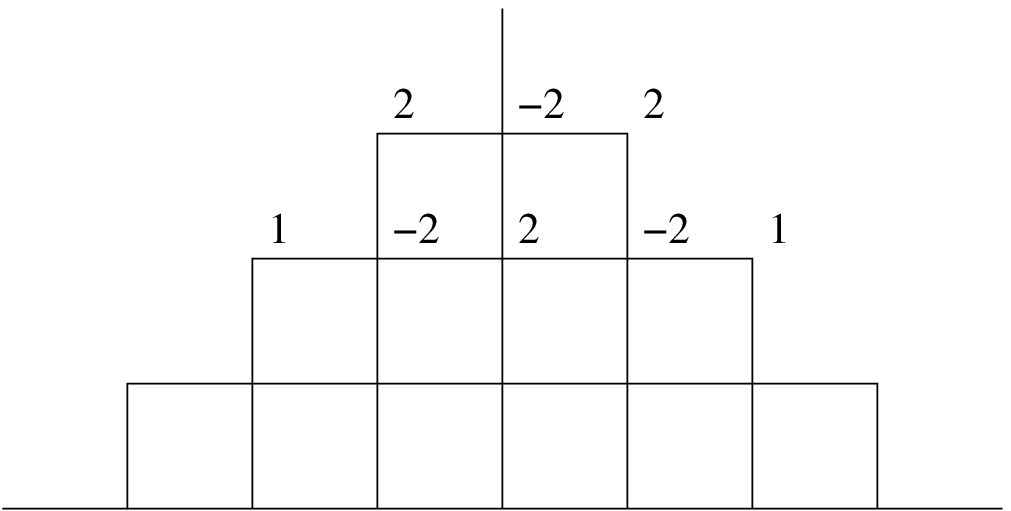}}{\mypic{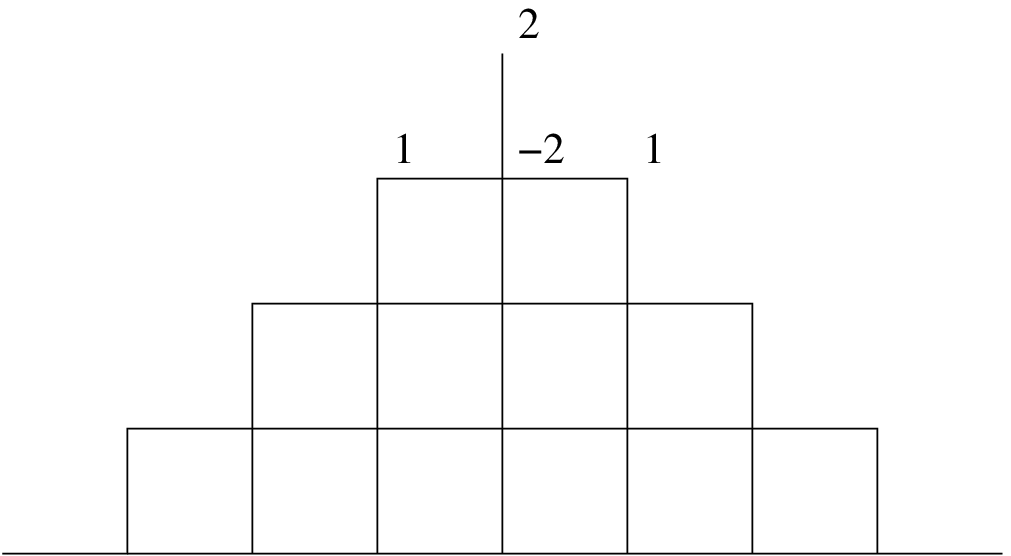}}
\bigskip
\centerline{\smc Figure 3.14{\rm (b). The $\tilde{u}'_k$'s for $n=4$.}}

\endinsert

\pf Let $\wt_n$ be the function on the vertices of of $HOD_n$ given by the pattern

$$
\spreadmatrixlines{4\jot}
\ \ \ \ 
\matrix 
&& && && && && && &&\mms{\scriptstyle 0}&& && && && && && &&
\\
&& && && && && &&\mms{\scriptstyle-(n-1)}&&\mms{\scriptstyle0}&&\mms{\scriptstyle n-1}&& && && && && &&
\\
&& && && && &&\mms{\scriptstyle -2(n-2)}&&\mms{\scriptstyle -(n-2)}&&\mms{\scriptstyle 0}&&\mms{\scriptstyle n-2}&&
\mms{\scriptstyle 2(n-2)}&& && && && &&
\\
&& && && && &&\mms{\scriptstyle \vdots}&&\mms{\scriptstyle \vdots}&&\mms{\scriptstyle \vdots}&&\mms{\scriptstyle \vdots}&&
\mms{\scriptstyle \vdots}&& && && && &&
\\
&& && &&\mms{\scriptstyle -3(n-3)}&&\mms{\scriptstyle \cdots}&&\mms{\scriptstyle -6}&&\mms{\scriptstyle -3}&&\mms{\scriptstyle 0}&&
\mms{\scriptstyle 3}&&\mms{\scriptstyle 6}&&\mms{\scriptstyle \cdots}&&\mms{\scriptstyle 3(n-3)}&& && &&
\\
&& &&\mms{\scriptstyle -2(n-2)}&&\mms{\scriptstyle \cdots}&&\mms{\scriptstyle -6}&&\mms{\scriptstyle -4}&&\mms{\scriptstyle -2}&&
\mms{\scriptstyle 0}&&\mms{\scriptstyle 2}&&\mms{\scriptstyle 4}&&\mms{\scriptstyle 6}&&\mms{\scriptstyle \cdots}&&
\mms{\scriptstyle 2(n-2)}&& &&
\\
&&\mms{\scriptstyle -(n-1)}&&\mms{\scriptstyle \cdots}&&\mms{\scriptstyle -4}&&\mms{\scriptstyle -3}&&\mms{\scriptstyle -2}&&
\mms{\scriptstyle -1}&&\mms{\scriptstyle 0}&&\mms{\scriptstyle 1}&&\mms{\scriptstyle 2}&&\mms{\scriptstyle 3}&&\mms{\scriptstyle 4}&&
\mms{\scriptstyle \cdots}&&\mms{\scriptstyle n-1}&& 
\\
\mms{\scriptstyle 0}&&\mms{\scriptstyle \cdots}&&\mms{\scriptstyle 0}&&\mms{\scriptstyle 0}&&\mms{\scriptstyle 0}&&
\mms{\scriptstyle 0}&&\mms{\scriptstyle 0}&&\mms{\scriptstyle 0}&&\mms{\scriptstyle 0}&&\mms{\scriptstyle 0}&&
\mms{\scriptstyle 0}&&\mms{\scriptstyle 0}&&\mms{\scriptstyle 0}&&
\mms{\scriptstyle \cdots}&&\mms{\scriptstyle 0}
\endmatrix
$$

\topinsert
\centerline{\mypic{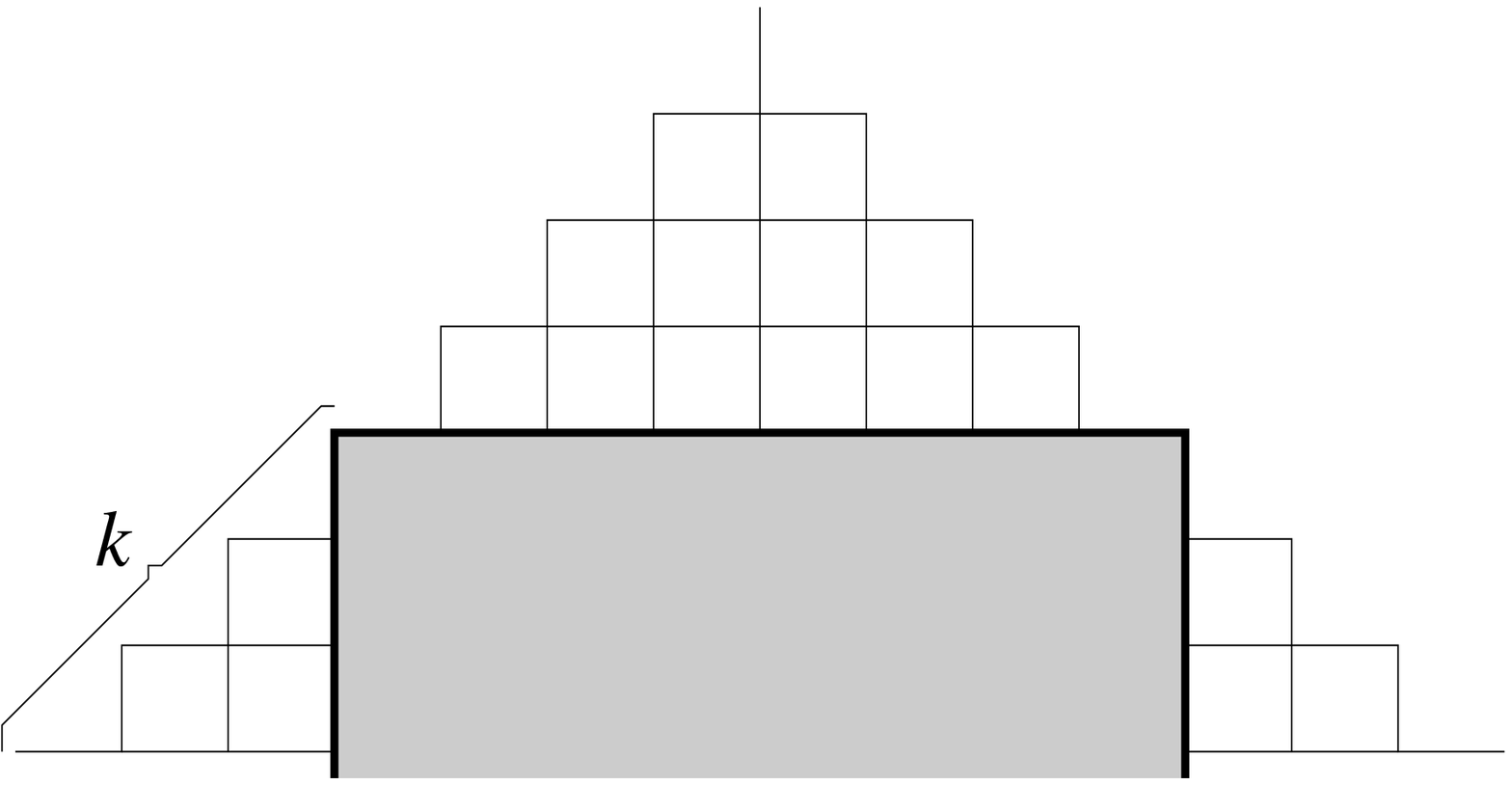}}
\centerline{{\smc Figure~3.15.} The region $C_4$ for $HOD_7$.}
\endinsert

Let $C_k$ be the closed shaded region illustrated in Figure 3.15. It is straightforward to check that
$$
\sum_{v\in V_1\cap C_k} \wt_n(v)f_v+\left[-(n-k+1)g_{k-1}+(n-k+2)g_k\right]=(n+1)\tilde{u}_k,
$$ 
for all $k=1,\dotsc,n$.

Similarly, if $\wt'_n$ is the weight function on the vertices of $HOD_n$ that equals $\wt_n$ on the black vertices and $-\!\wt_n$
on the white vertices, one checks that
$$
\sum_{v\in V_1\cap C_k} \wt'_n(v)f_v+(-1)^{k-1}\left[(n-k+1)g'_{k-1}+(n-k+2)g'_k\right]=(n+1)\tilde{u}'_k,
$$ 
for $k=1,\dotsc,n$. \epf

\proclaim{Lemma 3.5} Let $t$ be the topmost vertex of $HOD_n$. Then $e_t$ is in the span of the vectors $\{f_v:v\in V_1\}$, 
$\{\tilde{u}_k:k=1,\dotsc,n\}$, $\{\tilde{u}'_k:k=1,\dotsc,n\}$, and $h$.

\endproclaim  

\pf Let $\alpha_{n}$ be the weight function on the subset $V_1$ of the vertices of $\overline{HOD}_n$ obtained by weighting the 
$i$th topmost row, from left to right, by the successive coordinates of the vector
$$
\frac{(-1)^i(n-i+1)}{2i}(2i-1,0,2i-3,0,\dotsc,3,0,1),
$$
for $i=1,\dotsc,n$ (recall that by definition $HOD_n$ has $n+1$ rows).

Let 
$$
d_k=\frac{(-1)^{n+k}(2k(n-k+2)-n-1)}{4(n-k+1)(n-k+2)},\ \ \ k=1,\dotsc,n.
$$
Then it is straightforward to check that
$$
\sum_{v\in V_1}\alpha_{n}(v)f_v + \sum_{k=1}^n d_k\left(\tilde{u}_k+\tilde{u}'_k\right) + h =(n+1)e_t.
$$
This proves the claim. \epf

\proclaim{Corollary 3.6} The number of spanning trees of the halved Aztec diamond $HMD_n$ is
given by
$$
\te(HMD_n)=\prod_{{\scriptstyle 1\leq j<k\leq 2n-1}\atop {\scriptstyle j+k\leq 2n-1}}
\left(4-4\cos\frac{j\pi}{2n}\cos\frac{k\pi}{2n}\right).\tag3.22
$$
\endproclaim

\pf The graph obtained from the planar dual of $HMD_n$ by removing the vertex corresponding to the
infinite face is readily seen to be isomorphic to $HOD_{n-2}$ (see Figures 3.3 and 3.4). 
Thus, Theorem 2.2 implies
$$
\te(HMD_n)=\Pe(HOD_{n-2};4).\tag3.23
$$
By \cite{\Kn} we have
$$
\Pe(OD_n;x)=x\prod_{j=1}^n\left(x-4\cos^2\frac{j\pi}{2n+2}\right)
\prod_{1\leq j<k\leq 2n+1}\left(x-4\cos\frac{j\pi}{2n+2}\cos\frac{k\pi}{2n+2}\right).\tag3.24
$$
By (2.2), the eigenvalues of $Q_n$ are $2+2\cos\frac{j\pi}{n+1}=4\cos^2\frac{j\pi}{2n+2}$, $j=1,\dotsc,n$,
while those of $Q'_n$ are $-2+2\cos\frac{j\pi}{n+1}$, $j=1,\dotsc,n$; the latter listed in reversed
order thus equal $-4\cos^2\frac{j\pi}{2n+2}$, $j=1,\dotsc,n$. 
Theorem 3.3 and (3.24) imply then,
after some manipulation, that
$$
\Pe(HOD_{n-1};x)=\prod_{{\scriptstyle 1\leq j<k\leq 2n+1}\atop {\scriptstyle j+k\leq 2n+1}}
\left(x-4\cos\frac{j\pi}{2n+2}\cos\frac{k\pi}{2n+2}\right).\tag3.25
$$
The statement follows now by (3.23).\epf

\proclaim{Corollary 3.7} The number of spanning trees of the halved Aztec diamond $HOD_n$ is
given by
$$
\te(HOD_n)=\prod_{{\scriptstyle 1\leq j<k\leq 2n}\atop {\scriptstyle j+k\leq 2n}}
\left(4-4\cos\frac{j\pi}{2n+1}\cos\frac{k\pi}{2n+1}\right).\tag3.26
$$
\endproclaim

\pf The graph obtained from the planar dual of $HOD_n$ by removing the vertex corresponding to the
infinite face is readily seen to be isomorphic to $HMD_{n-1}$. Therefore, by Theorem 2.2 we obtain
$$
\te(HOD_n)=\Pe(HMD_{n-1};4).\tag3.27
$$
By \cite{\Kn} we have
$$
\Pe(MD_n;x)=\prod_{j=1}^n\left(x-4\cos^2\frac{j\pi}{2n+1}\right)
\prod_{1\leq j<k\leq 2n}\left(x-4\cos\frac{j\pi}{2n+1}\cos\frac{k\pi}{2n+1}\right).\tag3.28
$$
It is not hard to show that the eigenvalues of the adjacency matrix of $R_n$ are 
$2+2\cos\frac{2j\pi}{2n+1}=4\cos^2\frac{j\pi}{2n+1}$, $j=1,\dotsc,n$, while the eigenvalues of
the adjacency matrix of $R'_n$ are just the negatives of the former. 
Using this, one sees after some manipulation
that
$$
\Pe(HMD_{n-1};x)=\prod_{{\scriptstyle 1\leq j<k\leq 2n}\atop {\scriptstyle j+k\leq 2n}}
\left(x-4\cos\frac{j\pi}{2n+1}\cos\frac{k\pi}{2n+1}\right).\tag3.29
$$
Together with (3.27), this implies (3.26).\epf

\medskip
\flushpar
{\smc Remark 3.8.} It is amusing to note that the range for $(j,k)$ in the product giving
the characteristic polynomial of $HMD_n$ can naturally be regarded as the vertex set of $HOD_{n}$,
and vice-versa.



\mysec{4. Symmetry classes of spanning trees of Aztec diamonds}

It is easy to see that no spanning tree of an Aztec diamond $AD_n$ can be symmetric with respect to a
a symmetry axis $\ell$ of the diamond that makes a $45^\circ$ angle with the horizontal. Indeed, suppose 
$AD_n$ had a spanning tree $T$ symmetric about $\ell$. Let $a$ and $b$ be two distinct vertices of
$AD_n$ on $\ell$. Since $T$ is connected, there exists a path $P$ in $T$ connecting $a$ to $b$. Since
$T$ is symmetric about $\ell$, the reflection $P'$ of $P$ across $\ell$ is also contained in $T$. But
no path in $AD_n$ from $a$ to $b$ is invariant under this reflection, so $P$ and $P'$ are two distinct 
paths connecting two vertices of the tree $T$, a contradiction.

The same argument works for the odd Aztec diamonds $OD_n$ when $n\geq2$. The mixed Aztec diamonds 
$MD_n$ do not possess a symmetry axis at a $45^\circ$ angle from the horizontal.

Let $h$, $v$ and $r$ denote the symmetry across the horizontal, the symmetry across the vertical, and
the rotation by $90^\circ$, respectively; view them as elements of the symmetry group of $AD_n$ or $OD_n$.
The above two paragraphs
imply that there are a total of five inequivalent symmetry classes of spanning trees of 
$AD_n$ and $OD_n$: the base case---trees with no symmetry requirement; 
horizontally symmetric trees---invariant under the action of $h$;
horizontally and vertically symmetric trees---invariant under the action of $\langle h,v\rangle$;
trees invariant under the action of $r^2$; and trees invariant under the action of $r$.
Since $MD_n$ is not $r$-invariant, there are only four symmetry classes for its spanning trees.

We denote the number of spanning trees of a graph $G$ that are invariant under the action of the
group of symmetries $H$ by $\te_H(G)$.


The base case was done by Knuth \cite{\Kn}. In this section we provide product formulas for all but three of the remaining
cases. The latter are phrased as open problems. Our arguments will express these numbers in terms of the number of perfect 
matchings of three
families of subgraphs of the infinite grid $\Z^2$. Explicit formulas for the latter are deduced using
results from the previous sections and the factorization theorem for perfect matchings of \cite{\Cone}.

To define our families of graphs, it will be convenient to denote by $NE(i,j)$ the infinite zig-zag 
lattice path
in $\Z^2$ starting at the lattice point $(i,j)$ and taking alternately two unit steps 
north and two east; 
$NW(i,j)$ is defined analogously, alternating between two steps north and two steps west.

Let $G_{2n}$ be the subgraph of the infinite grid $\Z^2$ induced by the vertices $(i,j)$ with
$0\leq i,j\leq 2n-1$. We define $A_n$ to be the subgraph of $G_{2n}$ induced by its vertices on
or above the path $NE(1,0)$. $B_n$ is defined analogously,
but using the path $NE(0,0)$ instead of $NE(0,1)$ ($A_5$ and $B_5$ are illustrated in Figures 4.1 and
4.2, respectively; the two dots in the latter emphasize that there are two vertices at the indicated
positions).

\topinsert
\twoline{\mypic{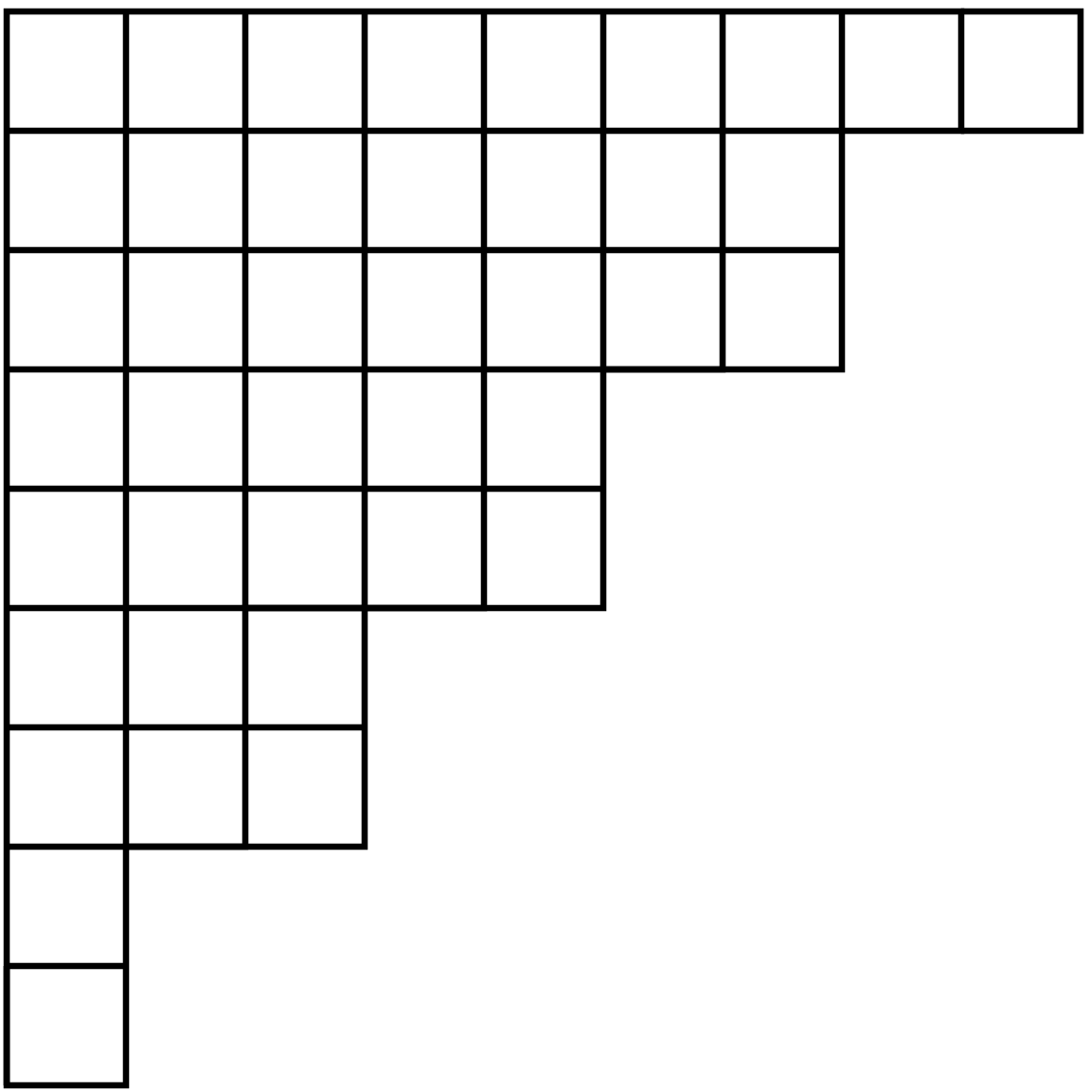}}{\mypic{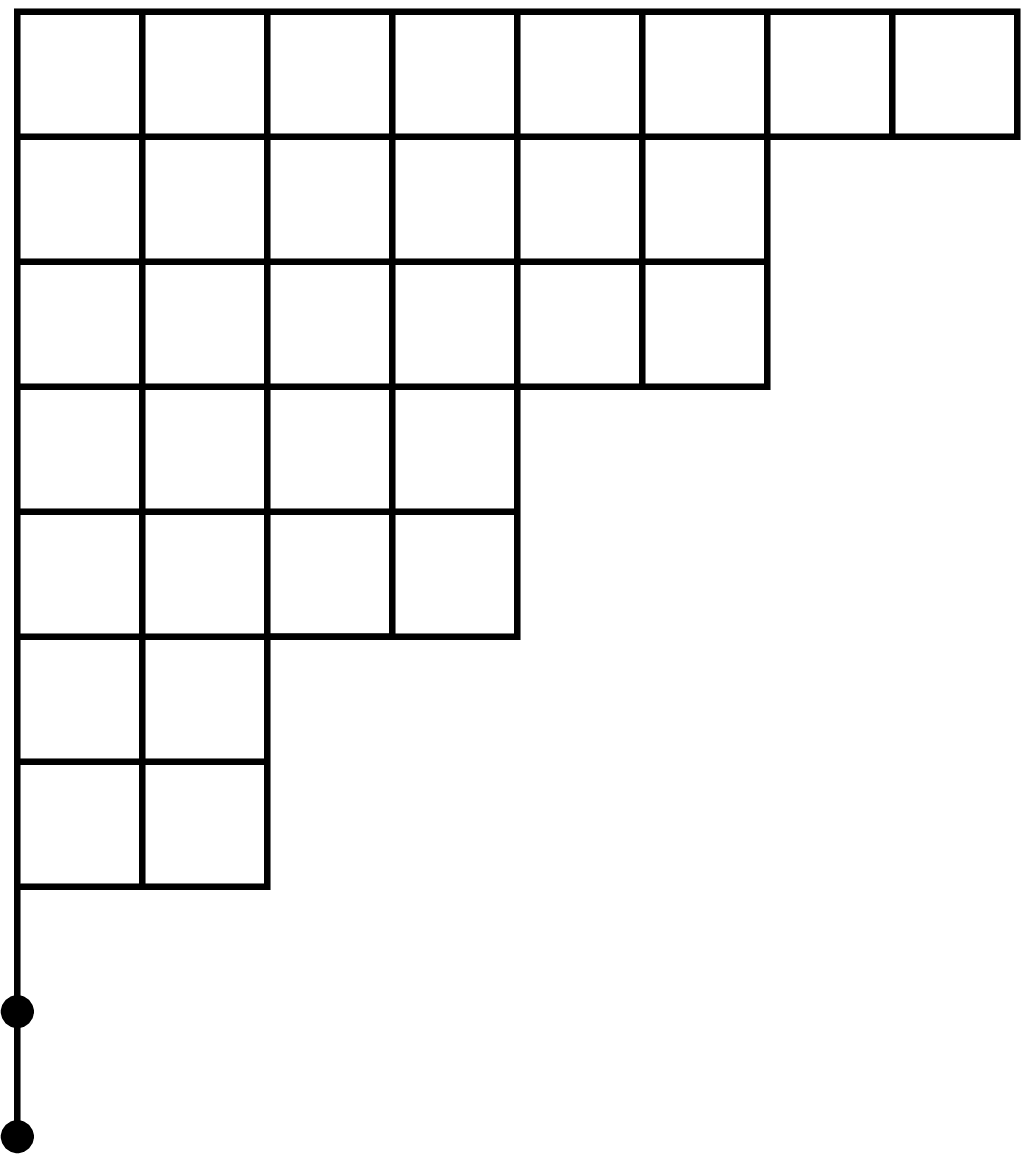}}
\twoline{\smc Figure 4.1{\rm . $A_5$.}}{\smc Figure 4.2{\rm . $B_5$.}}

\bigskip
\twoline{\mypic{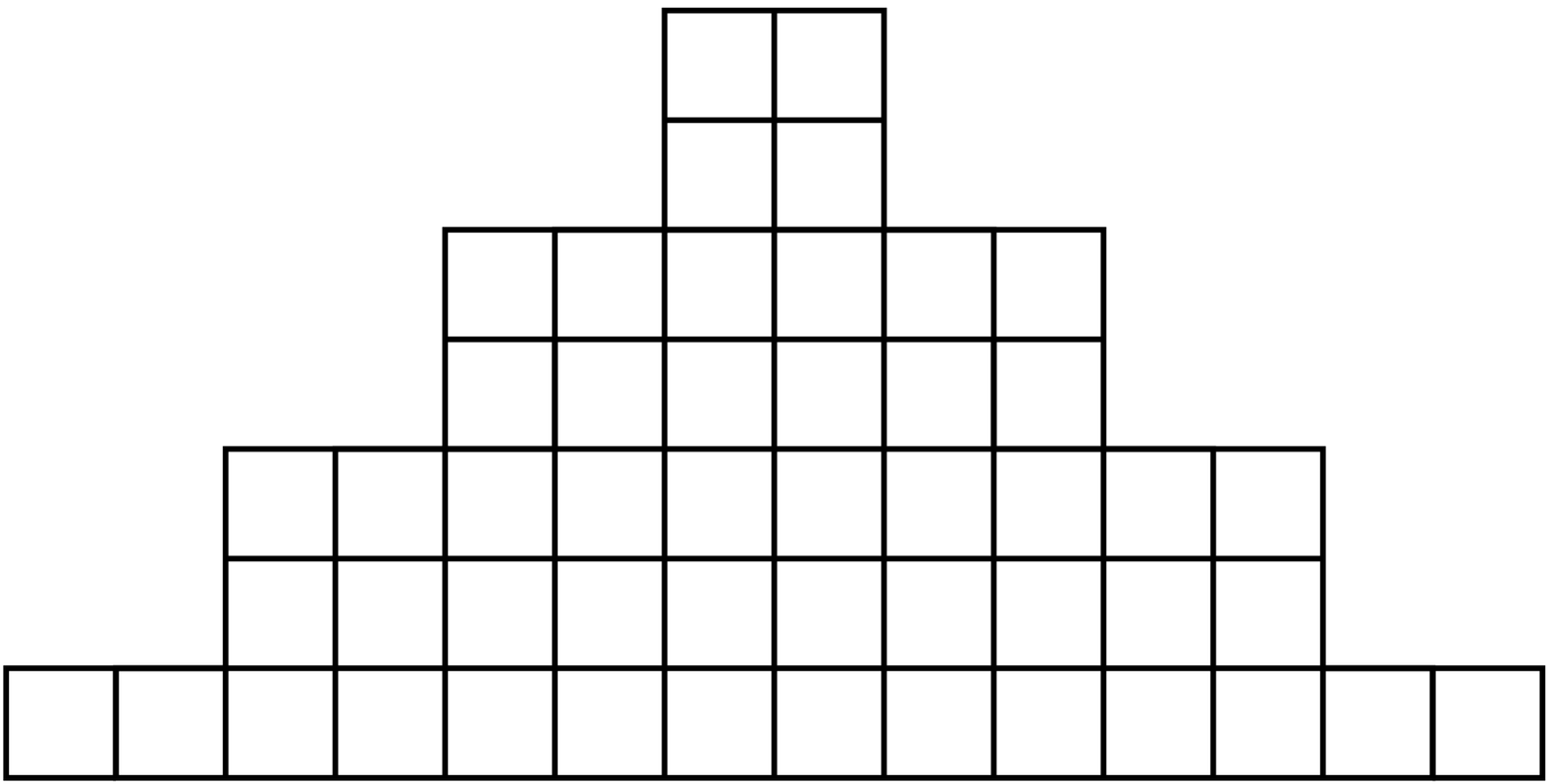}}{\mypic{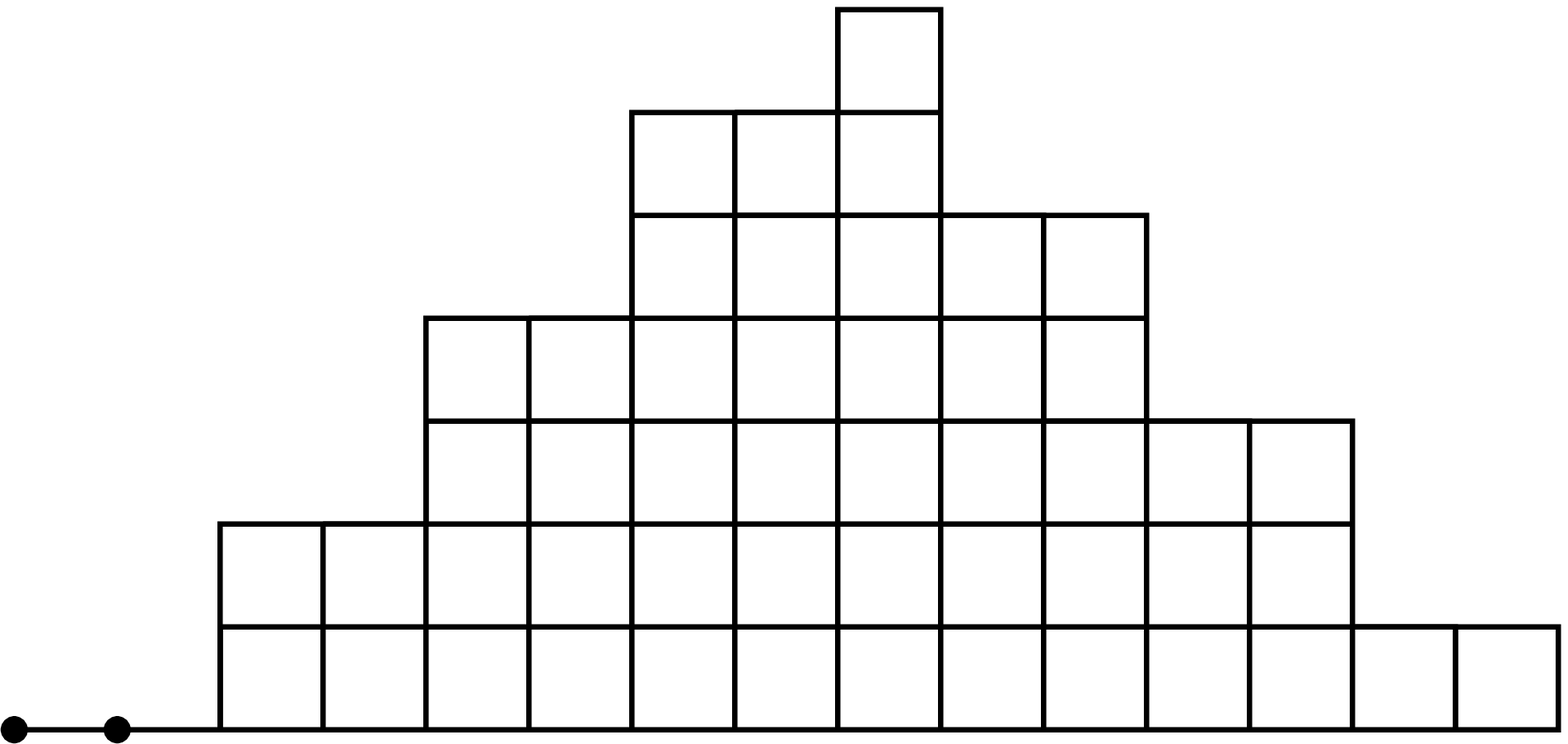}}
\medskip
\twoline{\smc Figure 4.3{\rm . $C_8$.}}{\smc Figure 4.4{\rm . $D_8$.}}
\endinsert

Define $C_n$ to be the subgraph of $\Z^2$ induced by the lattice points $(i,j)$, $j\geq0$, that are
on or below both $NE(-n+1,-1)$ and $NW(n-1,-1)$; $C_8$ is shown in Figure 4.3.

A fourth family will be relevant in the proof of (4.3); it will also be used in Section 5.
We define $D_n$ to be the subgraph of $\Z^2$ induced by the lattice points $(i,j)$, $j\geq0$, that are
on or below both $NE(-n,-2)$ and $NW(n-1,-1)$; $D_8$ is shown in Figure 4.4.

For a graph $G$ we denote the number of its perfect matchings by $\M(G)$ (if $G$ is weighted, $\M(G)$ denotes
the sum of the weights of all its perfect matchings, the weight of a matching being the product of the
weights of its constituent edges).

\proclaim{Lemma 4.1} We have
$$
\M(A_n)=\frac{1}{2^n}\prod_{1\leq j\leq k\leq n}
\left(4-2\cos\frac{j\pi}{n+1}-2\cos\frac{k\pi}{n+1}\right),\tag4.1
$$
$$
\M(B_n)=
\prod_{1\leq j<k\leq n}\left(4\cos^2\frac{j\pi}{2n+1}+4\cos^2\frac{k\pi}{2n+1}\right),
\tag4.2
$$
$$
\M(C_n)=\frac{1}{ 2^{ 2\left\lfloor \frac{n+1}{2} \right\rfloor } }
\prod_{{\scriptstyle 1\leq j\leq k\leq n}\atop{\scriptstyle j+k \,\text{\rm even}}  }
\left(4-2\cos\frac{j\pi}{n+1}-2\cos\frac{k\pi}{n+1}\right),
\tag4.3
$$
and
$$
\M(D_n)=\frac{1}{ 2^{ \left\lfloor \frac{n}{2} \right\rfloor } }
\prod_{{\scriptstyle 1\leq j<k\leq n}\atop{\scriptstyle j+k \,\text{\rm odd}}  }
\left(4-2\cos\frac{j\pi}{n+1}-2\cos\frac{k\pi}{n+1}\right).
\tag4.4
$$

\endproclaim

The following result of Temperley \cite{\Temp} will be useful in our proof. Let $G=(V,E)$ be a finite
connected planar 
graph; denote the set of its bounded faces by $F$. Pick a point $v(e)$ in the interior of each edge $e$ 
of $G$; let $V_e:=\{v(e):e\in E\}$. 
Pick also a point $v(f)$ in the interior of each bounded face $f$ of $G$, and let $V_f:=\{v(f):f\in F\}$.
Define $\Te(G)$ to be the graph whose vertex set is $V\cup V_e\cup V_f$, with edge set consisting of
the pairs $\{v,v(e)\}$, $v\in V$, $e\in E$, $v$ incident to $e$, and  
$\{v(e),v(f)\}$, $e\in E$, $f\in F$, $e$ incident to $f$. In particular, if $G_n$ is the $n\times n$
grid graph, $\Te(G_n)$ is the refined grid $G_{2n-1}$.

\proclaim{Theorem 4.2 (Temperley \cite{\Temp}\cite{\Lov, Problem 4.30})} 
For any $v\in V$ incident to the infinite face of $G$ we have
$$
\te(G)=\M(\Te(G)\setminus v).
$$

\endproclaim

We will also employ the following ``Factorization Theorem'' we obtained in \cite{\Cone}.

\proclaim{Theorem 4.3 (Factorization Theorem \cite{\Cone,Theorem\,2.1})} 
Suppose the weighted bipartite planar graph $G$ is drawn on the plane so that it is symmetric with respect to a 
horizontal symmetry axis $\ell$. Label the vertices on $\ell$ from left to right by $a_1,b_1,a_2,b_2,\dotsc,a_{k},b_{k}$,
and assume they form a cut set. Properly color the vertices of $G$ white and black, so that $a_1$ is white.

Let $G'$ be the graph obtained from $G$ by 
\medskip
$(i)$ removing the edges above $\ell$ incident to any white $a_i$ or any black $b_j$;

$(ii)$ removing the edges below $\ell$ incident to any black $a_i$ or any white $b_j$; and

$(iii)$ reducing by $1/2$ the weight of any edge of $G$ lying along $\ell$.

\medskip
Let $G^+$ be the subgraph of $G'$ induced by the vertices above $\ell$, the black $a_i$'s, and white $b_j$'s.
Let $G^-$ be the subgraph of $G'$ induced by the vertices below $\ell$, the white $a_i$'s, and black $b_j$'s.
Then
$$
\M(G)=2^k \M(G^+)\M(G^-).
$$

\endproclaim

{\it Proof of Lemma 4.1.} Apply the factorization theorem for perfect matchings \cite{\Cone, Theorem 2.1} to
the graph $H$ obtained from the square grid $G_{2n+1}$ by removing its top right vertex $v$, 
with respect to the symmetry axis of $H$ (see Figure 4.5). 
One of the resulting subgraphs of $H$ is then isomorphic
to $A_n$ (the boundary of this subgraph is indicated by the upper bold outline in Figure 4.5), while the other 
can be regarded as $\Te(QAD_{n+1})\setminus v$ (the boundary of the latter is traced out by the lower bold outline
in Figure 4.5). Since $H$ itself can be
regarded as $\Te(G_{n+1})\setminus v$, we obtain by the factorization theorem that
$$
\M(\Te(G_{n+1})\setminus v)=2^n\M(A_n)\M(\Te(QAD_{n+1})\setminus v).\tag4.5
$$
By Theorem 4.2, the left hand side equals $\te(G_{n+1})$, which in turn by Theorem 2.2 equals 
$\Pe(G_n;4)$; (2.3) provides an explicit formula for the latter. On the other hand, by Theorem~4.2,
$\M(\Te(QAD_{n+1})\setminus v)=\te(QAD_{n+1})$, for which Equation (2.4) gives a product expression. Solving
for $\M(A_n)$ in (4.5) yields then (4.1).

\topinsert
\twoline{\mypic{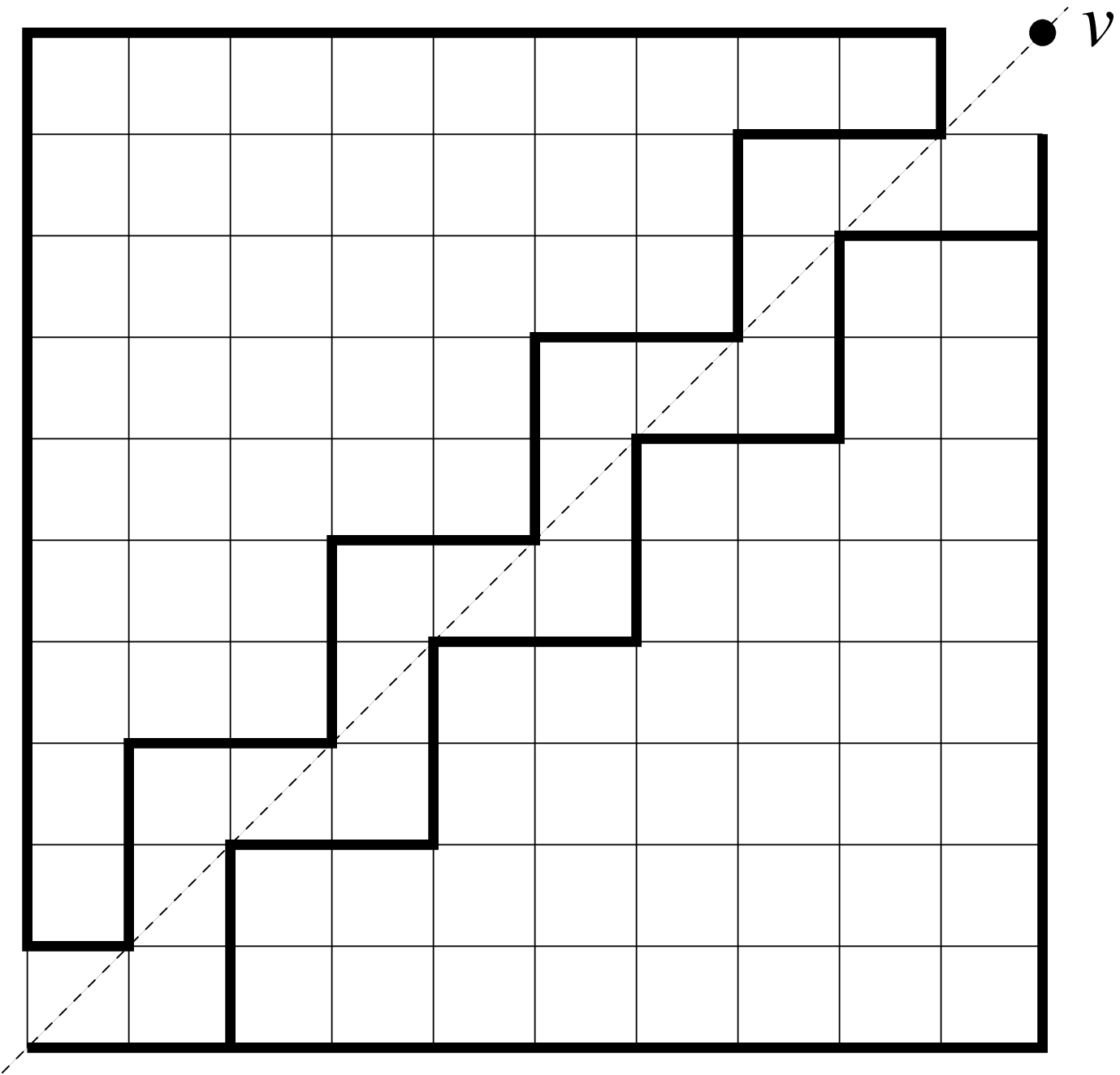}}{\mypic{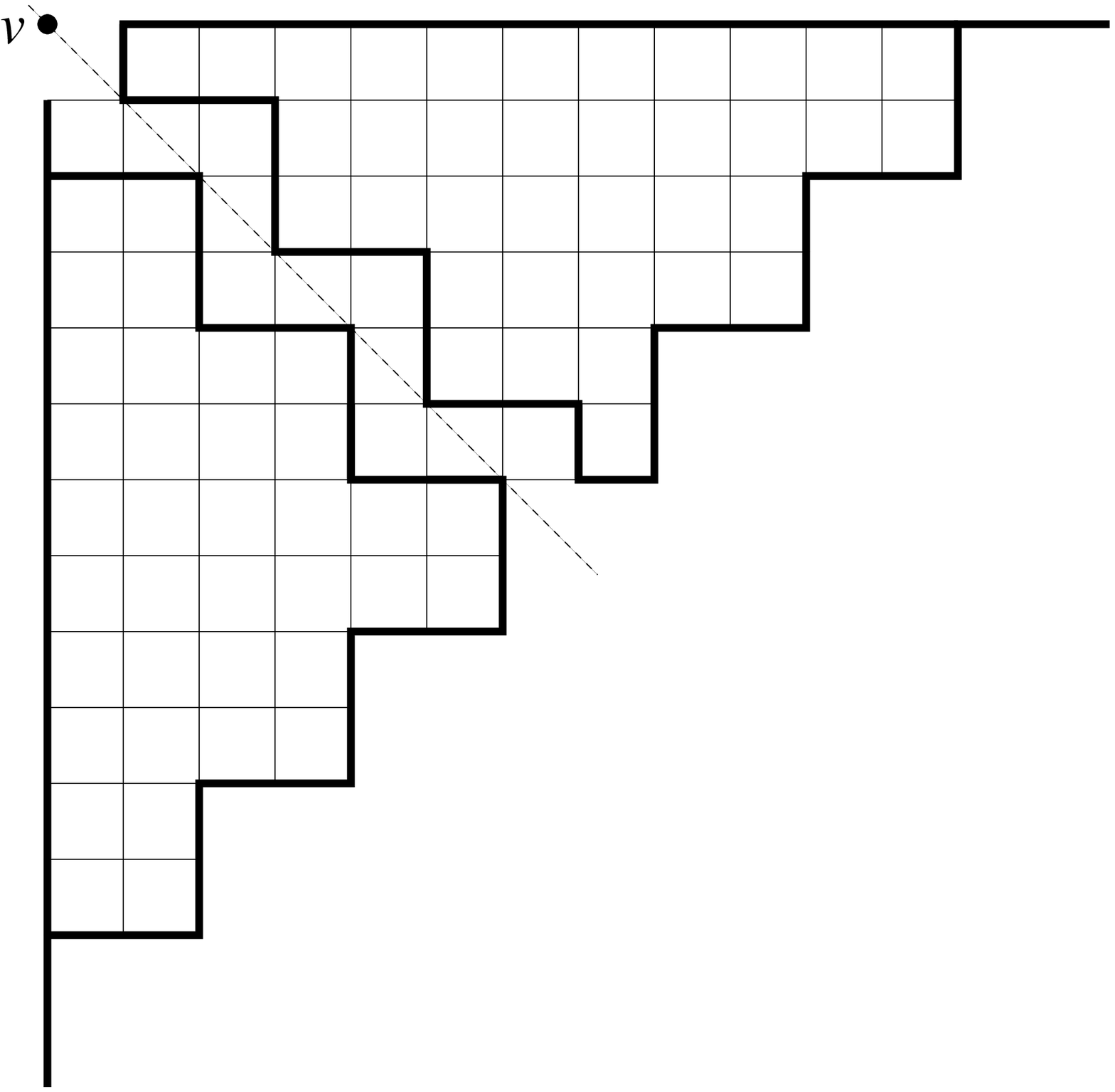}}
\twoline{\smc Figure 4.5.{\rm }}{\smc Figure 4.6.{\rm }}
\endinsert

Equality (4.2) follows by applying the factorization theorem to the square grid $G_{2n}$ and its
diagonal symmetry axis. The resulting subgraphs are both isomorphic to $B_n$. We obtain
$$
\M(G_{2n})=2^n\,[\M(B_n)]^2.\tag4.6
$$
A classical result of Kasteleyn, Temperley and Fisher (see e.g. \cite{\Lov, Problem 4.29}) gives
$$
\M(G_{2n})=2^{2n^2}\prod_{1\leq j,k\leq n}\left(\cos^2\frac{j\pi}{2n+1}+\cos^2\frac{k\pi}{2n+1}\right).
\tag4.7
$$
Plugging this into (4.6) and using the well known identity $\prod_{j=1}^n 2\cos\frac{j\pi}{2n+1}=1$
we obtain (4.2).

Next, we prove (4.4). Apply the factorization theorem to the graph $\Te(QAD_{n+1})\setminus v$, with 
$v$ chosen to be both on the symmetry axis and on the infinite face of $\Te(QAD_{n+1})$ (see Figure 4.6). 
One of the
resulting subgraphs is then $D_n$, while the other is $\Te(HOD_{n/2})\setminus v$, when $n$ is even,
and $\Te(HMD_{(n+1)/2})\setminus v$, when $n$ is odd. The formulas provided by Corollaries 2.3, 3.3 
and 3.4 can then be used to obtain an expression for $\M(D_n)$. After some manipulation one obtains
(4.4).

Now $\M(C_n)$ can be worked out by applying the factorization theorem to the graph $A_{n}$. It is readily
seen that one of the resulting subgraphs is $D_{n}$, while the other is $C_n$. We obtain
$$
\M(A_{n})=2^{ \left\lfloor \frac{n+1}{2} \right\rfloor }\M(D_{n})\M(C_n).
$$
Using (4.1) and (4.4) this implies (4.3). \epf

\proclaim{Theorem 4.4} The only nonempty non-trivial symmetry class of spanning trees of $AD_n$
is that of horizontally symmetric trees. We have
$$
\te_{\langle h\rangle}(AD_n)=
2n \te(HMD_n)=2n \prod_{{\scriptstyle 1\leq j<k\leq 2n-1}\atop {\scriptstyle j+k\leq 2n-1}}
\left(4-4\cos\frac{j\pi}{2n}\cos\frac{k\pi}{2n}\right).\tag4.8
$$

\endproclaim

\pf Suppose $T$ is a horizontally and vertically symmetric spanning connected subgraph of $AD_n$. 
Since $T$ is connected,
it contains some edge $e$ crossed by the horizontal symmetry axis $\ell_h$. The mirror image $e'$ of $e$ 
across the vertical symmetry axis $\ell_v$ is also in $T$, as $T$ is vertically symmetric. Let $a$ be the 
top vertex of $e$, and $a'$ its mirror image across $\ell_v$. Since $T$ is connected, it contains a path
$P$ from $a$ to $a'$. But then the union of the reflection of $P$ across $\ell_h$ with the edges $e$ and
$e'$ yields another path in $T$, distinct from $P$, connecting $a$ to $a'$; so $T$ cannot be tree.

A similar argument shows that, if $T$ is a spanning connected subgraph of $AD_n$ which is invariant under 
rotation by $180^\circ$, and $a$ is say one of the four centralmost vertices of $AD_n$, then
$T$ must contain two distinct paths connecting $a$ to its image through rotation by $180^\circ$. Since
quarter turn invariant spanning trees are obviously also half turn invariant, to conclude the proof it suffices
to prove (4.8).

Let $T$ be a horizontally symmetric spanning tree of $AD_n$. The argument in the first paragraph of this
proof implies that $T$ contains precisely one edge $e$ that crosses $\ell_h$. The removal of this edge from
$T$ leaves two mirror images of a spanning tree of $HMD_n$, and defines a $2n$-to-$1$ map from the set of
horizontally symmetric spanning trees of $AD_n$ to the set of spanning trees of $HMD_n$. This proves
the first equality of (4.8); the second follows by Corollary 3.6. \epf

\proclaim{Theorem 4.5} 
The non-trivial symmetry classes of spanning trees of $OD_n$ not involving rotations
are enumerated by
$$
\te_{\langle h\rangle}(OD_n)=
\M(C_{2n-1})=\frac{1}{ 2^{ 2n } }
\prod_{{\scriptstyle 1\leq j\leq k\leq 2n-1}\atop{\scriptstyle j+k \,\text{\rm even}}  }
\left(4-2\cos\frac{j\pi}{2n}-2\cos\frac{k\pi}{2n}\right)\tag4.9
$$
and
$$
\te_{\langle h,v\rangle}(OD_n)=
\M(A_{n-1})=\frac{1}{2^{n-1}}\prod_{1\leq j\leq k\leq n-1}
\left(4-2\cos\frac{j\pi}{n}-2\cos\frac{k\pi}{n}\right).\tag4.10
$$

\endproclaim

\pf
Let $T$ be an $h$-invariant spanning tree of $OD_n$, and let $a$ and $b$ be the two vertices of $OD_n$
that are both on its infinite face and on its horizontal symmetry axis $\ell_h$. Let $P$ be the path in $T$
connecting $a$ to $b$. By the uniqueness of such a path it follows that $P$ is $h$-invariant. Clearly,
this implies that $P$ consists of the $2n$ edges of $OD_n$ along $\ell_h$.

It follows that the $h$-invariant spanning trees of $OD_n$ are in bijection with the spanning trees of $HOD_n$
that contain all the bottom $2n$ edges. In turn, the latter can be identified with the spanning trees of 
the graph $G$ obtained from $HOD_{n-1}$ by including a new vertex $w$ connected  to its $2n-1$ bottommost
vertices. However, $\Te(G)\setminus w$ is readily seen to be isomorphic to $C_{2n-1}$. Then (4.9) follows
from Theorem 4.2 and (4.3).

An argument similar to the one in the first paragraph of this proof shows that 
the $\langle h,v\rangle$-invariant spanning trees can be identified with the spanning trees of
$QAD_n$ which contain all the $2n$ edges along the straight line portions of its boundary. In analogy
to the previous paragraph, by Theorem 4.2 the number of the latter is seen to equal $\M(A_{n-1})$. Thus
(4.10) follows by~(4.1).~\epf

\proclaim{Theorem 4.6} The non-trivial symmetry classes of spanning trees of $MD_n$
not involving rotations are enumerated by
$$
\te_{\langle h\rangle}(MD_n)=
\M(C_{2n-2})=\frac{1}{ 2^{ 2n-2} }
\prod_{{\scriptstyle 1\leq j\leq k\leq 2n-2}\atop{\scriptstyle j+k \,\text{\rm even}}  }
\left(4-2\cos\frac{j\pi}{2n-1}-2\cos\frac{k\pi}{2n-1}\right)\tag 4.11
$$
and
$$
\te_{\langle h,v\rangle}(MD_n)=
\M(B_{n-1})=2^{n-1}\prod_{j=1}^{n-1}\cos\frac{j\pi}{2n-1}
\prod_{1\leq j<k\leq n-1}\left(4\cos^2\frac{j\pi}{2n-1}+4\cos^2\frac{k\pi}{2n-1}\right).
\tag4.12
$$

\endproclaim

\pf Arguments analogous to the ones in the proof of Theorem 4.5 show that 
$\te_{\langle h\rangle}(MD_n)=\M(C_{2n-2})$, and
$\te_{\langle h,v\rangle}(MD_n)=\M(B_{n-1})$. Then (4.11) and (4.12) follow by (4.3) and~(4.2). \epf

\proclaim{Open Problem 4.7} Find formulas for $\te_{\langle r\rangle}(OD_n)$, $\te_{\langle r^2\rangle}(OD_n)$,
and $\te_{\langle r^2\rangle}(MD_n)$.

\endproclaim

\mysec{5. Symmetry classes of perfect matchings of odd squares with a central unit hole}

The existence of formula (4.7) invites one to look for an analog in the case of
odd square grids. One very natural candidate is the graph $H_n$ obtained from the $(2n+1)\times(2n+1)$ 
square grid by removing its central vertex (see Figure 5.1).

As it turns out, the numbers $\M(H_n)$ do not seem to have many factors in their prime factorization (see
Section 7).
This makes it implausible for them to possess a ``nice'' product expression (e.g. similar in style to
the product formulas that occur in this paper).

However, all non-trivial symmetry classes of perfect matchings of $H_n$ turn out to be enumerated
by variations of the product formulas we have encountered in the previous sections.

To present our results we will need, besides the graphs introduced in the previous section, to enumerate
the perfect matchings of certain weighted versions of $A_n$ and $B_n$. 

Define $\tilde{A}_n$ to be the weighted graph obtained from $A_n$ by weighting its top $2n-1$ edges
by $1/2$, and keeping weight 1 for all its other edges. Define $\tilde{B}_n$ to be the weighted graph 
obtained from $B_n$ by weighting its top $2n-2$ edges by $1/2$, and keeping weight 1 for all its other 
edges. For a weighted graph $G$, $\M(G)$ denotes the sum of the weights of all its perfect matchings, the
weight of a matching being equal to the product of the weights on its constituent edges. 

\proclaim{Lemma 5.1}
$$
\M(\tilde{A}_n)=\frac{1}{2^{3n}}\prod_{j=1}^{2n}\left(4-4\cos\frac{j\pi}{2n+1}\right)
\prod_{1\leq j<k\leq n}\left(4-2\cos\frac{2j\pi}{2n+1}-2\cos\frac{2k\pi}{2n+1}\right),\tag 5.1
$$
and
$$
\M(\tilde{B}_n)=\frac{1}{2^{n-1}}
\prod_{{\scriptstyle 1\leq j < k \leq 2n-1}\atop{\scriptstyle j+k\leq 2n-1,\ j+k \,\text{\rm odd}}}
\left(4-4\cos\frac{j\pi}{2n}\cos\frac{k\pi}{2n}\right).\tag5.2
$$

\endproclaim

Before giving the proof of the above lemma, we show how it implies the following result.

Clearly, no perfect matching of $H_n$ can be symmetric with respect to the diagonal. Thus the reflection
$h$ across the horizontal, the reflection $v$ across the vertical, and the rotation $r$ by $90^\circ$
generate the relevant symmetry group. As is Section 4, there are four inequivalent non-trivial 
symmetry classes, corresponding to the subgroups $\langle h \rangle$, $\langle h,v \rangle$, 
$\langle r^2 \rangle$, and $\langle r \rangle$, respectively. Denote by $\M_K(G)$ the number of
perfect matchings of $G$ that are invariant under the symmetry group~$K$.

\proclaim{Theorem 5.2} We have
$$
\M_{\langle h \rangle}(H_{2n})=\prod_{j=0}^{n-1}\prod_{k=0}^{2n}
\left(4-2\cos\frac{(2j+1)\pi}{2n+1}-2\cos\frac{k\pi}{2n+1}\right),\tag5.3
$$
$$
\M_{\langle h,v \rangle}(H_{2n})=\prod_{1\leq j,k\leq n}
\left(4\cos^2\frac{j\pi}{2n+1}+4\cos^2\frac{k\pi}{2n+1}\right),\tag5.4
$$
$$
\M_{\langle r^2 \rangle}(H_n)=2^{n-2\left\lfloor \frac{n}{2} \right\rfloor}
\prod_{{\scriptstyle 1\leq j<k\leq n}\atop{\scriptstyle j+k \,\text{\rm odd}}  }
\left(4-2\cos\frac{j\pi}{n+1}-2\cos\frac{k\pi}{n+1}\right)^2,\tag5.5
$$
$$
\M_{\langle r \rangle}(H_{2n-1})=2\prod_{1\leq j\leq k \leq n-1}
\left(4-2\cos\frac{j\pi}{n}-2\cos\frac{k\pi}{n}\right)
\!\!\prod_{{\scriptstyle 1\leq j < k \leq 2n-1}\atop{\scriptstyle j+k\leq 2n-1,\ j+k \,\text{\rm odd}}}
\left(4-4\cos\frac{j\pi}{2n}\cos\frac{k\pi}{2n}\right),\tag5.6
$$
and
$$
\align
&
\M_{\langle r \rangle}(H_{2n})=\frac{1}{2^n}
\prod_{j=1}^{2n}\left(4-4\cos\frac{j\pi}{2n+1}\right)
\\
&\ \ \ \ \ \ \ \ \
\times
\prod_{1\leq j < k \leq n}
\left(4-2\cos\frac{2j\pi}{2n+1}-2\cos\frac{2k\pi}{2n+1}\right)
\left(4\cos^2\frac{j\pi}{2n+1}+4\cos^2\frac{k\pi}{2n+1}\right).\tag5.7
\endalign
$$

\endproclaim

\pf For a perfect matching of $H_n$ to be $h$-invariant, all vertices along the horizontal symmetry axis
$\ell_h$ must be matched by edges along $\ell_h$; this is possible only when $n$ is even. The
$h$-invariant perfect matchings of $H_{2n}$ can then be identified with perfect matchings of the 
rectangular $2n\times (4n+1)$ grid graph. The latter were enumerated by Kasteleyn, Fisher and Temperley;
(5.3) follows from the form given in \cite{\KPW}.

By a similar argument, the $\langle h,v \rangle$-invariant perfect matchings of $H_{2n}$ can be 
identified with perfect matchings of the square grid graph $G_{2n}$. Thus (5.4) follows by (4.7).

To prove (5.5), note that $\M_{\langle r^2 \rangle}(H_n)=\M(H_n/\langle r^2 \rangle)$, where 
$H_n/\langle r^2 \rangle$ is the orbit graph of the action of $r^2$ on $H_n$. This orbit graph
can be regarded as being obtained from the subgraph of $H_n$ induced by its vertices on or under a 
diagonal $d$ by identifying the pairs of vertices on $d$ that are at the same distance from the center;
Figure 5.1 shows $H_6/\langle r^2 \rangle$.

\topinsert
\centerline{\mypic{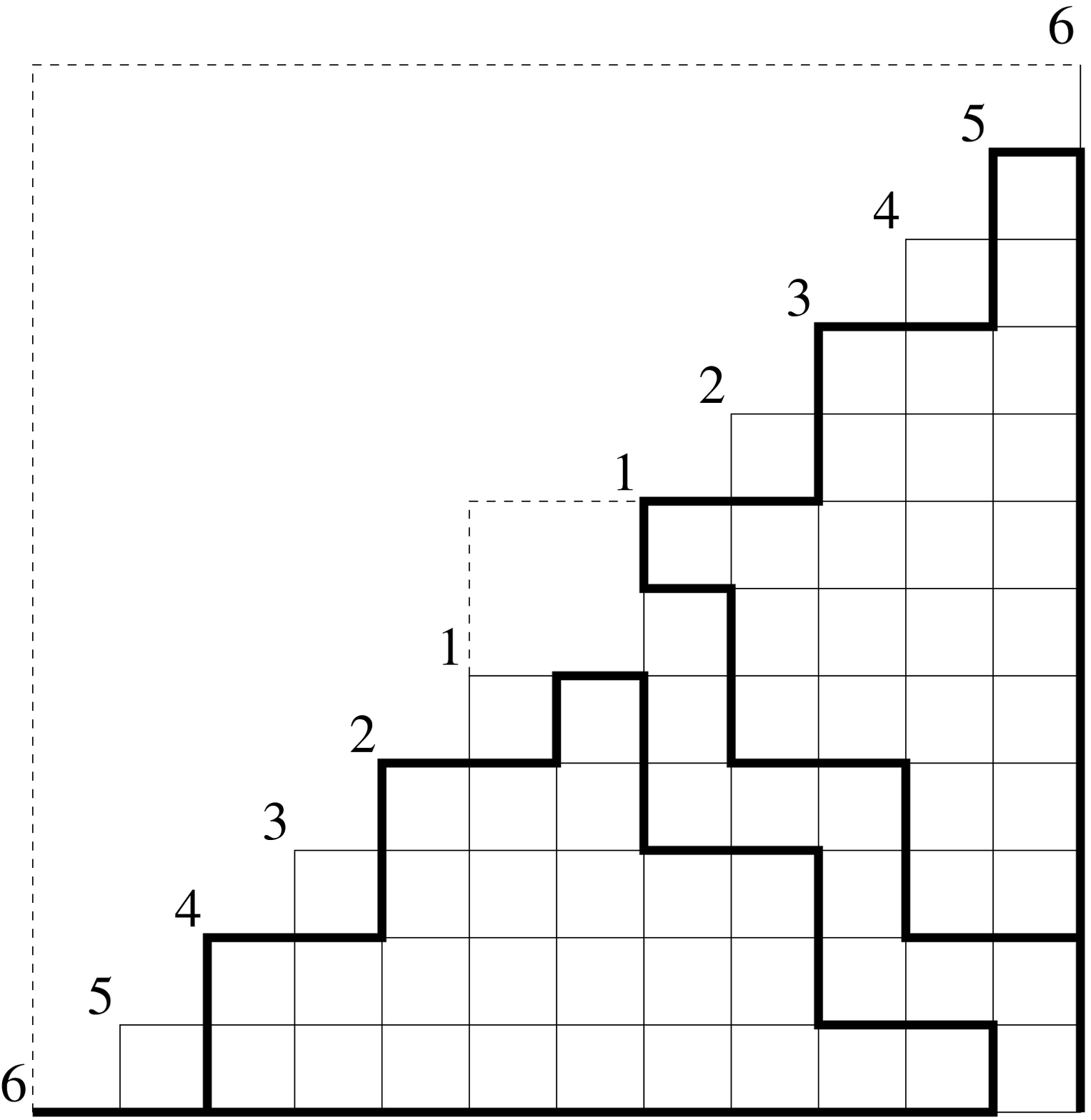}}
\medskip
\centerline{{\smc Figure~5.1.} {\rm Factorization applied to $H_6/\langle r^2 \rangle$.}}
\bigskip
\endinsert

Apply the factorization theorem of \cite{\Cone} to $H_n/\langle r^2 \rangle$ with respect to the diagonal
perpendicular to $d$. The resulting graphs both turn out to be isomorphic to $D_n$. We obtain
$$
\M_{\langle r^2 \rangle}(H_n)=2^n[\M(D_n)]^2.
$$
Formula (4.4) implies then (5.5).

We now turn to the proof of (5.6). Let $H_n/\langle r \rangle$ be the orbit graph of the action of $r$ 
on $H_n$. By definition, $\M_{\langle r \rangle}(H_n)=\M(H_n/\langle r \rangle)$. The graph
$H_n/\langle r \rangle$ is obtained from the subgraph of $H_n$ induced by its vertices on or under 
the union of its two diagonals 
by identifying the pairs of vertices on this union that are at the same distance from the center;
Figure 5.2 shows $H_7/\langle r \rangle$.

Apply the factorization theorem of \cite{\Cone} to $H_n/\langle r \rangle$ with respect to the vertical
symmetry axis. The shape of the resulting subgraphs is different for different parities of $n$.

For odd $n$, one of the graphs is $A_{(n-1)/2}$, while the other is $\tilde{B}_{(n+1)/2}$ (see Figure 5.2).
We obtain
$$
\M_{\langle r \rangle}(H_{2n-1})=2^{2n-1}\M(A_{n-1})\M(\tilde{B}_{n}).
$$
Using (4.1) and Lemma 5.1, this implies (5.6).

On the other hand, when $n$ is even, one of the subgraphs obtained by applying the factorization theorem to
$H_n/\langle r \rangle$ is $B_{n/2}$, while the other is $\tilde{A}_{n/2}$ (see Figure 5.3). 
Thus
$$
\M_{\langle r \rangle}(H_{2n})=2^{2n}\M(B_{n})\M(\tilde{A}_{n}).
$$
Formula (5.7) follows then by (4.2) and Lemma 5.1. \epf

\topinsert
\twoline{\mypic{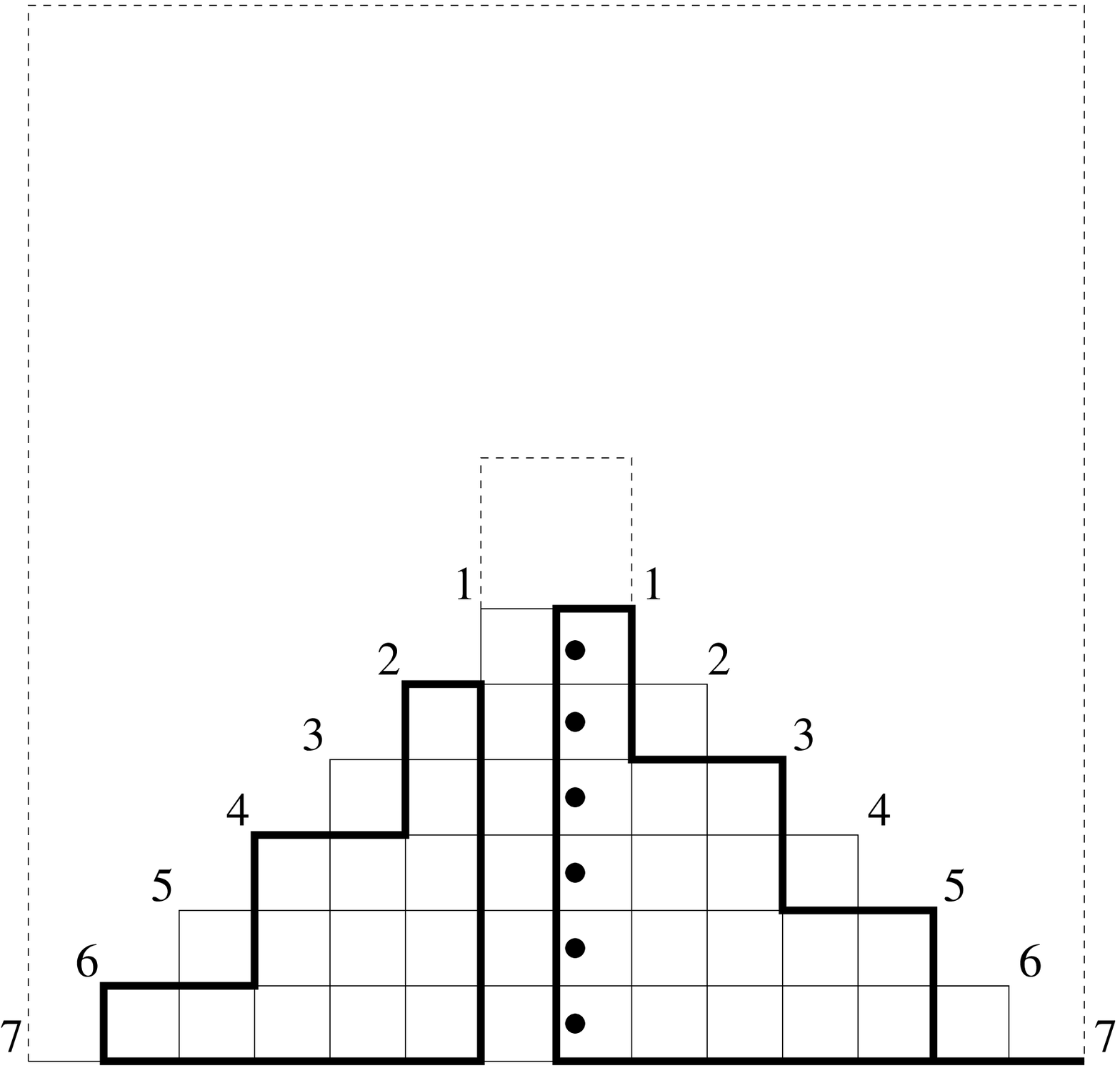}}{\mypic{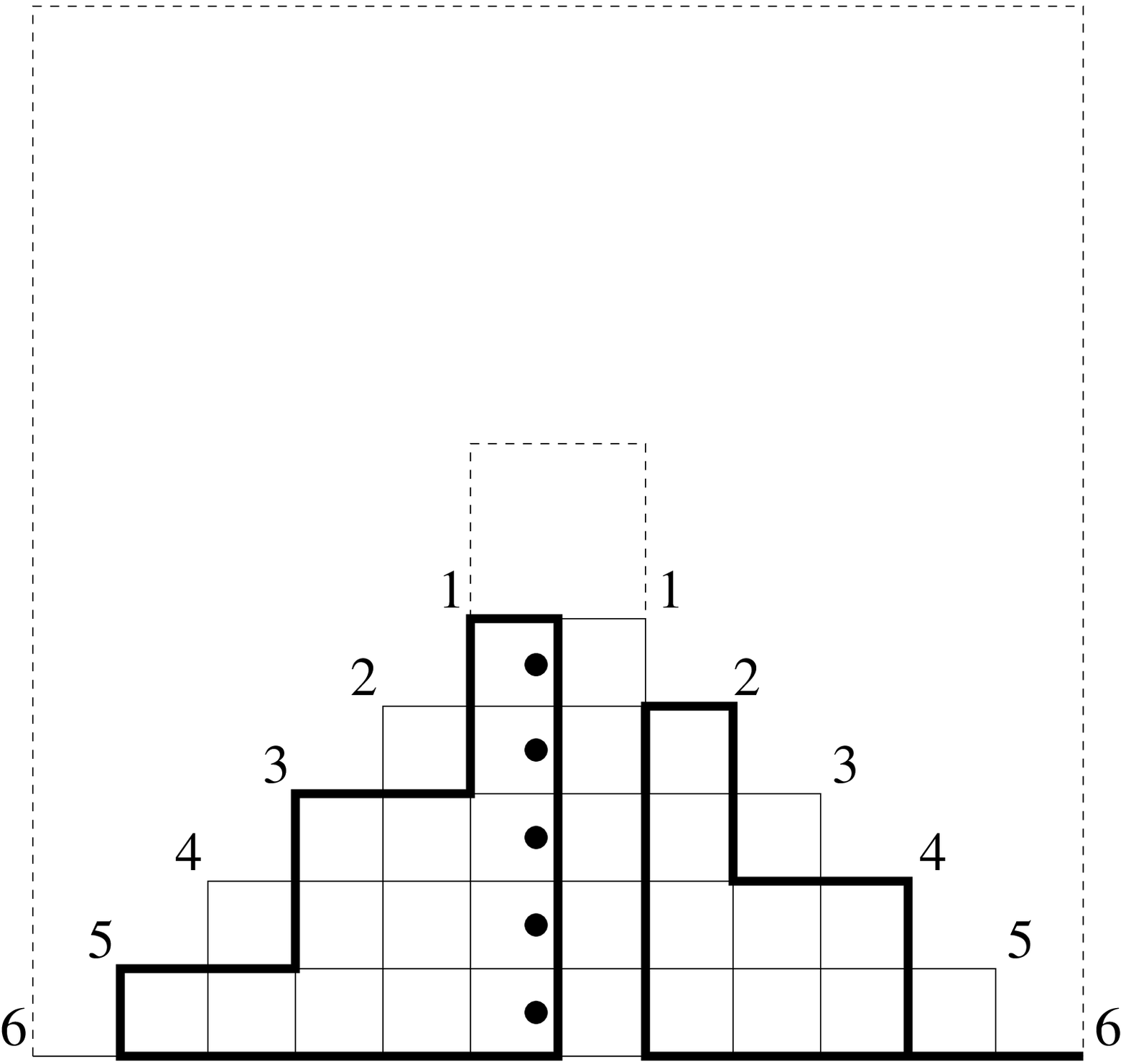}}
\medskip
\twoline{\smc Figure 5.2. {\rm Factorization for $H_7/\langle r \rangle$.}}
{\smc Figure 5.3. {\rm Factorization for $H_6/\langle r \rangle$.}}
\bigskip
\endinsert

{\it Proof of Lemma 5.1}. Apply the factorization theorem of \cite{\Cone} to the graph $C_{2n}$ 
(see Figure 5.4; the dotted edges are weighted $1/2$). 
One of the resulting subgraphs is then
$\tilde{A}_n$, while the other is $B_n$. We obtain
$$
\M(C_{2n})=2^n\M(\tilde{A}_n)\M(B_n).
$$
Plug in the expressions for $\M(B_{n})$ and $\M(C_{2n})$ given by (4.2) and (4.3). 
We obtain
$$
\align
&
\frac{1}{2^{2n}}\prod_{{\scriptstyle 1\leq j\leq k\leq 2n}\atop{\scriptstyle j+k \,\text{\rm even}}  }
\left(4-2\cos\frac{j\pi}{2n+1}-2\cos\frac{k\pi}{2n+1}\right)
\\
&\ \ \ \ \ \ \ \ \ \ \ \ \ \ \ \ \ \ \ \ \ \ \ \ \ \ \ \ \ \ \ \ \ \ \ \ \ \ \ \ 
=2^n\M(\tilde{A}_n)
\prod_{1\leq j<k\leq n}\left(4\cos^2\frac{j\pi}{2n+1}+4\cos^2\frac{k\pi}{2n+1}\right).\tag5.8
\endalign
$$
Use the identity $2\cos^2 (x) = 1+\cos(2x)$ to write
$$
\spreadlines{3\jot}
\align
4\cos^2\frac{j\pi}{2n+1}+4\cos^2\frac{k\pi}{2n+1}
&=
4+2\cos\frac{2j\pi}{2n+1}+2\cos\frac{2k\pi}{2n+1}
\\
&=4-2\cos\frac{(2n+1-2j)\pi}{2n+1}-2\cos\frac{(2n+1-2k)\pi}{2n+1}.
\endalign
$$
Rewrite the factors on the right hand side of (5.8) using the above equality, and simplify out their
occurrences in the product on the left hand side of (5.8). Solving for $\M(\tilde{A}_n)$ in the resulting
equation yields (5.1).


\topinsert
\twoline{\mypic{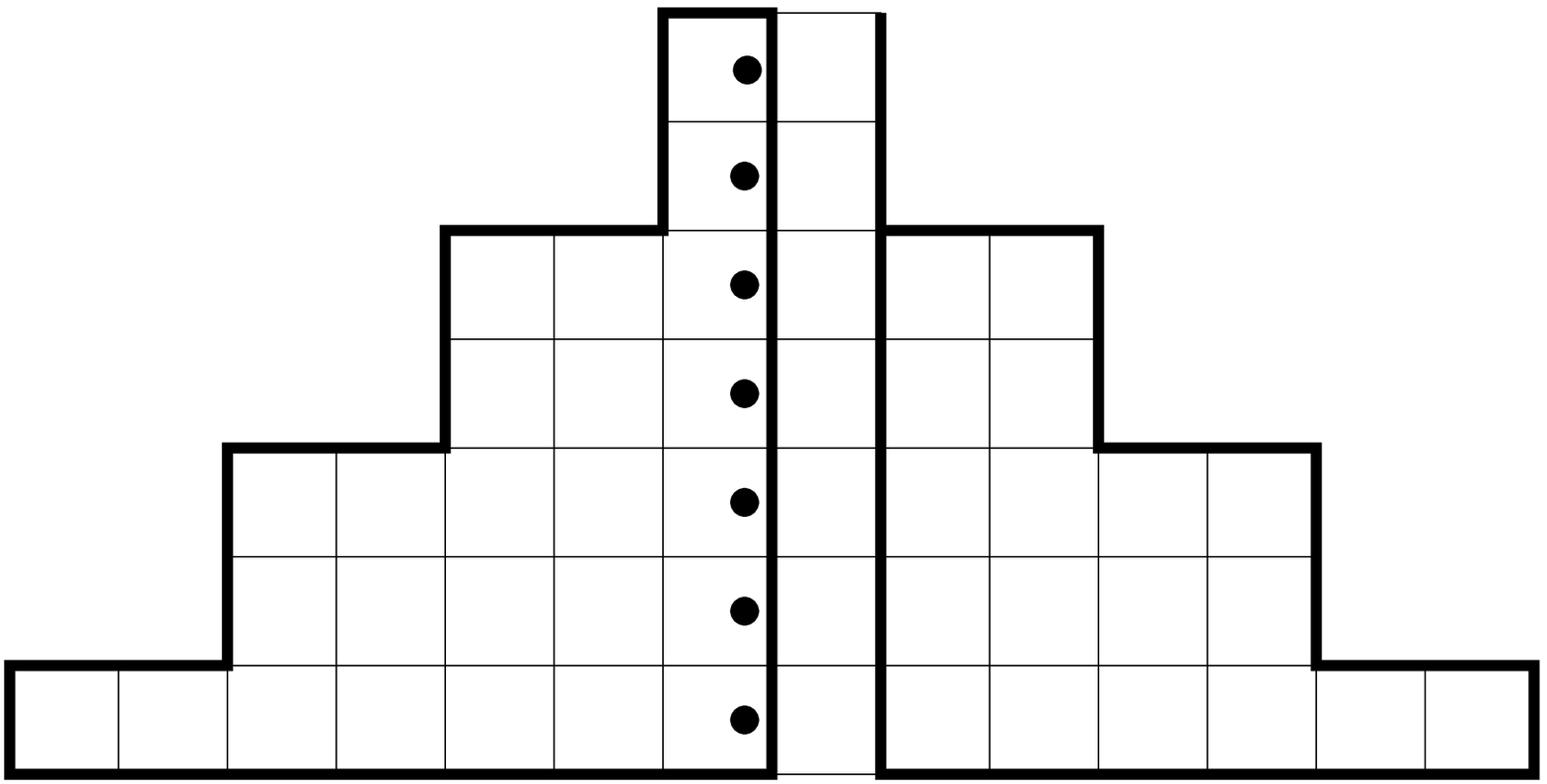}}{\mypic{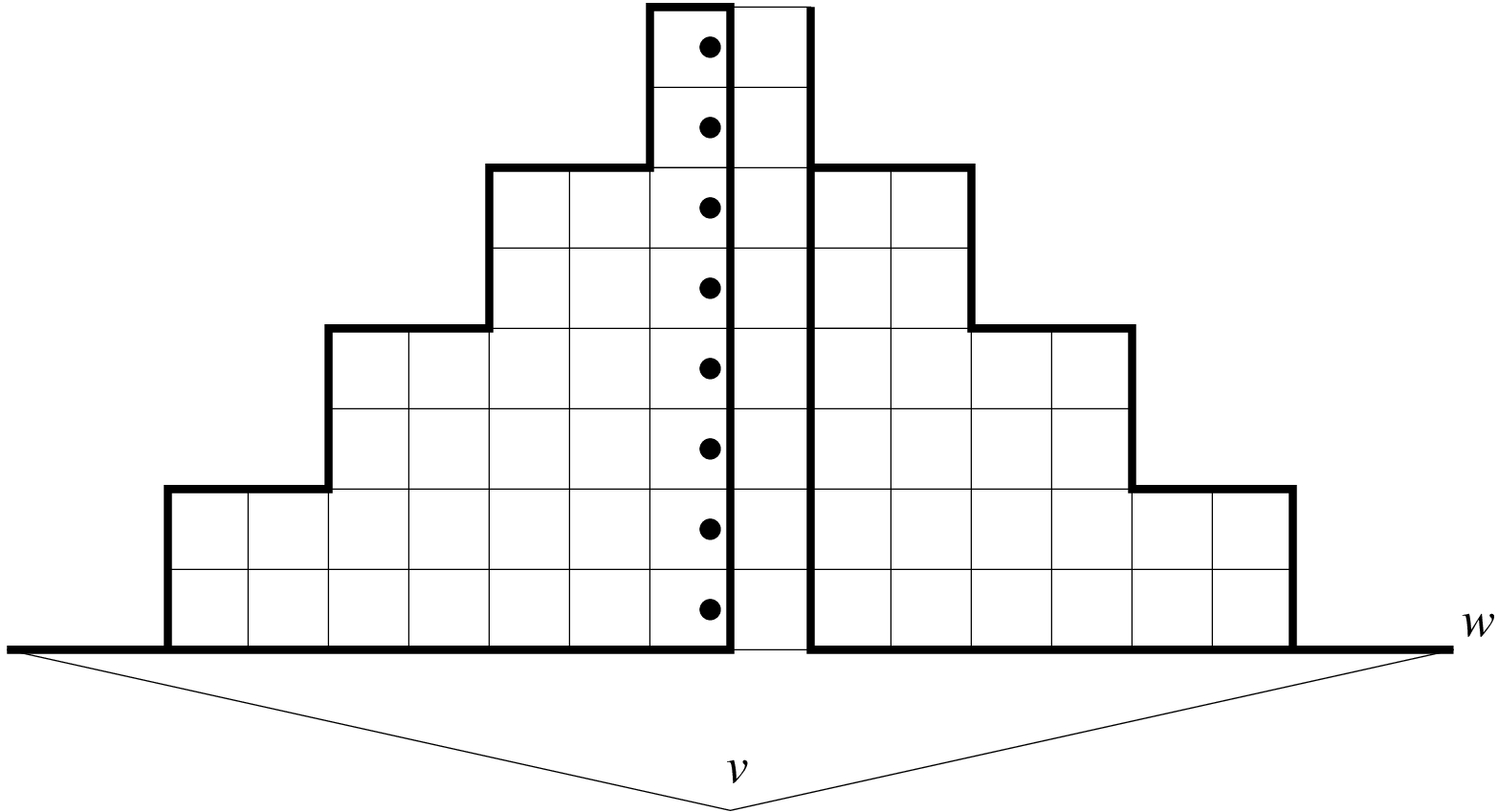}}
\twoline{\smc Figure 5.4. {\rm Factorization applied to $C_8$.}}
{\smc Figure 5.5. {\rm Factorization applied to $G_5$.}}
\endinsert

To deduce (5.2), consider the graph $G_n$ obtained from
$\Te(HMD_n)$ by including a new vertex $v$ on its vertical symmetry axis $\ell$ and 
joining $v$ by an edge to the two bottommost vertices of $HMD_n$ on its infinite face (see Figure 5.5). 
Apply the factorization theorem of \cite{\Cone} to the graph $G$ with respect to $\ell$. 
One of the resulting subgraphs is then
$\tilde{B}_n$, while the other is $\Te(QAD_n)\setminus w$, where $w$ is a leaf of $QAD_n$. We obtain
$$
\M(G_n)=2^n\M(\tilde{B}_n)\M(\Te(QAD_n)\setminus w).\tag5.9
$$
However, since $v$ can be matched in $G_n$ only to $w$ or its mirror image across $\ell$, we have by
symmetry $\M(G_n)=2\M(\Te(HMD_n)\setminus w)$. Using Theorem 4.2, (5.9) becomes
$$
2\te(HMD_n)=2^n\M(\tilde{B}_n)\te(QAD_n).
$$
Combining this with formulas (2.5) and (3.22) we obtain 
$$
\prod_{{\scriptstyle 1\leq j<k\leq 2n-1}\atop {\scriptstyle j+k\leq 2n-1}}
\left(4-4\cos\frac{j\pi}{2n}\cos\frac{k\pi}{2n}\right)
=
2^n\M(\tilde{B}_n)\prod_{1\leq j<k\leq n-1}\left(4-2\cos\frac{j\pi}{n}-2\cos\frac{k\pi}{n}\right).\tag5.10
$$
For same parity indices $j$ and $k$ in the product on the left hand side, write the corresponding factor as
$$
4-4\cos\frac{j\pi}{2n}\cos\frac{k\pi}{2n}
=4-2\cos\frac{(k-j)\pi}{2n}-2\cos\frac{(k+j)\pi}{2n}.
$$
The resulting factors cancel all the factors in the product on the right hand side of (5.10). Solving for
$\M(\tilde{B}_n)$ in the resulting equality yields (5.2).


\epf 

\mysec{6. 
Aztec pillowcases}

Theorem 2.1 shows that the square grid graph $G_n$ is similar to
the disjoint union of two copies of the same graph (namely, $QAD_{n-1}$) with a graph having a
relatively small number of vertices (namely, $P_n^{(2)}$).

Analogous statements about the Aztec diamond graphs $MD_n$ and $OD_n$ follow by Theorems~3.1 and 3.2.

These three results are instances of finite graphs $G$ with an automorphism $T$ of order 2
so that if $G_1$ is the subgraph of the orbit graph $G/T$ induced by the orbits of the vertices
not fixed by $T$, then $G$ is similar to the disjoint union of two copies of $G_1$
with some simple weighted graph on the vertices fixed by $T$. If in addition the number of vertices
fixed by $T$ is relatively small compared to the total number of vertices,
we refer to such a pair $(G,T)$ as {\it linearly squarish}.

{\it Which pairs $(G,T)$ are linearly squarish?}

\topinsert
\centerline{\mypic{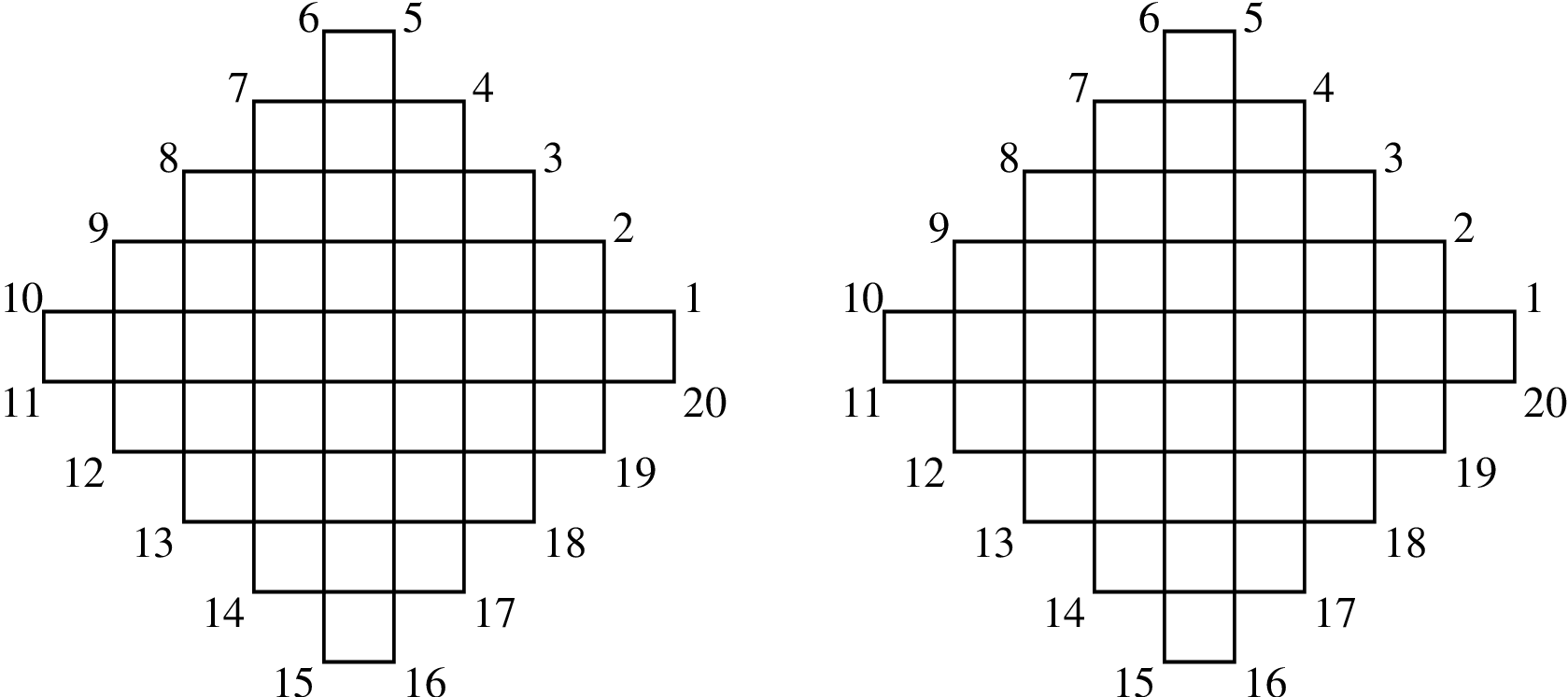}}
\medskip
\centerline{{\smc Figure~6.1.} {\rm The Aztec pillowcase graph $AP_5$.}}
\endinsert

This seems to be a difficult question to answer. Lemma 1.1 applies whenever $T^2=1$ and the fixed
points of $T$ form a cut set. Its proof shows that $G$ is similar to the disjoint union
of $G_1$ with a directed weighted graph constructed from the orbit graph $G/T$. The challenging part is to 
decide when the latter is similar
to the disjoint union of $G_1$ with some other graph. 

The question seems hard even when we restrict to the special case of graphs that are
dual to simply connected lattice regions in $\Z^2$. If one uses the natural choice (2.8) (or (3.3)) 
to define part of the desired new basis, the required action (2.9) (or (3.4)) of $F$ will hold 
as long as the corresponding vertex has its entire $\Z^2$-second neighborhood contained in $G$; but
in general problems arise near the boundary. 
We may regard the above examples of the square grid and mixed and odd Aztec diamonds as instances when we managed to
solve these problems at the boundary.


In this section we give another example of a linearly squarish graph, whose being so is based on a variation 
of the change of basis (2.8).

Define the {\it Aztec pillowcase} graph $AP_n$ to be the graph obtained from two copies of $AD_n$ by
identifying corresponding vertices on their convex hull (this leads to four pairs
of parallel edges; $AP_5$ is pictured in Figure 6.1). The {\it odd Aztec pillowcase} $OP_n$ is defined 
analogously using two copies of $OD_n$ (Figure 6.2 shows $OP_4$; the vertex identifications have been
carried out). Note that $OP_n$ is precisely the planar dual of $AP_n$. 

\topinsert
\centerline{\mypic{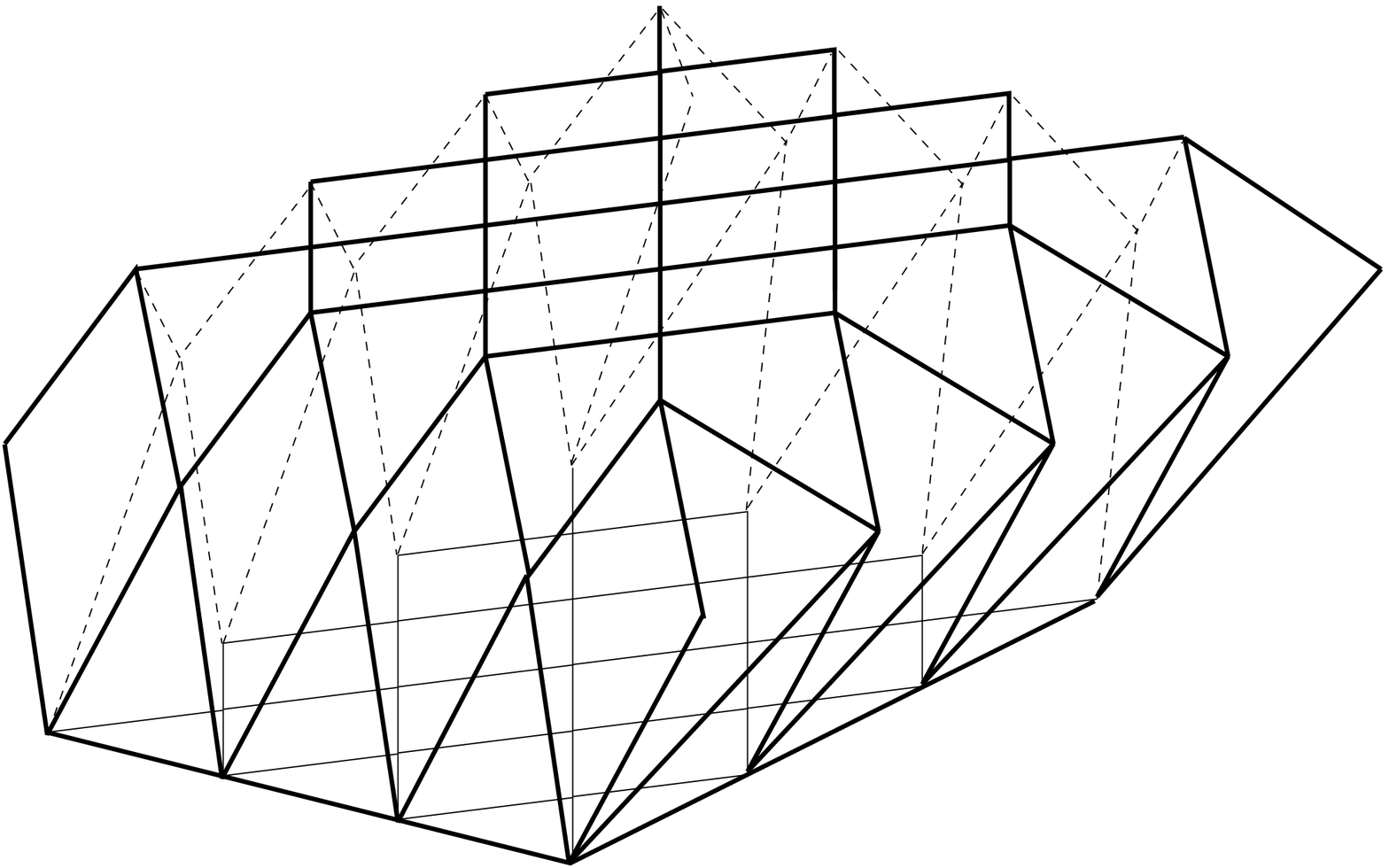}}
\medskip
\centerline{{\smc Figure~6.2.} {\rm The odd Aztec pillowcase $OP_4$.}}
\endinsert

Let $P_n^{(q)}$ be the graph obtained from the path $P_n$ by weighting each of its edges by~$q$.

Let $S_{n}$ be the directed graph obtained from $P_{n-1}^{(2)}$ by regarding each edge
as a pair of anti-parallel arcs of weight 2, and including an extra vertex $v$, an arc of weight 2 directed from $v$ 
to the rightmost vertex of $P_{n-1}^{(2)}$, and a loop at $v$ of weight 4 ($S_5$ is pictured in Figure 6.3). 
Let $S'_{n}$ be constructed analogously from $P_{n-1}^{(-2)}$, with the loop at the new vertex having weight $-4$
($S'_5$ is illustrated in Figure 6.4).

\topinsert
\twoline{\mypic{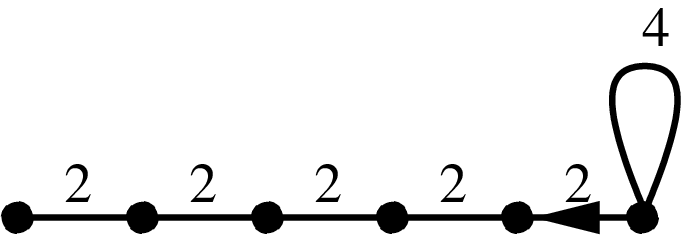}}{\mypic{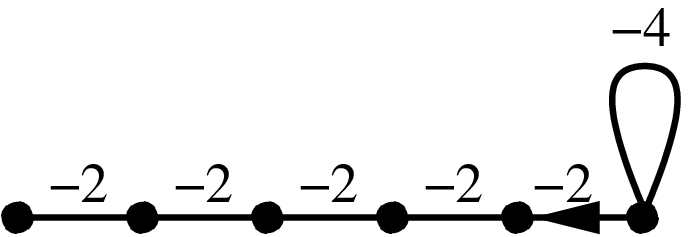}}
\twoline{Figure~6.3. {\rm $S_5$.}}{Figure~6.4. {\rm $S'_5$.}}
\endinsert

\proclaim{Theorem 6.1} $AP_n$ is similar to the disjoint union of two copies of $AD_{n-1}$,
one copy of $S_{2n-1}$, and one copy of $S'_{2n-1}$.


\endproclaim

\pf Apply Lemma 1.1 to obtain that $AP_n$ is similar to the disjoint union of $AD_{n-1}$ and
the weighted directed graph $\overline{AD}_n$ obtained from $AD_n$ by marking the vertices on its
convex hull, replacing each edge between two unmarked vertices by a pair of anti-parallel arcs of weight 1,
and replacing each edge between an unmarked vertex $v$ and a marked vertex $w$ by an arc $(v,w)$ weighted 1
and an arc $(w,v)$ weighted 2.

Define the operators $\No$, $\So$, $\Ea$ and $\We$ on the set of vertices of $\overline{AD}_n$ to be the nearest 
neighbors in specified cardinal directions, in analogy to the
operators $\NE$, $\NW$, $\SW$ and $\SE$ from the proof of Theorem 3.1. Let $V$ be the set of vertices of 
$\overline{AD}_n$, and $V_1$ the subset of vertices that are not on the boundary of its convex hull. Clearly, 
the subgraph of $\overline{AD}_n$
induced by $V_1$ of that copy is isomorphic to $AD_{n-1}$.
 
For $v\in V_1$, define
$$
f_v:=d_{\No(v)}e_{\No(v)}-d_{\We(v)}e_{\We(v)}+d_{\So(v)}e_{\So(v)}-d_{\Ea(v)}e_{\Ea(v)},\tag6.1
$$
where $d_w=1$ if $w\in V_1$, $d_w=2$ if $w\in V\setminus V_1$, and $d_\infty=0$.  

A case analysis similar to that in the proof of Lemma~2.3 checks that $\{f_v:v\in V_1\}$ spans an 
$F$-invariant subspace $\bold V$ of 
${\bold U}={\bold U}_{\overline{AD}_n}$, and that the matrix of the restriction of $F$ to $\bold V$
in the basis defined by (6.1) is just the adjacency matrix of $AD_{n-1}$. 

Let $u_1,\dotsc,u_{2n}$ be the vectors in ${\bold U}$ that are linear combinations with coefficients 0 or 1 of the basis vectors 
$\{e_v:v\in V\}$, and whose supports are given by the patterns in Figure 6.5. Let $u'_k$ be obtained from $u_k$ by multiplying
the coefficient of $e_v$ by $c_v$, where $c_v$ equals 1 for black vertices and $-1$ for whites.

\topinsert
\medskip
\centerline{\mypic{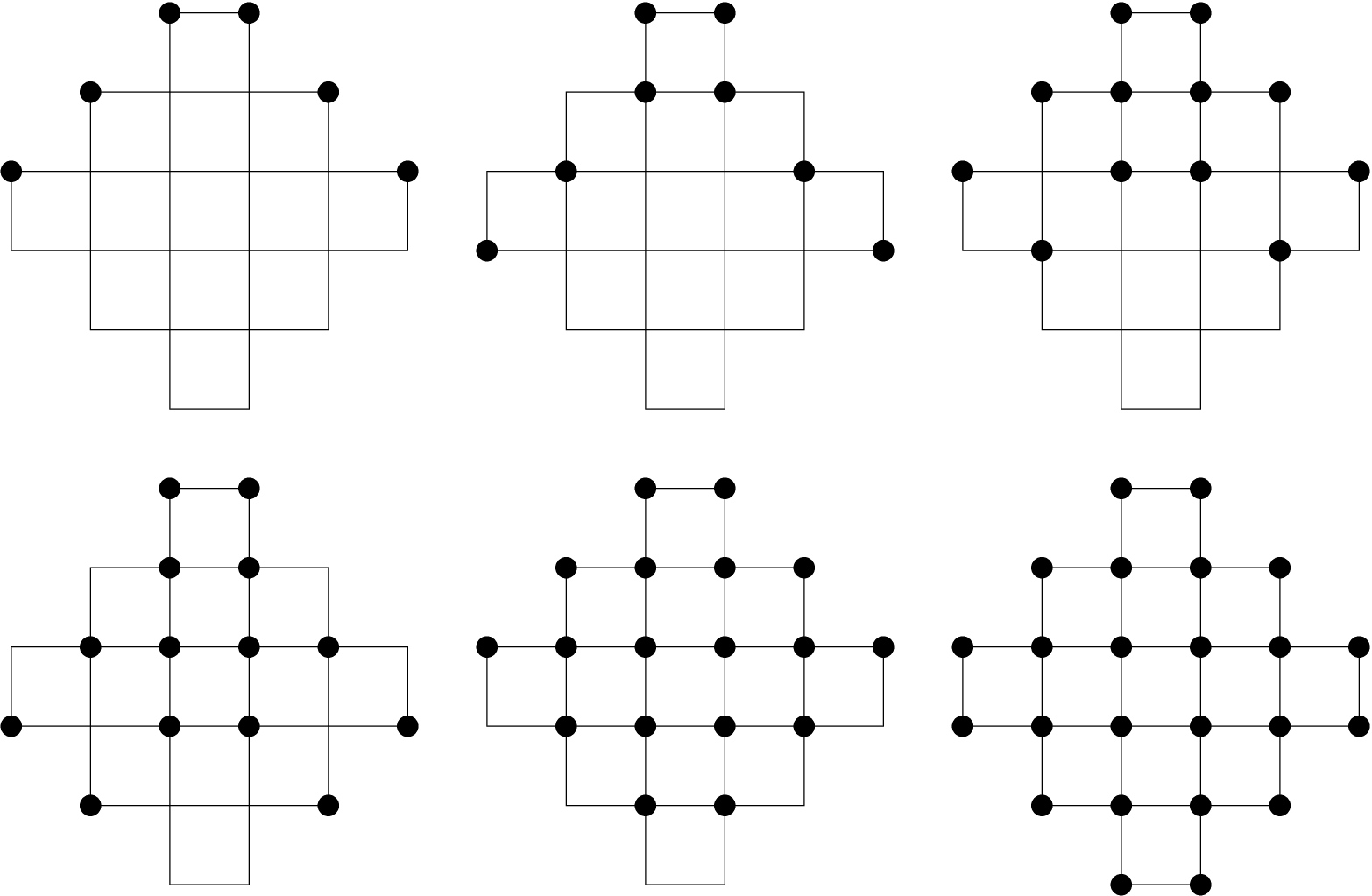}}
\centerline{\smc{Figure~6.5.} {\rm The vectors $u_k$ for $n=3$.}}
\endinsert

One readily checks that $\{u_1,\dotsc,u_{2n}\}$ are linearly independent; denote their span by ${\bold U}_1$. It is easy to verify
that ${\bold U}_1$ is $F$-invariant, and that the matrix of the restriction of $F$ to ${\bold U}_1$ in the basis 
$\{u_1,\dotsc,u_{2n}\}$ is
$$
\spreadmatrixlines{2\jot}
\left[\matrix
0&&2&&0&&0&&\cdots&&0&&0&&0\\
2&&0&&2&&0&&\cdots&&0&&0&&0\\
0&&2&&0&&2&&\cdots&&0&&0&&0\\
0&&0&&2&&0&&\cdots&&0&&0&&0\\
\vdots&&\vdots&&\vdots&&\vdots&& &&\vdots&&\vdots&&\vdots\\
0&&0&&0&&0&&\cdots&&0&&2&&0\\
0&&0&&0&&0&&\cdots&&2&&0&&0\\
0&&0&&0&&0&&\cdots&&0&&2&&4
\endmatrix\right].\tag6.2
$$
A similar statement holds for the vectors $\{u'_1,\dotsc,u'_{2n}\}$, with the restriction of $F$ to the subspace they span
being given by the negative of the above matrix. An approach similar to the one detailed for three previous examples in this paper
shows that the vectors $\{f_v:v\in V_1\}$, $\{u_1,\dotsc,u_{2n}\}$ and $\{u'_1,\dotsc,u'_{2n}\}$ span the entire space ${\bold U}$,
and thus form a basis of it. The matrix of $F$ in this basis has the desired block diagonal form. \epf


\proclaim{Corollary 6.2} The number of spanning trees of the Aztec pillowcase graphs is given by
$$
\te(AP_n)=\te(OP_n)=\frac{1}{2n^2}\prod_{j=1}^{2n-1}\left(4-4\cos\frac{j\pi}{2n}\right)^2
\prod_{1\leq j,k\leq 2n-1}\left(4-4\cos\frac{j\pi}{2n}\cos\frac{k\pi}{2n}\right).
\tag6.3
$$

\endproclaim

\pf The first equality follows from the fact that $OP_n$ is the planar dual of $AP_n$. 
The spectrum of the matrix (6.2) is readily determined. 
It follows from Theorem 6.1 and \cite{\Kn, Equations\,(5)--(6)} that
$$
\Pe(AP_n;x)=(x-4)(x+4)\prod_{j=1}^{2n-1}\left(x-4\cos\frac{j\pi}{2n}\right)^2
\prod_{1\leq j,k\leq 2n-1}\left(x-4\cos\frac{j\pi}{2n}\cos\frac{k\pi}{2n}\right).\tag6.4
$$
Since all vertex degrees in $AP_n$ are equal to 4, the eigenvalues of its negative Laplacian are
obtained by adding 4 to the eigenvalues of the adjacency matrix of $AD_n$, which in turn are
immediately read off from (6.4). Then (6.3) follows by a well known result 
(see e.g. \cite{\Kr, Theorem~VI}) which states that the number of spanning trees of a connected
graph on $m$ vertices is equal to $1/m$ times the product of the non-zero eigenvalues of its negative
Laplacian.~\epf

\medskip
\flushpar
{\smc Remark 6.3.} Embed $AP_n$ on the pillowcase manifold (topologically equivalent to $S^2$).
Let $AOP_n$ be the graph obtained by superimposing upon $AP_n$ its planar dual (which
as we pointed out is isomorphic to $OP_n$), regarding all points where edges cross or meet as vertices,
and all lines connecting such points as edges. Then Theorem 4.2 implies that 
$$
\M(AOP_n\setminus \{v,w\})=\te(AP_n)\tag6.5
$$
for any vertex $v$ of $AP_n$ and any vertex $w$ of $OP_n$ on the same face of $AOP_n$ as $v$. By (6.3) we
obtain a product formula for $\M(AOP_n\setminus \{v,w\})$. In particular, this number is independent
of $v$ and $w$, a fact not apparent a priori.

\mysec{7. Some open problems}

%

The number of perfect matchings of the $(2n+1)\times (2n+1)$ grid $H_n$ with a central unit hole 
does not seem to factor into a product of small primes. We have the following prime factorizations:
$$
\spreadlines{2\jot}
\align
\M(H_1)&=2
\\
\M(H_2)&=2^2\cdot 7^2
\\
\M(H_3)&=2^3\cdot 97^2
\\
\M(H_4)&=2^4\cdot 6121^2
\\
\M(H_5)&=2^5\cdot 31^2\cdot 113^2\cdot 271^2
\\
\M(H_6)&=2^6\cdot 592442159^2
\\
\M(H_7)&=2^7\cdot 7417^2\cdot 132605129^2 
\\
\M(H_8)&=2^8\cdot 4481^2\cdot 8513^2\cdot 9929^2\cdot 16361^2
\\
\M(H_9)&=2^9\cdot 4639^2\cdot 23357676333902111^2
\\
\M(H_{10})&=2^{10}\cdot 7^2\cdot 73^2\cdot 191^2\cdot 479^2\cdot 51151^2\cdot 2905610745223^2
\\
\M(H_{11})&=2^{11}\cdot 1033^2\cdot 1049^2\cdot 1663^2\cdot 166151^2\cdot 4241286739685449^2
\\
\M(H_{12})&=2^{12}\cdot 41^2\cdot 137^2\cdot 7057^2\cdot 20992575527970355281835400921^2.
\endalign
$$
As pointed out in \cite{\Kn}, the factors in a product formula of the kind presented in the previous 
sections can be grouped into small groups so that the product over each group is an integer. This shows that
such products have a fairly large number of prime factors. In this light, the above data suggests that
it is unlikely for $M(H_n)$ to possess a similar product expression.
%
%
It would be
interesting to see whether $\M(H_n)$ can nevertheless be expressed in terms of such products, e.g. as a 
sum of a small number of them.

We conclude with an open problem that suggests a possible generalization of the linear squarishness of
the square grid to grids in higher dimensions.

Let $G_n^{(d)}$ be the subgraph of the $d$-dimensional cubical grid $\Z^d$ induced by the lattice
points with positive coordinates less or equal to $d$ (thus $G_n^{(2)}$ is just the square grid of 
Section 2). Calculations show that for $d=3$ their characteristic polynomials start out as
$$
\spreadlines{2\jot}
\align
\Pe(G_1^{(3)};x)&=x
\\
\Pe(G_2^{(3)};x)&=(x-3)(x+3)\,[(x-1)(x+1)]^3
\\
\Pe(G_3^{(3)};x)&=x(x^2-18)\,[x^2(x^2-8)(x^2-2)^2]^3
\\
\Pe(G_4^{(3)};x)&=(x^2-3x-9)(x^2+3x-9)
\\
&\,
\times
[(x^2-3x+1)(x^2+3x+1)(x^2-x-11)(x^2+x-11)(x^2-x-1)^3(x^2+x-1)^3]^3
\\
\Pe(G_5^{(3)};x)&=x(x-3)(x+3)(x^2-27)
\\
&\,
\times
[(x-2)(x+2)(x^2-12)(x^2-4x+1)(x^2+4x+1)(x^2-2x-11)
\\
&\ \ \ \ \ 
\times
(x^2+2x-11)
(x^2-2x-2)^2(x^2+2x-2)^2(x-1)^4(x+1)^4(x^2-3)^4]^3.
\endalign
$$
This suggests that $G_n^{(3)}$ may be similar to the disjoint union of three copies of
some graph on $(n^3-n)/3$ vertices with some other graph on $n$ vertices. It would be interesting 
if one could prove this, and identify the smaller graphs.

The pattern for $d=4$ is equally suggestive. We have
$$
\spreadlines{2\jot}
\align
\Pe(G_1^{(4)};x)&=x
\\
\Pe(G_2^{(4)};x)&=x^2(x-4)(x+4)\,[x(x-2)(x+2)]^4
\\
\Pe(G_3^{(4)};x)&=x^3(x^2-32)(x^2-8)^2\,[x^4(x^2-18)(x^2-8)^2(x^2-2)^4]^4
\\
\Pe(G_4^{(4)};x)&=x^4(x-2)^2(x+2)^2(x^2-4x-16)(x^2+4x-16)(x^2-20)^2\,
\\
&\,
\times
[x^8(x-2)(x+2)(x^2-20)(x^2-4x-1)(x^2+4x-1)(x^2-2x-19)(x^2+2x-19)
\\
&\ \ \ \ \ 
\times
(x^2-2x-4)^4(x^2+2x-4)^4(x-1)^6(x+1)^6(x^2-5)^6]^4
\\
\Pe(G_5^{(4)};x)&=x^5(x-4)(x+4)(x^2-48)(x-2)^2(x+2)^2(x^2-12)^2(x^2-4x-8)^2(x^2+4x-8)^2
\\
&\,
\times
[x^{14}(x-3)(x+3)(x^2-6x+6)(x^2+6x+6)(x^2-2x-26)(x^2+2x-26)(x^2-27)
\\
&\ \ \ \ \ 
\times
(x^2-4x-8)(x^2+4x-8)(x^2-4x+1)^3(x^2+4x+1)^3
\\
&\ \ \ \ \ 
\times
(x^2-2x-11)^3(x^2+2x-11)^3
(x-2)^5(x+2)^5(x^2-12)^5
\\
&\ \ \ \ \ 
\times
(x^2-2x-2)^9(x^2+2x-2)^9(x-1)^{10}(x+1)^{10}(x^2-3)^{10}]^4.
\endalign
$$
This raises the problem of understanding whether $G_n^{(4)}$ is similar to the
disjoint union of four copies of some graph on $(n^4-n^2)/4$ vertices with some other graph on $n^2$ 
vertices. More generally, it would be interesting to see if an analogous phenomenon occurs for general $d$.

Another open problem (suggested by one of the referees) is to consider mixed Aztec pillowcase graphs and 
see if they admit a decomposition similar to the one given in Theorem 6.1.

\bigskip
{\bf Acknowledgments.} I would like to thank the two anonymous referees for their careful reading of the manuscript and
their helpful suggestions. The numbers $\M(H_n)$ shown at the beginning of Section 7 were
computed using David Wilson's program Vaxmacs for counting perfect matchings.

%
%

\bigskip

\mysec{References}
{\openup 1\jot \frenchspacing\raggedbottom
\roster

\myref{\Biggs}
  N. Biggs, ``Algebraic graph theory,'' Second edition,  
Cambridge University Press, Cambridge, 1993.  

\myref{\Cone}
  M. Ciucu, Enumeration of perfect matchings in graphs with reflective symmetry,  
{\it J. Combin. Theory Ser. A}, {\bf 77} (1997), 67--97. 

\myref{\Ctwo}
  M. Ciucu, Enumeration of spanning trees of the quartered Aztec diamond, Special Session on 
Combinatorics and Enumerative Geometry, AMS Meeting \#931, University of Louisville, Louisville, KY,
March 21, 1998 (Abstract Issue: 19/2).

\myref{\CDGT}
  D. Cvetkovi\'c, M. Doob, I. Gutman and A. Torga\v sev, ``Recent results in the theory of graph spectra,''
Annals of Discrete Mathematics 36, North-Holland Publishing Co., Amsterdam, 1988.


\myref{\KPW}
 R. W. Kenyon, J. G. Propp and D. B. Wilson, Trees and matchings,  
{\it Electron. J. Combin.}, {\bf 7} (2000), Research Paper 25, 34 pp.

\myref{\Kn}
  D. E. Knuth,
Aztec diamonds, checkerboard graphs, and spanning trees, 
{\it J. Algebraic Combin.}, {\bf 6} (1997), 253--257.

\myref{\Kr} 
  G. Kreweras, Complexit\'e et circuits eul\'eriens dans les sommes tensorielles de graphes, 
{\it J. Combin. Theory Ser. B}, {\bf  24}  (1978), 202--212.

\myref{\Lov}
  L. Lov\'asz, ``Combinatorial problems and exercises,'' Second edition,
North-Holland Publishing Co., Amsterdam, 1993. 

\myref{\Stan}
  R. P. Stanley, Unpublished manuscript, 1994.

\myref{\Temp}
  H. N. V. Temperley, ``Combinatorics,'' {\it London Math. Soc. Lecture Notes Series}, {\bf 13} (1974), 202--204.

\endroster\par}

\enddocument